\documentclass[11pt,reqno]{article}
\usepackage{amsfonts}
\usepackage{graphicx}

\begin{document}


\def\R{{\mathbb R}}
\def\T{{\mathbb T}}
\def\S{{\mathbb S}}
\def\C{{\mathbb C}}
\def\Z{{\mathbb Z}}
\def\N{{\mathbb N}}
\def\H{{\mathbb H}}
\def\B{{\mathbb B}}
\def\diam{\mbox{\rm diam}}
\def\sn{\S^{n-1}}
\def\rr{{\cal R}}
\def\mt{{\Lambda}}
\def\e{\emptyset}
\def\dQ{\partial Q}
\def\dk{\partial K}
\def\endofproof{{\rule{6pt}{6pt}}}
\def\di{\displaystyle}
\def\dist{\mbox{\rm dist}}
\def\sa+{\Sigma_A^+}
\def\du{\frac{\partial}{\partial u}}
\def\dv{\frac{\partial}{\partial v}}
\def\dt{\frac{d}{d t}}
\def\dx{\frac{\partial}{\partial x}}
\def\con{\mbox{\rm const }}
\def\nn{{\cal N}}
\def\mm{{\cal M}}
\def\kk{{\cal K}}
\def\ll{{\cal L}}
\def\vv{{\cal V}}
\def\bb{{\cal B}}
\def\mat{\mm_{at}}
\def\ma{\mm_{a}}
\def\lab{L_{ab}}
\def\labt{L_{abt}}
\def\mabn{\mm_{a}^N}
\def\man{\mm_a^N}
\def\labn{L_{ab}^N}
\def\fa{f^{(a)}}
\def\ff{{\cal F}}
\def\i{{\bf i}}
\def\gge{{\cal G}_\epsilon}
\def\gej{\chi^{(j)}_\mu}
\def\ge{\chi_\epsilon}
\def\chio{\chi^{(1)}}
\def\chit{\chi^{(2)}}
\def\chii{\chi^{(i)}}
\def\chil{\chi^{(\ell)}}
\def\gett{\chi^{(2)}_{\mu}}
\def\geol{\chi^{(1)}_{\ell}}
\def\getl{\chi^{(2)}_{\ell}}
\def\geil{\chi^{(i)}_{\ell}}
\def\gee{\chi_{\ell}}
\def\tt{{\cal T}}
\def\uu{{\cal U}}
\def\wloc{W_{\epsilon}}
\def\Int{\mbox{\rm Int}}
\def\dist{\mbox{\rm dist}}
\def\pr{\mbox{\rm pr}}
\def\pp{{\cal P}}
\def\tpp{\widetilde{\pp}}
\def\aa{{\cal A}}
\def\cc{{\cal C}}
\def\supp{\mbox{\rm supp}}
\def\Arg{\mbox{\rm Arg}}
\def\In{\mbox{\rm Int}}
\def\con{\mbox{\rm const}\;}
\def\Re{\mbox{\rm Re}}
\def\li{\mbox{\rm li}} 
\def\Seo{S^*_\epsilon(\Omega)}
\def\sdk{S^*_{\dk}(\Omega)}
\def\lae{\Lambda_{\epsilon}}
\def\ep{\epsilon}
\def\oo{{\cal O}}
\def\be{\begin{equation}}
\def\ee{\end{equation}}
\def\beqn{\begin{eqnarray*}}
\def\eeqn{\end{eqnarray*}}
\def\Pr{\mbox{\rm Pr}}

\def\gi{\gamma^{(i)}}
\def\ii{{\imath }}
\def\jj{{\jmath }}
\def\II{{\cal I}}
\def\ccij{ \cc_{i'_0,j'_0}[\eta]}
\def\dd{{\cal D}}
\def\la{\langle}
\def\ra{\rangle}
\def\bs{\bigskip}
\def\xio{\xi^{(0)}}
\def\xo{x^{(0)}}
\def\zo{z^{(0)}}
\def\Con{\mbox{\rm Const}\;}
\def\do{\partial \Omega}
\def\dk{\partial K}
\def\dl{\partial L}
\def\ll{{\cal L}}
\def\kk{{\cal K}}
\def\kk{{\cal K}}
\def\pr{{\rm pr}}
\def\ff{{\cal F}}
\def\G{{\cal G}}
\def\C{{\bf C}}
\def\dist{{\rm dist}}
\def\dds{\frac{d}{ds}}
\def\con{{\rm const}\;}
\def\Con{{\rm Const}\;}
\def\di{\displaystyle}
\def\oo{\mbox{\rm O}}
\def\hess{\mbox{\rm Hess}}
\def\gi{\gamma^{(i)}}
\def\endofproof{{\rule{6pt}{6pt}}}
\def\xm{x^{(m)}}
\def\vm{\varphi^{(m)}}
\def\km{k^{(m)}}
\def\dm{d^{(m)}}
\def\kam{\kappa^{(m)}}
\def\dem{\delta^{(m)}}
\def\xim{\xi^{(m)}}
\def\ep{\epsilon}
\def\ms{\medskip}
\def\ex{\mbox{\rm extd}}

\def\clip{C^{\mbox{\footnotesize \rm Lip}}}
\def\wlocs{W^s_{\mbox{\footnote\rm loc}}}
\def\Lip{\mbox{\rm Lip}}

\def\Xr{X^{(r)}}
\def\lip{\mbox{{\footnotesize\rm Lip}}}
\def\Vol{\mbox{\rm Vol}}

\def\naf{\nabla f(z)}
\def\so{\sigma_0}
\def\Xo{X^{(0)}}
\def\z1{z^{(1)}}
\def\Vo{V^{(0)}}
\def\Yo{Y{(0)}}

\def\uo{u^{(0)}}
\def\vo{v^{(0)}}
\def\no{\nu^{(0)}}
\def\psa{\partial^{(s)}_a}
\def\hcd{\hc^{(\delta)}}
\def\Md{M^{(\delta)}}
\def\Uo{U^{(1)}}
\def\Ut{U^{(2)}}
\def\Uj{U^{(j)}}
\def\no{n^{(1)}}
\def\nt{n^{(2)}}
\def\nj{n^{(j)}}
\def\ccm{\cc^{(m)}}

\def\ooo{\oo^{(1)}}
\def\oot{\oo^{(2)}}
\def\ooj{\oo^{(j)}}
\def\fo{f^{(1)}}
\def\ft{f^{(2)}}
\def\fj{f^{(j)}}
\def\wo{w^{(1)}}
\def\wt{w^{(2)}}
\def\wj{w^{(j)}}
\def\Vo{V^{(1)}}
\def\Vt{V^{(2)}}
\def\Vj{V^{(j)}}

\def\Ul{U^{(\ell)}}
\def\Uj{U^{(j)}}
\def\wl{w^{(\ell)}}
\def\Vl{V^{(\ell)}}
\def\Ujj{U^{(j+1)}}
\def\wjj{w^{(j+1)}}
\def\Vjj{V^{(j+1)}}
\def\Ujo{U^{(j_0)}}
\def\wjo{w^{(j_0)}}
\def\Vjo{V^{(j_0)}}
\def\vj{v^{(j)}}
\def\vl{v^{(\ell)}}

\def\f0{f^{(0)}}

\def\gl{\gamma_\ell}
\def\id{\mbox{\rm id}}
\def\piU{\pi^{(U)}}

\def\cca{C^{(a)}}
\def\bba{B^{(a)}}
\def\saa{\Sigma^+_A}
\def\sa{\Sigma_A}

\def\Int{\mbox{\rm Int}}
\def\epo{\ep^{(0)}}
\def\pH{\partial \H^{n+1}}
\def\sh{S^*(\H^{n+1})}
\def\zoo{z^{(1)}}
\def\yoo{y^{(1)}}
\def\xoo{x^{(1)}}


\def\supp{\mbox{\rm supp}}
\def\Arg{\mbox{\rm Arg}}
\def\In{\mbox{\rm Int}}
\def\diam{\mbox{\rm diam}}
\def\e{\emptyset}
\def\endofproof{{\rule{6pt}{6pt}}}
\def\di{\displaystyle}
\def\dist{\mbox{\rm dist}}
\def\con{\mbox{\rm const }}
\def\Box{\spadesuit}
\def\Int{\mbox{\rm Int}}
\def\dist{\mbox{\rm dist}}
\def\pr{\mbox{\rm pr}}
\def\be{\begin{equation}}
\def\ee{\end{equation}}
\def\beqn{\begin{eqnarray*}}
\def\eeqn{\end{eqnarray*}}
\def\la{\langle}
\def\ra{\rangle}
\def\bs{\bigskip}
\def\Con{\mbox{\rm Const}\;}
\def\clip{C^{\alpha}}
\def\wlocs{W^s_{\mbox{\footnote\rm loc}}}
\def\Lip{\mbox{\rm Lip}}
\def\lip{\mbox{\footnotesize\rm Lip}}
\def\Re{\mbox{\rm Re}}
\def\li{\mbox{\rm li}} 
\def\ep{\epsilon}
\def\ms{\medskip}
\def\dds{\frac{d}{ds}}
\def\oo{\mbox{\rm O}}
\def\hess{\mbox{\rm Hess}}
\def\id{\mbox{\rm id}}
\def\ii{{\imath }}
\def\jj{{\jmath }}
\def\graph{\mbox{\rm graph}}
\def\span{\mbox{\rm span}}
\def\Intu{\Int^u}
\def\Ints{\Int^s}

\def\i{{\bf i}}
\def\C{{\bf C}}

\def\ss{{\cal S}}
\def\tt{{\cal T}}
\def\Ee{{\cal E}}
\def\tEe{{\widetilde{\cal E}}}
\def\rr{{\cal R}}
\def\nn{{\cal N}}
\def\mm{{\cal M}}
\def\kk{{\cal K}}
\def\ll{{\cal L}}
\def\vv{{\cal V}}
\def\ff{{\cal F}}
\def\hh{{\cal H}}
\def\tt{{\cal T}}
\def\uu{{\cal U}}
\def\cc{{\cal C}}
\def\pp{{\cal P}}
\def\aa{{\cal A}}
\def\oo{{\cal O}}
\def\II{{\cal I}}
\def\dd{{\cal D}}
\def\ll{{\cal L}}
\def\ff{{\cal F}}
\def\G{{\cal G}}

\def\hs{\hat{s}}
\def\hz{\hat{z}}
\def\hL{\widehat{L}}
\def\hl{\hat{l}}
\def\hl{\hat{l}}
\def\hc{\hat{\cc}}
\def\hbb{\widehat{\cal B}}
\def\hu{\hat{u}}
\def\hX{\hat{X}}
\def\hx{\hat{x}}
\def\hu{\hat{u}}
\def\hv{\hat{v}}
\def\hQ{\hat{Q}}
\def\hC{\widehat{C}}
\def\hF{\hat{F}}
\def\hf{\hat{f}}
\def\hii{\hat{\ii}}
\def\hr{\hat{r}}
\def\hq{\hat{q}}
\def\hy{\hat{y}}
\def\hZ{\widehat{Z}}
\def\hz{\hat{z}}
\def\hE{\widehat{E}}
\def\hR{\widehat{R}}
\def\hell{\hat{\ell}}
\def\hs{\hat{s}}
\def\hW{\widehat{W}}
\def\hS{\widehat{S}}
\def\hV{\widehat{V}}
\def\hB{\widehat{B}}
\def\hhh{\widehat{\cal H}}
\def\hK{\widehat{K}}
\def\hU{\widehat{U}}
\def\hhh{\widehat{\hh}}
\def\hdd{\widehat{\dd}}
\def\hZ{\widehat{Z}}
\def\hGa{\widehat{\Gamma}}

\def\hal{\hat{\alpha}}
\def\hbe{\hat{\beta}}
\def\hg{\hat{\gamma}}
\def\hrho{\hat{\rho}}
\def\hd{\hat{\delta}}
\def\hphi{\hat{\phi}}
\def\hmu{\hat{\mu}}
\def\hnu{\hat{\nu}}
\def\hsi{\hat{\sigma}}
\def\htau{\hat{\tau}}
\def\hpi{\hat{\pi}}
\def\hep{\hat{\epsilon}}
\def\hxi{\hat{\xi}}
\def\hLa{\widehat{\Lambda}^u}
\def\hPhi{\widehat{\Phi}}
\def\hPsi{\widehat{\Psi}}
\def\hPhii{\widehat{\Phi}^{(i)}}
\def\hath{\hat{h}}

\def\tc{\tilde{C}}
\def\tg{\tilde{\gamma}}  
\def\tV{\widetilde{V}}
\def\tC{\widetilde{\cc}}
\def\tb{\tilde{b}}
\def\tt{\tilde{t}}
\def\tx{\tilde{x}}
\def\tp{\tilde{p}}
\def\tz{\tilde{Z}}
\def\tZ{\tilde{Z}}
\def\tF{\tilde{F}}
\def\tf{\tilde{f}}
\def\tp{\tilde{p}}
\def\te{\tilde{e}}
\def\tv{\tilde{v}}
\def\tu{\tilde{u}}
\def\tw{\tilde{w}}
\def\ts{\tilde{\sigma}}
\def\talpha{\tilde{\alpha}}
\def\tr{\tilde{r}}
\def\tU{\widetilde{U}}
\def\tS{\tilde{S}}
\def\tP{\widetilde{P}}
\def\ttau{\tilde{\tau}}
\def\tLip{\widetilde{\Lip}}
\def\tz{\tilde{z}}
\def\tS{\tilde{S}}
\def\tts{\tilde{\sigma}}
\def\tVl{\widetilde{V}^{(\ell)}}
\def\tVj{\widetilde{V}^{(j)}}
\def\tVo{\widetilde{V}^{(1)}}
\def\tVj{\widetilde{V}^{(j)}}
\def\tPsi{\tilde{\Psi}}
 \def\tp{\tilde{p}}
 \def\tVjo{\widetilde{V}^{(j_0)}}
\def\tvj{\tilde{v}^{(j)}}
\def\tVjj{\widetilde{V}^{(j+1)}}
\def\tvl{\tilde{v}^{(\ell)}}
\def\tVt{\widetilde{V}^{(2)}}
\def\tR{\widetilde{R}}
\def\tQ{\widetilde{Q}}
\def\oL{\tilde{\Lambda}}
\def\tq{\tilde{q}}
\def\tk{\tilde{k}}
\def\tx{\tilde{x}}
\def\ty{\tilde{y}}
\def\tz{\tilde{z}}
\def\txo{\tilde{x}^{(0)}}
\def\tso{\tilde{\sigma}_0}
\def\tmt{\tilde{\Lambda}}
\def\tg{\tilde{g}}
\def\tsi{\tilde{\sigma}}
\def\tC{\tilde{C}}
\def\tell{\tilde{\ell}}
\def\trho{\tilde{\rho}}
\def\ts{\tilde{s}}
\def\tB{\widetilde{B}}
\def\thh{\widetilde{\cal H}}
\def\tV{\widetilde{V}}
\def\trr{\tilde{r}}
\def\te{\tilde{e}}
\def\tv{\tilde{v}}
\def\tu{\tilde{u}}
\def\tw{\tilde{w}}
\def\trho{\tilde{\rho}}
\def\tell{\tilde{\ell}}
\def\tz{\tilde{Z}}
\def\tF{\tilde{F}}
\def\tf{\tilde{f}}
\def\tp{\tilde{p}}
\def\ttau{\tilde{\tau}}
\def\tz{\tilde{z}}
\def\tg{\tilde{\gamma}}  
\def\tV{\widetilde{V}}
\def\tC{\widetilde{\cc}}
\def\tLa{\widetilde{\Lambda}^u}
\def\tR{\widetilde{R}}
\def\tr{\tilde{r}}
\def\tc{\widetilde{C}}
\def\tD{\widetilde{D}}
\def\tt{\tilde{t}}
\def\tp{\tilde{p}}
\def\tS{\tilde{S}}
\def\tts{\tilde{\sigma}}
\def\tZ{\widetilde{Z}}
\def\tdelta{\tilde{\delta}}
\def\th{\tilde{h}}
\def\tB{\widetilde{B}}
\def\thh{\widetilde{\hh}}
\def\tep{\tilde{\ep}}
\def\tE{\widetilde{E}}
\def\tu{\tilde{u}}
\def\txi{\tilde{\xi}}
\def\teta{\tilde{\eta}}
\def\tRR{\widetilde{\rr}}

\def\sr{{\sc r}}
\def\mt{{\Lambda}}
\def\do{\partial \Omega}
\def\dk{\partial K}
\def\dl{\partial L}
\def\wloc{W_{\epsilon}}
\def\piU{\pi^{(U)}}
\def\Rio{\R_{i_0}}
\def\Ri{\R_{i}}
\def\Rii{\R^{(i)}}
\def\Riii{\R^{(i-1)}}
\def\hRii{\widehat{\R}_i}
\def\hRiio{\widehat{\R}_{(i_0)}}
\def\Eii{E^{(i)}}
\def\Eio{E^{(i_0)}}
\def\Rj{\R_{j}}
\def\Vio{{\cal V}^{i_0}}
\def\Vi{{\cal V}^{i}}
\def\Wio{W^{i_0}}
\def\Wioo{W^{i_0-1}}
\def\hi{h^{(i)}}
\def\Psii{\Psi^{(i)}}
\def\pii{\pi^{(i)}}
\def\piii{\pi^{(i-1)}}
\def\gxyii{g_{x,y}^{i-1}}
\def\span{\mbox{\rm span}}
\def\Jac{\mbox{\rm Jac}}
\def\Vol{\mbox{\rm Vol}}
\def\limp{\lim_{p\to\infty}}
\def\hh{{\mathcal H}}

\def\xijl{X_{i,j}^{(\ell)}}
\def\xij{X_{i,j}}
\def\hyijl{\widehat{Y}_{i,j}^{(\ell)}}
\def\hxijl{\widehat{X}_{i,j}^{(\ell)}}
\def\hxij{\widehat{X}_{i,j}}
\def\eijl{\omega_{i,j}^{(\ell)}}
\def\eij{\omega_{i,j}}
\def\Gl{\Gamma_\ell}

\def\cB{\check{B}}
\def\tpi{\tilde{\pi}}
\def\J{{\sf J}}
\def\bJ{{\mathbb J}}

\def\hcc{\widehat{\cc}}
\def\hpp{\widehat{\pp}}
\def\ttP{\widetilde{\pp}}
\def\tP{\widetilde{P}}
\def\hP{\widehat{P}}
\def\hY{\widehat{Y}}

\def\diamtef{{\footnotesize \diam_\theta}}


\def\tc{\tilde{C}}
\def\tg{\tilde{\gamma}}  
\def\tV{\widetilde{V}}
\def\tW{\widetilde{W}}
\def\tC{\widetilde{\cc}}
\def\tKo{\widetilde{K_0}}
\def\tUKo{\widetilde{U\setminus K_0}}

\def\wo{w^{(1)}}
\def\vo{v^{(1)}}
\def\uo{u^{(1)}}
\def\wt{w^{(2)}}
\def\xio{\xi^{(1)}}
\def\xit{\xi^{(2)}}
\def\etao{\eta^{(1)}}
\def\etat{\eta^{(2)}}
\def\zetao{\zeta^{(1)}}
\def\zetat{\zeta^{(2)}}
\def\vt{v^{(2)}}
\def\ut{u^{(2)}}
\def\Wo{W^{(1)}}
\def\Vo{V^{(1)}}
\def\Uo{U^{(1)}}
\def\Wt{W^{(2)}}
\def\Vt{V^{(2)}}
\def\Ut{U^{(2)}}
\def\tmu{\tilde{\mu}}
\def\tla{\tilde{\lambda}}
\def\diamf{{\rm\footnotesize diam}}
\def\Intu{\mbox{\rm Int}^u}
\def\Ints{\mbox{\rm Int}^s}

\def\Bmt{\overline{B_{\ep_0}(\mt)}}
\def\Lye{L_{y,\eta}}
\def\Lyep{L^{(p)}_{y,\eta}}
\def\Fyp{F^{(p)}_y}
\def\Fxp{F^{(p)}_x}
\def\Lxx{L_{x,\xi}}
\def\Lxxp{L^{(p)}_{x,\xi}}

\def\Wuo{W^{u,1}}
\def\Wui{W^{u,i}}
\def\Wuj{W^{u,j}}
\def\Wut{\tW^{u,2}}
\def\Wuk{W^{u,k}}
\def\Wuh{\hW^{u}}
\def\tWuo{\tW^{u,1}}
\def\tWui{\tW^{u,i}}
\def\tWuj{\tW^{u,j}}
\def\tWuk{\tW^{u,k}}
\def\hWuo{\hW^{u,1}}
\def\hWui{\hW^{u,i}}
\def\hWuj{\hW^{u,j}}
\def\hWuk{\hW^{u,k}}
\def\dj{\delta^{(j)}}
\def\do{\delta^{(1)}}
\def\epj{\ep^{(j)}}
\def\epo{\ep^{(1)}}
\def\hSj{\widehat{S}^{(j)}}
\def\hSo{\widehat{S}^{(1)}}

\def\tmu{\tilde{\mu}}
\def\tla{\tilde{\lambda}}
\def\hE{\widehat{E}}
\def\uk{u^{(k)}}
\def\ui{u^{(i)}}
\def\uj{u^{(j)}}
\def\vk{v^{(k)}}
\def\vl{v^{(l)}}
\def\vi{v^{(i)}}
\def\vj{v^{(j)}}
\def\wk{w^{(k)}}
\def\wi{w^{(i)}}
\def\wj{w^{(j)}}
\def\etak{\eta^{(k)}}
\def\etai{\eta^{(i)}}
\def\etaj{\eta^{(j)}}
\def\zetak{\zeta^{(k)}}
\def\zetai{\zeta^{(i)}}
\def\zetaj{\zeta^{(j)}}


\def\yj{y^{(j)}}
\def\yi{y^{(i)}}
\def\tyi{\ty^{(i)}}
\def\yo{y^{(1)}}
\def\zj{z^{(j)}}
\def\zo{z^{(1)}}
\def\vj{v^{(j)}}
\def\vo{v^{(1)}}
\def\kaj{\kappa^{(j)}}
\def\kao{\kappa^{(1)}}

\def\tyj{\tilde{y}^{(j)}}
\def\yl{y^{(\ell)}}
\def\tyl{\tilde{y}^{(l)}}
\def\wo{w^{(1)}}
\def\vo{v^{(1)}}
\def\vi{v^{(i)}}
\def\vj{v^{(j)}}
\def\vk{v^{(k)}}
\def\uo{u^{(1)}}
\def\wt{w^{(2)}}
\def\xio{\xi^{(1)}}
\def\xit{\xi^{(2)}}
\def\xii{\xi^{(i)}}
\def\xij{\xi^{(j)}}
\def\hxio{\hxi^{(1)}}
\def\hxit{\hxi^{(2)}}
\def\hxii{\hxi^{(i)}}
\def\hxij{\hxi^{(j)}}

\def\cxi{\check{\xi}}
\def\cxio{\cxi^{(1)}}
\def\cxit{\cxi^{(2)}}
\def\cet{\check{\eta}}
\def\ceto{\cet^{(1)}}
\def\cett{\cet^{(2)}}
\def\cv{\check{v}}
\def\cvo{\cv^{(1)}}
\def\cvt{\cv^{(2)}}
\def\cu{\check{u}}
\def\cuo{\cu^{(1)}}
\def\cut{\cu^{(2)}}
\def\cj{c^{(j)}}
\def\fj{f^{(j)}}
\def\gji{g^{(j,i)}}
\def\tPsi{\widetilde{\Psi}}
\def\chU{\check{U}}


\def\Ulo{U^{(\ell_0)}}
\def\dte{D_\theta}
\def\diamte{\mbox{\rm diam}_{\theta}}
\def\Ial{I^{(\alpha)}}
\def\uml{u_m^{(\ell)}}
\def\yl{y^{(\ell)}}
\def\tyl{\tilde{y}^{(\ell)}}
\def\ool{\oo^{(\ell)}}
\def\fl{f^{(\ell)}}
\def\hep{\hat{\ep}}
\def\dl{d^{(\ell)}}
\def\dli{d_{\ell,i}}
\def\dlo{d_{\ell,1}}
\def\dlt{d_{\ell,2}}
\def\Lipt{{\Lip_\theta}}
\def\lipt{{\footnotesize \Lip_\theta}}
\def\tm{\tilde{m}}
\def\tj{\tilde{j}}
\def\lengthf{\mbox{\rm\footnotesize length}}
\def\length{\mbox{\rm length}}


\def\Xijl{X^{(\ell)}_{i,j}}
\def\hXijl{\widehat{X}^{(\ell)}_{i,j}}
\def\Wl{W^{(\ell)}}
\def\omijl{\omega^{(\ell)}_{i,j}}

\def\hXitl{\widehat{X}^{(\ell)}_{i,t}}
\def\Vl{V^{(\ell)}}
\def\omitl{\omega^{(\ell)}_{i,t}}
\def\Xisl{X^{(\ell)}_{i,s}}
\def\hXisl{\widehat{X}^{(\ell)}_{i,s}}
\def\omisl{\omega^{(\ell)}_{i,s}}
\def\hGa{\widehat{\Gamma}}
\def\hOm{\widehat{\Omega}}
\def\tGa{\widetilde{\Gamma}}
\def\hA{\widehat{A}}
\def\tnu{\tilde{\nu}}
\def\tX{\widetilde{X}}

\def\ww{{\mathcal W}}
\def\Zl{Z^{(\ell)}}
\def\hpp{\widehat{\pp}}
\def\tnn{\widetilde{\nn}}

\def\ftt{f^{(t)}}
\def\f0{f^{(0)}}
\def\fat{f^{(at)}}
\def\Fat{F^{(at)}}
\def\Fa{F^{(a)}}
\def\F0{F^{(0)}}
\def\tu{\tilde{u}}
\def\tD{\widetilde{D}}
\def\tchi{\tilde{\chi}}
\def\tC{\widetilde{C}}
\def\hC{\widehat{C}}
\def\hQ{\widehat{Q}}
\def\hF{\widehat{F}}
\def\hD{\widehat{D}}
\def\hr{\hat{r}}
\def\psid{\psi^\dag}
\def\taud{\tau^\dag}
\def\Omn{\Omega^{(n)}}
\def\Omm{\Omega^{(m)}}
\def\Omk{\Omega^{(k)}}
\def\Conf{{\mbox{\footnotesize\rm Const}}}
\def\hp{\hat{p}}

\def\tj{t^{(j)}}
\def\tyj{\tilde{y}^{(j)}}
\def\tyjo{\tilde{y}_{j,1}}
\def\tyjt{\tilde{y}_{j,2}}
\def\tyji{\tilde{y}_{j,i}}
\def\yjo{y_{j,1}}
\def\yjt{y_{j,2}}
\def\yji{y_{j,i}}
\def\tylo{\tilde{y}_{\ell,1}}
\def\tylt{\tilde{y}_{\ell,2}}
\def\tyli{\tilde{y}_{\ell,i}}
\def\ylo{y_{\ell,1}}
\def\ylt{y_{\ell,2}}
\def\yli{y_{\ell,i}}

\def\ulo{u_{\ell,1}}
\def\ult{u_{\ell,2}}
\def\uli{u_{\ell,i}}
\def\tulo{\tilde{u}_{\ell,1}}
\def\tult{\tilde{u}_{\ell,2}}
\def\tuli{\tilde{u}_{\ell,i}}

\def\tdlo{\tilde{d}_{\ell,1}}
\def\tdlt{\tilde{d}_{\ell,2}}
\def\tdli{\tilde{d}_{\ell,i}}
\def\dlo{d_{\ell,1}}
\def\dlt{d_{\ell,2}}
\def\dli{d_{\ell,i}}
\def\wjo{w_{j,1}}
\def\wjt{w_{j,2}}
\def\wji{w_{j,i}}
\def\sj{s^{(j)}}
\def\Yj{Y^{(j)}}
\def\Vj{V^{(j)}}
\def\Zj{Z^{(j)}}
\def\vj{v^{(j)}}
\def\wj{w^{(j)}}
\def\twj{\tilde{w}^{(j)}}
\def\gj{g^{(j)}}
\def\tgj{\tilde{g}^{(j)}}
\def\tg{\tilde{g}}
\def\hn{\hat{n}}

\def\hbeta{\hat{\beta}}
\def\hmu{\hat{\mu}}
\def\piS{\pi^{(S)}}
\def\hb{\hat{b}}
\def\shP{P^\sharp}
\def\tshP{\widetilde{P}^\sharp}
\def\T{\mathcal T}
\def\tut{\tu^{(2)}}
\def\twt{\tw^{(2)}}
\def\piS{\pi^{(S)}}
\def\tsigma{\tilde{\sigma}}
\def\td{\tilde{d}}
\def\m{{\sf m}}
\def\tXi{\widetilde{\Xi}}

\def\Omb{\Omega^{(\hb)}}
\def\Xib{\Xi^{(L\hb)}}





\begin{center}
{\large\bf Ruelle transfer operators for contact Anosov flows\\ and decay of correlations}
\end{center}


\begin{center}
{\sc by Luchezar Stoyanov}
\end{center}


\footnotesize

\noindent
{\bf Abstract.} We prove exponential decay of correlations for H\"older continuous observables with respect 
to any Gibbs measure for contact Anosov flows admitting Pesin sets with exponentially small tails.  
This is achieved by establishing strong  spectral estimates for certain Ruelle transfer operators for such flows.

\normalsize

\section{Introduction and Results}
\renewcommand{\theequation}{\arabic{section}.\arabic{equation}}

\subsection{Introduction}
The study of statistical properties of dynamical systems has a long history and has been the subject 
of a considerable interest due to their applications in statistical mechanics and thermodynamics. 
Many physical systems poses some kind of `strong hyperbolicity' and are known to have or expected 
to have strong mixing properties. For example in the 70's, due to works by Sinai, Bowen and Ruelle,
it was already known that  for Anosov diffeomorphisms exponential decay of correlations takes place 
for H\"older continuous observables (see e.g. the survey article \cite{ChY}). However the continuous 
case proved to be much more difficult and it took more than twenty years until the breakthrough work of 
Dolgopyat \cite{D}, where he established exponential decay of correlations for H\"older continuous 
potentials in two major cases: (i) geodesic flows on compact surfaces of negative curvature (with respect to any 
Gibbs measure);  (ii) transitive Anosov flows  on compact Riemann manifolds with $C^1$ jointly non-integrable local 
stable and unstable foliations (with respect to the Sinai-Bowen-Ruelle measure). 

Dolgopyat's work was followed by a considerable activity to establish exponential and other types of decay of
correlations for various  kinds of systems -- see  \cite{BaL} for more information and historical remarks. 
See also \cite{Ch1}, \cite{Ch2}, \cite{ChY}, \cite{BaG}, \cite{BaT}, \cite{DL}, \cite{FT}, 
\cite{GL}, \cite{L1}, \cite{M}, \cite{N}, \cite{OWi},\cite{Y1}, \cite{Y2}, \cite{T}, \cite{Wi},
 and the references there.
Liverani \cite{L1} proved exponential decay of correlations for $C^4$ contact Anosov flows with
respect to the measure determined by the Riemann volume. Some finer results were obtained later by 
Tsujii \cite{T} (for $C^3$ contact Anosov flows) and recently by Nonnenmacher and Zworski \cite{NZ}
(for a class of $C^\infty$ flows which includes the $C^\infty$ contact Anosov flows); both papers dealing 
with the measure determined by the Riemann volume.

In this paper, as a consequnece of the main result, we derive exponential decay of correlations 
for $C^5$ contact Anosov flows on Rieman manifolds $M$ of any dimension and with respect to any 
Gibbs measure on $M$ admitting Pesin sets with exponentially small tails.  

More recently the emphasis in studying decay of correlations appears to be in trying to establish such 
results for non-uniformly hyperbolic systems and systems with singularities, e.g. billiards. In a remarkable
recent paper Baladi, Demers and Liverani \cite{BDL} established exponential decay of correlations for Sinai
billiards with finite horizon on a two-dimensional torus. See also the historical comments in \cite{BDL} for
more information on this topic.

Many of the works mentioned above used some ideas from \cite{D}, however most of them  followed a different 
approach, namely the so called functional-analytic approach initiated by the work of Blank, Keller and Liverani
\cite{BKL} which involves the study of the so called Ruelle-Perron-Frobenius operators  
$\di \ll_tg = \frac{g\circ \phi_{-t}}{|(\det d\phi_t)| \circ \phi_{-t} }$, $t\in \R$ (see e.g. 
the lectures of Liverani \cite{L2} for a nice exposition of the main ideas).

A similar approach, however studying Ruelle-Perron-Frobenius operators acting on currents, 
was used in a very recent paper by Giulietti, Liverani and Pollicott \cite{GLP} where 
they proved some remarkable results. For example, they established that
for $C^\infty$ Anosov flows the Ruelle zeta function is meromorphic in the whole complex plane. 


In \cite{D} Dolgopyat used a different approach and established some statistical properties 
(for the flows he considered) that appear to be much stronger than exponential decay of correlations. 
Indeed,  using these properties, a certain technique developed in \cite{D}  involving estimates of  
Laplace transforms of  correlations functions (following previous works of Pollicott \cite{Po} 
and Ruelle \cite{R3}),  leads more or less automatically to  exponential decay of correlations for 
H\"older continuous potentials. The approach in \cite{D} involved studying spectral properties of the so 
called Ruelle transfer operators whose definition requires a Markov partition. Given an
Anosov flow $\phi_t : M \longrightarrow M$ on a  Riemann manifold $M$, consider a Markov partition 
consisting of rectangles  $R_i = [U_i ,S_i ]$, where $U_i$ and  $S_i$ are pieces of unstable/stable 
manifolds at some $z_i\in M$, the first return time function 
$\tau : R = \cup_{i=1}^{k_0} R_i \longrightarrow [0,\infty)$ 
and the standard Poincar\'e map $\pp: R \longrightarrow R$ (see Sect. 2 for detals).
The {\it shift map} $\sigma : U = \cup_{i=1}^{k_0} U_i\longrightarrow U$, given by $\sigma = \piU\circ \pp$, where 
$\piU : R \longrightarrow U$ is the projection along the leaves of local stable manifolds, 
defines a dynamical system which is essentially isomorphic to an one-sided Markov shift.  Given a bounded function 
$f \in B (U)$, one defines the {\it Ruelle transfer operator} $L_{f} : B (U) \longrightarrow B (U)$ by 
$(L_f h)(x) = \sum_{\sigma(y) = x} e^{f(y)} h(y)$. Assuming that $f$ is real-valued and H\"older continuous, 
let $P_f \in \R$ be such that the topological pressure of $f-P_f\tau$ with respect to $\sigma$ is zero 
(cf. e.g. \cite{PP}). Dolgopyat proved (for the type of flows he considered in \cite{D}) that for small 
$|a|$ and large $|b|$ the spectral radius of the Ruelle operator 
$L_{f-(P_f+ a+\i b)\tau} : C^\alpha (U) \longrightarrow C^\alpha (U)$ 
acting on $\alpha$-H\"older continuous functions 
($0 < \alpha \leq 1$) is uniformly bounded by a constant $\rho < 1$. 

More general results of this kind were proved in \cite{St2} for mixing Axiom A flows on basic sets under some 
additional regularity assumptions,  amongst them -- Lipschitzness of the so local stable holonomy 
maps\footnote{In general these are only H\"older continuous -- see  \cite{Ha1}, \cite{Ha2}.} (see Sect. 2). 
Further results in this direction were established in \cite{St3}.

Our main result in this paper is that for contact Anosov flows on a compact  Riemann manifolds $M$ correlations 
for H\"older continuous observables decay exponentially fast with respect to any  Gibbs measure on $M$ admitting a 
Pesin set with exponentially small tails (see the definition in Sect.1.2). 

It was proved recently in \cite{GS} that Pesin sets with exponentially small tails exist for Gibbs measures for Axiom A flows (and
diffeomorphism) satisfying a certain condition, called exponential large deviations for all Lyapunov exponents
(see Sect. 3 below). In fact, under such a condition,  Pesin sets with exponentially small tails exist for any
continuous linear cocycle over a transitive subshift of finite type (see Theorem 1.7 in \cite{GS}). And it turns
out that in this generality, exponential large deviations for all exponents is a generic condition (see Theorem 1.5
in \cite{GS}). 

The main results mentioned above  are in fact consequences of a more general result.
Given $\theta \in (0,1)$, the metric $\dte$ on $U$ is defined by
$\dte(x,y) = 0$ if $x = y$, $\dte(x,y) = 1$ if $x,y$ belong to different $U_i$'s and $\dte(x,y) = \theta^N$ if 
$\pp^j(x)$ and $\pp^j(y)$ belong to the same rectangle $R_{i_j}$ for all $j = 0,1, \ldots,N-1$, 
and $N$ is the largest integer with this property. Denote by $\ff_\theta(U)$ {\it the space of all functions 
$h : U \longrightarrow \C$ with Lipschitz constants }
$|h|_\theta = \sup \{ \frac{|h(x) - h(y)|}{\dte(x,y)} : x\neq y  \, , \, x,y \in U\} < \infty .$
The central Theorem 1.3 below says that
for sufficiently large $\theta \in (0,1)$ and any real-valued function $f \in \ff_\theta(U)$ the Ruelle transfer 
operators related to $f$  are eventually  contracting on $\ff_\theta (U)$. A similar result holds for H\"older 
continuous functions on $U$ --  see Corollary 1.4 below. 

In the proof of the Theorem 1.3  we use the general framework of the method of Dolgopyat \cite{D} and 
its development in \cite{St2}, however some significant new ideas have been implemented.
The main problem is to deal with the lack of regularity of the local stable/unstable manifolds and 
related local stable/unstable holonomy maps\footnote{E.g. the local stable holonomy maps are defined by 
sliding along local stable manifolds.}  -- as we mentioned earlier, in general these are only
H\"older continuous. In \cite{D} and \cite{St2} these were assumed 
to be $C^1$ and  Lipschitz, respectively. Since the definition of Ruelle operators itself
involves sliding along local stable manifolds, it appears to be a significant problem to overcome the 
lack of regularity in general.

There are several novelties in the approach we use in this article that allow to deal with this difficulty: 
(a) making use of Pesin's theory of Lyapunov exponents; 
(b) using Liverani's Lemma B.7 in \cite{L1}\footnote{See also Appendix D in \cite{GLP} for an improved 
version of this lemma.} to estimate the so called temporal distance function\footnote{Which is only H\"older 
continuos  in general.} over cylinders using the smooth symplectic form defined by the contact form on $M$;
(c) dealing with arbitrary Gibbs measures, as long as they admit a Pesin set with exponentially small tails.
These features are of fundamental importance in this article.  Sect. 1.3 below contains more comments on the 
proof of the main result.

Here is the plan of the paper. The main results are stated in Sect. 1.2. 
Sects. 2 and 3 contain some basic definitions and facts from hyperbolic dynamics and 
Pesin's theory of Lyapunov exponents, respectively. 
The starting point in the central part of the paper is the assumption that there exists a Pesin set $P_0$ with exponentially 
small tails (see the definition in Sect. 1.2). 
In Sect. 4 we  state some properties concerning diameters of cylinders intersecting  the set $P_0$. 
It turns out that cylinders  intersecting a Pesin set $P_0$ have similar properties to these established 
in \cite{St4} under some pinching conditions. These properties (Lemma 4.1) are proved in Sect. 9.
In Sect. 4 we also state the Main Lemma 4.4 which is a rather strong non-integrability property of the contact Anosov 
flows we consider. We prove it in Sect. 8 using Liverani's Lemma (Lemma 4.2). 
Sects. 5-7, which should be regarded as the central part of this article, are devoted to 
the proofs of Theorem 1.3 and Corollary 1.4. We believe that the scope of applicability of the arguments
developed in Sects. 5-7 is significantly wider than what is actually stated as results in this paper.

\subsection{Statement of results}

Let $\phi_t : M \longrightarrow M$ be a $C^2$ contact Anosov flow on a $C^2$ compact Riemann manifold $M$.

Let $\Phi = \phi_1$ be the time-one map of the flow, and let $\m$ be an $\Phi$-invariant probability measure on $M$.
A compact subset $P$ of $M$ will be called a {\it Pesin set with exponentially small tails} with respect to $\m$
if $P$ is a Pesin set with respect to $\m$ and for every $\delta > 0$ there exist $C > 0$ and $c > 0$ such that
$$\m \left( \left\{ x\in \ll : \:\: \sharp \, \left\{ j : 0 \leq j \leq n-1 \: 
\mbox{\rm and } \:  \Phi^j(x) \notin P \right\} \geq \delta n \right\}\right) \leq C e^{- c n} ,$$
for all $n\geq 1$. See Sect. 3 for the definition of a Pesin set and for a sufficient condition for the
existence of Pesin sets with exponentially small tails. As explained below this sufficient condition is
`generic' in a certain sense.

The main result in this paper is the following.

\bs

\noindent
{\bf Theorem 1.1.} {\it Let $\phi_t : M \longrightarrow M$ be a $C^5$ contact Anosov flow,
let $F$ be a H\"older continuous function on $M$ and let $\nu_F$ be the Gibbs measure
determined by $F$ on $M$. Assume in addition that 
there exists a Pesin set with exponentially small tails with respect to $\nu_F$.
For every $\alpha > 0$ there exist constants $C = C(\alpha) > 0$ and $c = c(\alpha) > 0$ such that 
$$\left| \int_{M} A(x) B(\phi_t(x))\; d\nu_F(x) - 
\left( \int_{M} A(x)\; d\nu_F(x)\right)\left(\int_{M} B(x) \; d\nu_F(x)\right)\right|
\leq C e^{-ct} \|A\|_\alpha \; \|B\|_\alpha \;$$
for any two functions $A, B\in C^\alpha(M)$.}

\bs

We obtain this a consequence of Theorem 1.3 below  and the procedure described in \cite{D}.

It appears that so far the only results concerning exponential decay of correlations for 
general Gibbs potentials have been that of Dolgopyat \cite{D} for  geodesic flows on
compact surfaces and the one in \cite{St2} for Axiom A flows on basic sets (under additional 
assumptions including Lipschitz regularity of stable/unstable holonomy maps).
 As we mentioned earlier, Liverani \cite{L1} proved exponential decay of 
correlations for $C^4$ contact Anosov flows, and finer results (which imply exponential decay of correlations)
were  established by Tsujii \cite{T} and Nonnenmacher and Zworski \cite{NZ} (for 
$C^3$ and $C^\infty$ contact Anosov flows, respectively), however all these three papers
deal with the measure determined by the Riemann volume. 
In a recent paper Giulietti, Liverani and Pollicott \cite{GLP} derived (amongst other things) 
exponential decay of correlations for contact Anosov flows with respect to the measure of maximal 
entropy (generated by the potential $F = 0$) under a bunching condition (which implies that the 
stable/unstable foliations are $\frac{2}{3}$-H\"older).

Next, consider the {\it Ruelle zeta function} 
$$\zeta(s) = \prod_{\gamma} (1- e^{-s\ell(\gamma)})^{-1} \quad, \quad s\in \C ,$$
where $\gamma$ runs over the set of primitive  closed orbits of $\phi_t: M \longrightarrow M$
and $\ell(\gamma)$ is the least period of $\gamma$.  Denote by $h_T$ the  {\it topological entropy} 
of $\phi_t$ on $M$.

Using Theorem 1.3 below  and an argument of Pollicott and Sharp \cite{PoS1}, one derives the 
following\footnote{Instead of using the norm $\|\cdot \|_{1,b}$ as in \cite{PoS1}, in the present 
case one has to work with $\|\cdot \|_{\theta,b}$ for some  $\theta \in (0,1)$, 
and then one has to use the so called Ruelle's Lemma in the form proved in  \cite{W}. 
This is enough to prove the estimate (2.3) 
for $\zeta(s)$ in \cite{PoS1}, and from there the arguments are the same.}.

\bs

\noindent
{\bf Theorem 1.2.} {\it Let $\phi_t : M \longrightarrow M$ be a $C^2$ contact Anosov flow on a $C^2$ compact  
Riemann manifold $M$. Assume that there exists a Pesin set with exponentially small tails with respect to the
Sinai-Bowen-Ruelle measure\footnote{This is known to be true under some standard pinching conditions -- see
e.g. the comments at the end of Sect. 1 in \cite{GS}. However  we expect that this condition should be satisfied in
much more general circumstances.}. Then:} 

\ms

(a) {\it The Ruele zeta function $\zeta(s)$ of the flow  $\phi_t: M \longrightarrow M$ 
has an analytic  and non-vanishing 
continuation in a half-plane $\Re(s) > c_0$ for some $c_0 < h_T$ except for a  simple pole at $s = h_T$.  }

\ms

(b) {\it There exists $c \in (0, h_T)$ such that
$\di \pi(\lambda) = \# \{ \gamma : \ell(\gamma) \leq \lambda\} = \li(e^{h_T \lambda}) + O(e^{c\lambda})$
as $\lambda\to \infty$, where 
$\di \li(x) = \int_2^x \frac{du}{\log u} \sim \frac{x}{\log x}$ as  $x \to \infty$. }

\bigskip

Parts (a) and (b) were first established by Pollicott and Sharp \cite{PoS1} for geodesic flows on 
compact surfaces of negative curvature (using \cite{D}), and then similar results were proved in 
\cite{St2} for mixing Axiom A flows on basic sets satisfying certain additional assumptions 
(as mentioned above). Recently, using different methods, it was proved in \cite{GLP} that: 
(i) for volume preserving three dimensional Anosov flows (a) holds, and moreover, in the case of 
$C^\infty$ flows,  the Ruelle zeta function $\zeta(s)$ is meromorphic in $\C$ and 
$\zeta(s) \neq 0$ for $\Re(s) > 0$; 
(ii) (b) holds for geodesic flows on $\frac{1}{9}$-pinched compact Riemann manifolds of negative 
curvature. These were obtained as consequences of more general results in \cite{GLP}, 
one of the most remarkable being that for $C^\infty$ Anosov flows
the Ruelle zeta function $\zeta(s)$ is meromorphic in $\C$.

Let $\rr = \{R_i\}_{i=1}^{k_0}$ be a (pseudo-) Markov partition for $\phi_t$  consisting of 
rectangles $R_i = [U_i ,S_i ]$, where $U_i$ (resp. $S_i$) are (admissible) subsets of  $W^u_{\ep}(z_i)$
(resp. $W^s_{\ep}(z_i)$) for some $\ep > 0$ and $z_i\in M$ (cf. Sect. 2 for details). 
The first return time function $\tau : R = \cup_{i=1}^{k_0} R_i  \longrightarrow [0,\infty)$ 
is essentially $\alpha_1$-H\"older continous on $R$  for some $\alpha_1 > 0$,
i.e. there exists a constant $L > 0$ such that if $x,y \in R_i \cap \pp^{-1}(R_j)$  for some $i,j$, 
where $\pp: R \longrightarrow R$ is the 
standard Poincar\'e map,  then $|\tau(x) - \tau(y)| \leq L\, (d(x,y))^{\alpha_1}$.
The {\it shift map} $\sigma : U = \cup_{i=1}^{k_0} U_i \longrightarrow U$ 
is defined by $\sigma = \piU\circ \pp$, where 
$\piU : R \longrightarrow U$ is the projection along the leaves of local stable manifolds.
Let $\hU$ be the set of all $x \in U$ whose orbits do not have common points with the boundary of $R$.
Given $\theta \in (0,1)$, recall the metric $\dte$ on $\hU$ from Sect. 1.1.
Denote by $\ff_\theta(\hU)$ {\it the space of all bounded functions $h : \hU \longrightarrow \C$ 
with Lipschitz constants }
$|h|_\theta = \sup \{ \frac{|h(x) - h(y)|}{\dte(x,y)} : x\neq y; \ ; x,y \in \hU\} < \infty .$
Define the norm $\|.\|_{\theta,b}$ on $\ff_\theta (\hU)$ by  
$\| h\|_{\theta,b} = \|h\|_0 + \frac{|h|_{\theta}}{|b|}$, where $\|h\|_0 = \sup_{x\in \hU} |h(x)|$.

Given a  real-valued function $f \in \ff_\theta (\hU)$, set $g = g_f = f - P_f\tau$, where  
$P_f\in \R$ is the unique 
number such that the topological pressure $\Pr_\sigma(g)$ of $g$ with respect to $\sigma$ 
is zero (cf. \cite{PP}). 

We say  that  {\it Ruelle transfer operators related to $f$ are eventually contracting 
on $\ff_\theta(\hU)$} if for every $\epsilon > 0$ there exist constants $0 < \rho < 1$, $a_0 > 0$, 
$b_0 \geq 1$, $D_0 \geq 1$ and  $C > 0$ such that if $a,b\in \R$  satisfy $|a| \leq a_0$ and $|b| \geq b_0$, 
then\footnote{Notice that this definition is a bit different from the one in \cite{St2}.} 
$$\|L_{f -(P_f+a+ \i b)\tau}^m h \|_{\theta,b} \leq C \;\rho^m \;
\| h\|_{\theta,b}$$
for any integer $m \geq D_0 \log |b|$ and any  $h\in \ff_\theta (\hU)$.
This implies that the spectral radius  of $L_{f-(P_f+ a+\i b)\tau}$ on 
$\ff_\theta (\hU)$  does not exceed  $\rho$.

\bs

\noindent
{\bf Theorem 1.3.} {\it  Let $\phi_t : M \longrightarrow M$ be a $C^2$ contact Anosov flow 
on a $C^2$ compact  Riemann manifold $M$, let  $\rr = \{R_i\}_{i=1}^{k_0}$ be a (pseudo-) Markov 
partition for $\phi_t$ as above and let $\sigma : U \longrightarrow U$ be the corresponding shift map. 
There exists a constant  $0 <  \hat{\theta} < 1$ such that for any  $\theta \in [\hat{\theta}, 1)$ and any real-valued function 
$f\in \ff_{\theta}(\hU)$ which is the restriction of a H\"older continuous function $F$ on $M$ so that 
there exists a Pesin subset of $M$ with exponentially small tails with respect to $\nu_F$, the 
Ruelle transfer operators related to $f$ are eventually contracting on $\ff_\theta(\hU)$}.

\bs

Here $\hat{\theta}$ is the minimal number in $(0,1)$ such that the first-return time function
$\tau \in \ff_{\hat{\theta}}(\hU)$.

A similar result for H\"older continuous functions (with respect to the Riemann metric) looks a 
bit more complicated, since in general Ruelle transfer operators do not preserve any of the spaces 
$C^\alpha(\hU)$. However, they preserve a certain `filtration'
$\cup_{0 < \alpha \leq \alpha_0} \clip (\hU)$.
Here  $\alpha > 0$ and $C^{\alpha} (\hU)$ is {\it the space of all $\alpha$-H\"older 
complex-valued functions} on $\hU$. 
Then $|h|_\alpha$ is the smallest non-negative number so that 
$|h(x) - h(y)| \leq |h|_\alpha (d(x,y))^\alpha$ for all $x,y \in \hU$.
Define the norm $\|.\|_{\alpha,b}$ on $\clip (\hU)$ by  $\| h\|_{\alpha,b} = \|h\|_0 + \frac{|h|_{\alpha}}{|b|}$.

\bs

\noindent
{\bf Corollary 1.4.} {\it  Under the assumptions of Theorem 1.3,  there exists a constant $\alpha_0 > 0$ 
such that for any real-valued function  $f \in \clip(\hU)$ the  Ruelle transfer operators related to  $f$  are 
eventually contracting on $\cup_{0 < \alpha \leq \alpha_0} \clip (\hU)$. More precisely, there 
exists a constant $\hbeta\in (0,1]$
and for each $\epsilon > 0$ there exist constants $0 < \rho < 1$, $a_0 > 0$, $b_0 \geq 1$ and  
$C > 0$ such that if $a,b\in \R$  
satisfy $|a| \leq a_0$ and $|b| \geq b_0$, then for every integer $m > 0$ and every 
$\alpha \in (0,\alpha_0]$ the operator
$L_{f -(P_f+a+ \i b)\tau}^m : C^\alpha (\hU) \longrightarrow C^{\alpha \hbeta}  (\hU)$ is well-defined and
$\|L_{f -(P_f+a+ \i b)\tau}^m h \|_{\alpha\hbeta,b} \leq C \;\rho^m \;|b|^{\ep}\; \| h\|_{\alpha,b}$ for
every  $h\in C^\alpha  (\hU)$.}

\ms

The maximal constant $\alpha_0 \in (0,1]$ that one can choose above (which is determined by 
the minimal $\hat{\theta}$ one can choose in Theorem 1.3) is 
related to the regularity of the local stable/unstable foliations. Estimates for this constant 
can be derived from certain bunching 
condition concerning the rates of expansion/contraction of the  flow along local unstable/stable 
manifolds (see \cite{Ha1},  \cite{Ha2}, \cite{PSW}). In the proof of Corollary 1.4 
in Sect. 7 below we give some rough estimate for $\alpha_0$.

The above  was first proved by Dolgopyat (\cite{D}) in the case of geodesic flows on compact 
surfaces of negative curvature with $\alpha_0 = 1$ (then one can choose $\hbeta = 1$ as well).
The second main result in \cite{D} concerns transitive Anosov flows  on compact Riemann 
manifolds with $C^1$ jointly  non-integrable local stable and unstable foliations. For such 
flows Dolgopyat proved that the conclusion of Theorem 1.3 with 
$\alpha_0 = 1$ holds for the Sinai-Bowen-Ruelle  potential  $f = \log \det (d\phi_\tau)_{|E^u}$. 
More general results were proved in \cite{St2}, \cite{St4} for mixing Axiom A flows on basic sets
(again for $\alpha_0 = 1$) under some additional regularity assumptions.
For example the latter results apply to 
$C^2$ mixing Axiom A flows on basic sets satisfying a certain pinching condition (similar to the
$1/4$-pinching condition for geodesic flows on manifolds of negative curvature).








\bs

Without going into details here, let us just mention that strong spectral estimates for Ruelle 
transfer operators as the ones described  in Theorem 1.3  lead to a variety of deep results of 
various kinds -- see e.g. \cite{An},  \cite{PoS1} - \cite{PoS4}, \cite{PeS1} - \cite{PeS3} 
for some applications of the estimates in \cite{D} and \cite{St2}. Using Theorem 1.3 above,
one can prove similar results for some other relatively general systems. 

\def\ho{h^{(1)}}
\def\Ho{H^{(1)}}
\def\tcc{\widetilde{\cc}}

\bs

\footnotesize
\noindent
{\bf Acknowledgements.} Thanks are due to Boris Hasselblatt, Yakov Pesin, Mark Pollicott and Amie Wilkinson 
for providing me with various useful information, and to Vesselin Petkov for constant encouragement. 
Thanks are also due to Sebastian Gou\"ezel who discovered errors in the initial version of the paper.
Useful discussions with Mark Pollicott, Vesselin Petkov, Dima Dolgopyat, Sebastian Gou\"ezel and
Carlangelo Liverani during the program "Hyperbolic dynamics, large deviations and fluctuations" at the
Centre Bernoulli, EPFL, Lausanne (January-June 2013) are gratefully acknowledged.

\normalsize

\section{Preliminaries}
\setcounter{equation}{0}

Throughout this paper $M$ denotes a $C^2$ compact Riemann manifold,  and 
$\phi_t : M \longrightarrow M$ ($t\in \R$) a  $C^2$ Anosov flow on $M$. That is, there exist  
constants $C > 0$ and $0 < \lambda < 1$ such that there exists a $d\phi_t$-invariant decomposition  
$T_xM = E^0(x) \oplus E^u(x) \oplus E^s(x)$ of $T_xM$ ($x \in M$) 
into a direct  sum of non-zero linear subspaces,
where $E^0(x)$ is the one-dimensional subspace determined by the direction of the flow
at $x$, $\| d\phi_t(u)\| \leq C\, \lambda^t\, \|u\|$ for all  $u\in E^s(x)$ and $t\geq 0$, and
$\| d\phi_t(u)\| \leq C\, \lambda^{-t}\, \|u\|$ for all $u\in E^u(x)$ and  $t\leq 0$.

For $x\in M$ and a sufficiently small $\epsilon > 0$ let 
$$\wloc^s(x) = \{ y\in M : d (\phi_t(x),\phi_t(y)) \leq \epsilon \: \mbox{\rm for all }
\: t \geq 0 \; , \: d (\phi_t(x),\phi_t(y)) \to_{t\to \infty} 0\: \}\; ,$$
$$\wloc^u(x) = \{ y\in M : d (\phi_t(x),\phi_t(y)) \leq \epsilon \: \mbox{\rm for all }
\: t \leq 0 \; , \: d (\phi_t(x),\phi_t(y)) \to_{t\to -\infty} 0\: \}$$
be the (strong) {\it stable} and {\it unstable manifolds} of size $\epsilon$. Then
$E^u(x) = T_x \wloc^u(x)$ and $E^s(x) = T_x \wloc^s(x)$. 
Given $\delta > 0$, set $E^u(x;\delta) = \{ u\in E^u(x) : \|u\| \leq \delta\}$;
$E^s(x;\delta)$ is defined similarly. 

It follows from the hyperbolicity of the flow on $M$  that if  $\epsilon_0 > 0$ is sufficiently small,
there exists $\ep_1 > 0$ such that if $x,y\in M$ and $d (x,y) < \ep_1$, 
then $W^s_{\ep_0}(x)$ and $\phi_{[-\ep_0,\ep_0]}(W^u_{\ep_0}(y))$ intersect at exactly 
one point $[x,y]$  (cf. \cite{KH}). That is, there exists a unique  $t\in [-\ep_0, \ep_0]$ such that
$\phi_t([x,y]) \in W^u_{\ep_0}(y)$. Setting $\Delta(x,y) = t$, defines the so called {\it temporal distance
function\footnote{In fact in \cite{D} and \cite{L1} a different definition for $\Delta$ is given, however
in the important case (the only one considered below) when $x\in W^u_\ep(z)$ and 
$y \in W^s_\ep(z)$ for some $z \in M$, 
these definitions coincide with  the present one.}} (\cite{KB},\cite{D}, \cite{Ch1}, \cite{L1}).
For $x, y\in M$ with $d (x,y) < \ep_1$, define
$\pi_y(x) = [x,y] = W^s_{\ep}(x) \cap \phi_{[-\ep_0,\ep_0]} (W^u_{\ep_0}(y))$.
Thus, for a fixed $y \in M$, $\pi_y : W \longrightarrow \phi_{[-\ep_0,\ep_0]} (W^u_{\ep_0}(y))$ is the
{\it projection} along local stable manifolds defined on a small open neighbourhood $W$ of $y$ in $M$.
Choosing $\ep_1 \in (0,\ep_0)$ sufficiently small,  the restriction 
$\pi_y: \phi_{[-\ep_1,\ep_1]} (W^u_{\ep_1}(x)) \longrightarrow \phi_{[-\ep_0,\ep_0]} (W^u_{\ep_0}(y))$
is called a {\it local stable holonomy map\footnote{In a similar way one can define
holonomy maps between any two sufficiently close local transversals to stable laminations; see e.g.
\cite{PSW}.}.} Combining such a map with a shift along the flow we get another local stable holonomy  map
$\hh_x^y : W^u_{\ep_1}(x)  \longrightarrow W^u_{\ep_0}(y)$.
In a similar way one defines local holonomy maps along unstable laminations.

We will say that $A$ is an {\it admissible subset} of $W^u_{\ep}(z)$
if $A$ coincides with the closure of its interior in $W^u_\ep(z)$. Admissible subsets of 
$W^s_\ep(z)$ are defined similarly.

Let $D$ be a submanifold of $M$ of codimension one such that $\diam(D) \leq \ep$ and
$D$ is transversal to the flow $\phi_t$. Assuming that $\ep > 0$ is sufficiently small,
the projection $\pr_D : \phi_{[-\ep,\ep]}(D) \longrightarrow D$ along the flow is well-defined
and smooth. Given $x,y\in D$, set $\la x, y\ra_D = \pr_D([x,y])$. 
A subset $\tR$ of $D$ is called a {\it rectangle} if $\la x, y\ra_D \in \tR$ for all
$x,y\in \tR$.  The rectangle $\tR$ is called {\it proper} if $\tR$ coincides with the closure of its
interior in $D$. For any $x\in \tR$ define the stable and unstable leaves through $x$ in $\tR$
by $W^s_{\tR}(x) = \pr_D(W^s_\ep(x)\cap \phi_{[-\ep,\ep]}(D)) \cap \tR$ and
$W^u_{\tR}(x) = \pr_D(W^u_\ep(x)\cap \phi_{[-\ep,\ep]}(D)) \cap \tR$. For a subset $A$ of $D$ we will denote
by $\Int_D(A)$ the {\it interior} of $A$ in $D$. 

Let $\tRR = \{ \tR_i\}_{i=1}^{k_0}$ be a family of proper rectangles, where each $\tR_i$ is contained
in a submanifold $D_i$ of $M$ of codimension one. We may assume that each $\tR_i$ has the form 
$$\tR_i = \la U_i  , S_i \ra_{D_i} = \{ \la x,y\ra_{D_i} : x\in U_i, y\in S_i\}\;,$$
where $U_i \subset \wloc^u(z_i)$ and $S_i \subset \wloc^s(z_i)$, respectively, 
for some $z_i\in M$.  Moreover, we can take $D_i$ so that $U_i \cup S_i \subset D_i$. 
Set $\tR =  \cup_{i=1}^{k_0} \tR_i$. We will denote by $\Int(\tR_i)$ the {\it interior} of the set $\tR_i$ 
in the topology of the disk $D_i$.
The family $\tRR$ is called {\it complete} if  there exists $\chi > 0$ such that for every $x \in M$,
$\phi_{t}(x) \in \tR$ for some  $t \in (0,\chi]$.  The {\it Poincar\'e map} $\tpp: \tR \longrightarrow \tR$
related to a complete family $\tRR$ is defined by $\tpp(x) = \phi_{\ttau(x)}(x) \in \tR$, where
$\ttau(x) > 0$ is the smallest positive time with $\phi_{\ttau(x)}(x) \in \tR$.
The function $\ttau$  is called the {\it first return time}  associated with $\trr$. 
A complete family $\tRR = \{ \tR_i\}_{i=1}^{k_0}$ of rectangles in $M$ is called a 
{\it Markov family} of size $\chi > 0$ for the  flow $\phi_t$ if: (a) $\diam(\tR_i) < \chi$ for all $i$; 
(b)  for any $i\neq j$ and any $x\in \Int_{D}(\tR_i) \cap \tpp^{-1}(\Int_{D}(\tR_j))$ we have   
$W_{\tR_i}^s(x) \subset \overline{ \tpp^{-1}(W_{\tR_j}^s(\tpp(x)))}$ and 
$\overline{\tpp(W_{\tR_i}^u(x))} \supset W_{\tR_j}^u(\tpp(x))$;
(c) for any $i\neq j$ at least one of the sets $\tR_i \cap \phi_{[0,\chi]}(\tR_j)$ and
$\tR_j \cap \phi_{[0,\chi]}(\tR_i)$ is empty. 

The existence of a Markov family $\tRR$ of an arbitrarily small size $\chi > 0$ for $\phi_t$
follows from the construction of Bowen \cite{B}.

Following  \cite{R2} and \cite{D}, we will now slightly change the Markov family $\tRR$ to a 
{\it pseudo-Markov partition} $\rr = \{ R_i\}_{i=1}^{k_0}$ of {\it pseudo-rectangles} 
$R_i = [U_i  , S_i ] =  \{ [x,y] : x\in U_i, y\in S_i\}\,.$
where $U_i$ and $S_i$ are as above. Set $R =  \cup_{i=1}^{k_0} R_i$. Notice that 
$\pr_{D_i} (R_i) = \tR_i$ for all $i$.
Given $\xi = [x,y] \in R_i$, set $W^u_{R}(\xi) = W^u_{R_i}(\xi) = [U,y] = \{ [x',y] : x'\in U_i\}$ and
$W^s_{R}(\xi) = W^s_{R_i}(\xi) = [x,S_i] = \{[x,y'] : y'\in S_i\} \subset W^s_{\ep_0}(x)$.
The corresponding {\it Poincar\'e map} $\pp: R \longrightarrow R$ is defined by  
$\pp(x) = \phi_{\tau(x)}(x) \in R$, where
$\tau(x) > 0$ is the smallest positive time with $\phi_{\tau(x)}(x) \in R$. 
The function $\tau$  is the {\it first return time}  associated with $\rr$. 
The {\it interior} $\Int(R_i)$ of a rectangle $R_i$ is defined by 
$\pr_D(\Int(R_i)) = \Int_D(\tR_i)$. In a similar way one can define
$\Intu(A)$ for a subset $A$ of some  $W^u_{R_i}(x)$ and $\Ints(A)$ for a subset $A$ of some  $W^s_{R_i}(x)$.

We may and will assume that the family $\rr = \{ R_i\}_{i=1}^{k_0}$ has the same properties as $\tRR$,
namely: $(a')$ $\diam(R_i) < \chi$ for all $i$; $(b')$ for any $i\neq j$ and any  
$x\in \Int(R_i) \cap \pp^{-1}(\Int(R_j))$ we have   
$\pp(\Int (W_{R_i}^s(x)) ) \subset \Ints (W_{R_j}^s(\pp(x)))$ and 
$\pp(\Int(W_{R_i}^u(x))) \supset \Int(W_{R_j}^u(\pp(x)))$; $(c')$
for any $i\neq j$ at least one of the sets $R_i \cap \phi_{[0,\chi]}(R_j)$ and 
$R_j \cap \phi_{[0,\chi]}(R_i)$ is empty.
Define the matrix $A = (A_{ij})_{i,j=1}^k$  by $A_{ij} = 1$ if $\pp(\Int (R_i)) \cap \Int(R_j) \neq  \e$ 
and $A_{ij} = 0$ otherwise.  According to \cite{BR} (see section 2 there), we may assume that $\rr$ 
is chosen in such a way that $A^{M_0} > 0$ (all entries of the $M_0$-fold product of $A$  by itself are 
positive) for some integer $M_0 > 0$. In what follows we assume that the matrix $A$ has this property.

Notice that in general $\pp$ and $\tau$ are only (essentially) H\"older continuous. However there is an
obvious relationship between $\pp$ and the (essentially) Lipschitz map $\tpp$, and this will be used below.

From now on we will assume that $\tRR = \{ \tR_i\}_{i=1}^{k_0}$ is a fixed Markov family for  
$\phi_t$ of size $\chi < \ep_0/2 < 1$
and that  $\rr = \{ R_i\}_{i=1}^{k_0}$ is the related pseudo-Markov family. Set  
$$U = \cup_{i=1}^{k_0} U_i $$
and $\Intu (U) = \cup_{j=1}^{k_0} \Intu(U_j)$.

The {\it shift map} $\sigma : U   \longrightarrow U$ is given by
$\sigma  = \piU \circ \pp$, where $\piU : R \longrightarrow U$ is the {\it projection} along stable leaves. 
Notice that  $\tau$ is constant on each stable leaf $W_{R_i}^s(x) = W^s_{\ep_0}(x) \cap R_i$. 
For any integer $m \geq 1$
and any function $h : U \longrightarrow \C$ define $h_m : U \longrightarrow \C$ by
$$h_m(u) = h(u) + h(\sigma(u)) + \ldots + h(\sigma^{m-1}(u)) .$$

Denote by $\widehat{U}$ (or $\hR$) the {\it core} of  $U$ (resp. $R$), i.e. the set 
of those $x\in U$ (resp. $x \in R$)
such that  $\pp^m(x) \in \Int(R) = \cup_{i=1}^k \Int(R_i)$  for all $m \in \Z$. 
It is well-known (see \cite{B}) that $\hU$ is a residual subset 
of $U$ (resp. $R$) and has full measure with respect to any Gibbs measure on $U$ (resp. $R$).
Clearly in general $\tau$ is not continuous on $U$, however $\tau$ is {\it essentially H\"older} on $\hU$. 
The same applies to $\sigma : U \longrightarrow U$.  Throughout we will mainly 
work with the restrictions of $\tau$ and $\sigma$ to $\hU$. Set $\hU_i = U_i \cap \hU$.
For any $A \subset M$, let $\hA $  be {\it the set of all $x \in A$ whose trajectories 
do not pass through boundary points} of $R$.

Let $B(\hU)$ be the {\it space of  bounded functions} $g : \hU \longrightarrow \C$ with its standard norm  
$\|g\|_0 = \sup_{x\in \hU} |g(x)|$. Given a function $g \in B(\hU)$, the  {\it Ruelle transfer operator } 
$L_g : B(\hU) \longrightarrow B(\hU)$ is defined by 
$$\di (L_gh)(u) = \sum_{\sigma(v) = u} e^{g(v)} h(v) .$$
Given $\alpha > 0$, let $\clip(\hU)$ denote the {\it space of essentially $\alpha$-H\"older continuous functions}
$h :\hU \longrightarrow \C$, i.e. such that there exists $L \geq 0$ with 
$|h(x) -  h(y) | \leq L \, (d(x,y))^\alpha$
for all $i = 1, \ldots, k_0$ and all $x,y\in \hU_i$. The smallest $L > 0$ 
with this property is called the $\alpha$-H\"older
exponent of $h$ and is denoted $|h|_\alpha$.  
Set $\|g\|_\alpha = \|g\|_0 + |g|_\alpha$.

The hyperbolicity of the flow  implies the existence of
constants $c_0 \in (0,1]$ and $\gamma_1 > \gamma > 1$ such that
\be
c_0 \gamma^m\; d (x,y) \leq  d (\tpp^m(x)), \tpp^m(y)) \leq \frac{\gamma_1^m}{c_0} d (x,y)
\ee
for all $x,y\in \tR$ such that $\tpp^j(x), \tpp^j(y)$ belong to the same $\tR_{i_j}$ for all $j = 0,1, \ldots, m$. 

Throughout this paper $\alpha_1 \in (0,1]$ will denote the largest constant 
such that  $\tau\in C^{\alpha_1}(\hU)$ and the local stable/unstable
holonomy maps are uniformly $\alpha_1$-H\"older. We will also need to fix a constant $\talpha_1 \in (0,1)$
(take e.g. the largest again)  such that the projection $\tPsi : R \longrightarrow \tR$ along stable
leaves is $\talpha_1$-H\"older.

\section{Lyapunov exponents and Lyapunov regularity functions}
\setcounter{equation}{0}


Let $M$ be a $C^2$ Riemann manifold, and let $\phi_t$ be a $C^2$ Anosov  flow on $M$.
Let $\Phi$ be a H\"older continuous real-valued function on $M$ and let $\m$ be the Gibbs measure 
generated by $\Phi$ on $M$. Then $\m(\ll') = 1$, where $\ll'$ is the set of all {\it Lyapunov regular 
points} of $f = \phi_1$ (see \cite{P1} or section 2.1 in  \cite{BP}).  
There exists a subset $\ll$ of $\ll'$ with $\m(\ll) = 1$ such that the positive Lyapunov exponents
$\chi_1  < \chi_2 < \ldots < \chi_{\tk}$
of $f$ are constant on $\ll$. For $x\in \ll$, let 
$E^u(x) = E^u_1(x) \oplus E^u_2(x) \oplus \ldots \oplus E^u_{\tk}(x)$
be the $d\phi_t$-invariant decomposition of $E^u(x)$ into subspaces of constant dimensions 
$n_1, \ldots, n_{\tk}$ with $n_1 + n_2 + \ldots + n_{\tk} = n^u = \dim(E^u(x))$.
We have a similar decomposition for $E^s(x)$, $x\in \ll$. If the flow is contact, we have 
$n^s = \dim(E^s(x))  = n^u$.

Set  $\lambda_i = e^{\chi_i}$ for all $ i = 1, \ldots,\tk$. {\bf Fix an arbitrary constant} $\beta \in (0,1]$ such that 
$\lambda_j^\beta < \lambda_{j+1} $ for all $\leq j < \tk .$
Take $\hep > 0$ so small that 
\begin{eqnarray}
e^{8\hep} < \lambda_1  \:\: , \:\: e^{8\hep} < \lambda_{j}/\lambda_{j-1} \:\: (j = 2, \ldots,\tk) .
\end{eqnarray}
Some further assumptions about $\hep$ will be made later.
Set
\be
1 < \nu_0 = \lambda_1 e^{-8\hep} <  \mu_{j} = \lambda_{j} e^{-\hep} < \lambda_{j} < \nu_{j} = \lambda_{j} e^{\hep}
\ee
for all $j = 1, \ldots,\tk$.

{\bf Fix $\hep > 0$ with the above properties and set} $\ep = \hep/4$. There exists a 
{\it Lyapunov $\epsilon$-regularity function}  $R = R_{\ep} : \ll \longrightarrow (1,\infty)$, i.e. a function with
\be
e^{-\ep} \leq \frac{R (f(x))}{R (x)} \leq  e^{\ep} \quad , \quad x\in \ll ,
\ee
such that
\be
\frac{1}{R (x)\, e^{n\ep}} \leq \frac{\|df^n(x)\cdot v\|}{\lambda_i^n\|v\|} 
\leq R (x)\, e^{n\ep} \quad , \quad x\in \ll \;, \; v\in E^u_i(x)\setminus \{0\} \;, \; n \geq 0 .
\ee
We will discuss these functions in more details in Sect. 3.2.

\def\utk{u^{(\tk)}}

For $x \in \ll$ and $1 \leq j \leq d$ set 
$$\hE^u_j(x) = E^u_1(x) \oplus \ldots \oplus E^u_{j-1}(x) \quad , \quad \tE^u_{j} 
= E^u_j(x) \oplus \ldots \oplus E^u_{\tk}(x) .$$
Also set $\hE^u_1(x) = \{0\}$ and $\hE^u_{\tk+1}(x) = E^u(x)$.
For any $x\in \ll$ and any $u\in E^u(x)$ we will write $u = (\uo,\ut , \ldots, \utk)$, 
where $\ui \in E^u_i(x)$ for all $i$.
We will denote by $\|\cdot\|$ the norm on $E^u(x)$ generated by the Riemann metric.

It follows from the general theory of non-uniform hyperbolicity (see \cite{P1}, \cite{BP}) 
that for any $j = 1, \ldots,\tk$ the invariant bundle  $\{\tE^u_{j}(x)\}_{x\in \ll}$ is uniquely integrable over $\ll$, i.e. 
there exists a continuous $f$-invariant family $\{ \tWuj_{\tr(x)}(x)\}_{x\in \ll}$  of $C^2$ submanifolds
$\tWuj(x) = \tWuj_{\tr(x)}(x)$ of $M$ tangent to the bundle $\tE^u_{j}$ for some 
Lyapunov {\it $\hep/2$-regularity function} $\tr = \tr_{\hep/2} : \ll \longrightarrow (0,1)$.
Moreover, with $\beta \in (0,1]$ as in the beginning of this section, for $j > 1$ it follows from Theorem 6.6 in 
\cite{PS} and (3.1) that there exists an $f$-invariant family $\{ \hWuj_{\tr(x)}(x)\}_{x\in \ll}$ of
$C^{1+\beta}$ submanifolds $\hWuj(x) = \hWuj_{\tr(x)}(x)$ of $M$ tangent to the bundle $\hE^u_j$.
(However this family is not unique in general.) For each $x\in \ll$ and each $j = 2, \ldots, \tk$ fix
an $f$-invariant family $\{ \hWuj_{\tr(x)}(x)\}_{x\in \ll}$ with the latter properties.
Then we can find a Lyapunov $\hep$-regularity function  $r = r_{\hep}: \ll \longrightarrow (0,1)$ 
and for any $x\in \ll$ a $C^{1+\beta}$ diffeomorphism
$$\Phi_x^u: E^u (x; r(x)) \longrightarrow \Phi_x (E^u (x; r (x)) \subset W^u_{\tr (x)}(x)$$
such that 
\be
\Phi^u_x(\hE^u_j(x; r(x))) \subset \hWuj_{\tr(x)}(x)\quad , \quad
\Phi^u_x(\tE^u_{j}(x; r(x))) \subset \tWuj_{\tr (x)}(x) 
\ee
for all $x \in \ll$ and $j = 2, \ldots , \tk$. Moreover, since for each $j > 1$ the submanifolds
$ \hWuj_{r(x)}(x)$ and $\exp^u_x(\hE^u_j(x;r(x)))$ of $W^u_{\tr(x)}(x)$ are tangent at $x$ of order 
$1+\beta$, we can choose  $\Phi^u_x$ so that the diffeomorphism
$$\Psi^u_x =  (\exp^u_x)^{-1} \circ \Phi^u_x   : E^u(x: r(x)) \longrightarrow \Psi^u_x (E^u(x: r(x))) \subset E^u(x; \tr(x))$$
is $C^{1+\beta}$-close to identity. Thus, replacing $R(x)$ with a larger regularity function if  necessary, we may assume that
\be
\| \Psi^u_x(u) - u\| \leq R(x) \|u\|^{1+\beta} \quad , \quad \| (\Psi^u)^{-1}_x(u) - u\| \leq R(x) \|u\|^{1+\beta}
\ee
for all $ x\in \ll$ and  $u \in E^u(x;\tr(x))$, and also  that 
\be
\|d \Phi^u_x(u)\| \leq R(x) \quad ,\quad  \|(d\Phi^u_x(u))^{-1}\| \leq R(x) \quad,  
\quad x\in \ll \:, \: u \in E^u (x ; r(x)) .
\ee
Finally, again replacing $R(x)$ with a larger regularity function if  necessary, we may assume that
\be
\quad\|\Phi^u_x(v) - \Phi^u_x(u) - d\Phi^u_x(u) \cdot (v-u)\| 
\leq R(x) \, \|v-u\|^{1+\beta} \:\:,\:  x\in \ll\:, \: u,v \in E^u (x ; r(x)) ,
\ee
and
\be
\|d\Phi^u_x(u) - \id\| \leq R (x)\, \|u\|^\beta \quad , \quad x\in \ll \:, \: u \in E^u (x; r(x)) .
\ee
In a similar way one defines the maps $\Phi^s_x$ and we will assume that $r(x)$ is chosen so that 
these maps satisfy the analogues of the above properties.

For any $x\in \ll$ consider the $C^{1+\beta}$ map (defined locally near $0$)
$$\hf_x = (\Phi^u_{f(x)})^{-1} \circ f\circ \Phi^u_x : E^u (x)  \longrightarrow E^u (f(x))\;.$$
It is important to notice that
$$\hf^{-1}_x(\hE^u_j (f(x) ; r(f(x))) \subset \hE^u_j (x; r(x)) \quad , \quad 
\hf^{-1}_x(\tE^u_j (f(x) ; r(f(x))) \subset \tE^u_j (x; r(x)) \;$$
for all $x \in \ll$ and $j > 1$.

Given $y \in \ll$ and any integer $j \geq 1$ we will use the notation
$$\hf_y^j = \hf_{f^{j-1}(y)} \circ \ldots \circ \hf_{f(y)} \circ \hf_y\quad,
\quad \hf_y^{-j} = (\hf_{f^{-j}(y)})^{-1} \circ \ldots \circ (\hf_{f^{-2}(y)})^{-1} 
\circ (\hf_{f^{-1}(y)})^{-1} \;,$$
at any point where these sequences of maps are well-defined.

It is well known (see  e.g. the Appendix in \cite{LY1} or Sect. 3 in \cite{PS})  that there exists a 
Lyapunov $\hep$-regularity functions
$\Gamma = \Gamma_{\hep} : \ll \longrightarrow [1,\infty)$ and $r = r_{\hep} : \ll \longrightarrow (0,1)$ and 
for each $x\in \ll$ a norm $\| \cdot \|'_x$ on $T_xM$ such that
\be
\|v\| \leq \|v\|'_x \leq \Gamma (x) \|v\| \quad ,\quad x\in \ll\:,\: v \in T_xM ,
\ee
$$\angle(\hE^u_j (x), \tE^u_{j}(x)) \geq \frac{1}{R(x)} \quad , \quad x\in \ll\;, \; 2 \leq j \leq d ,$$
and for any $x \in \ll$ and any integer $m \geq 0$, assuming 
$\hf_x^j(u), \hf_x^j(v) \in E^u(f^j(x), r (f^j(x)))$ are well-defined for all 
$j =1, \ldots,m$, the following hold:
\be
\mu^m_j\, \|u-v\|'_x  \leq \|\hf_x^m(u) - \hf_x^m(v)\|'_{f^m(x)} \quad , \quad u,v \in \tE^u_j(x; r(x)) ,
\ee
\be
\mu^m_1\, \|u-v\|'_x  \leq \|\hf_x^m(u) - \hf_x^m(v)\|'_{f^m(x)}  \quad , \quad u,v \in E^u(x; r(x))  ,
\ee
\be 
\mu_1^m \, \|v\|'_x \leq \|d\hf_x^m(u)\cdot v\|'_{f^m(x)} \leq \nu_d^m \, \|v\|'_x
\quad , \quad x\in \ll \:,\: u \in E^u(x;r(x))\;,\; v \in E^u(x) ,
\ee
\be 
\mu_j^m \, \|v\|'_x \leq \|d\hf_x^m(0)\cdot v\|'_{f^m(x)} \leq \nu_j^m \, \|v\|'_x
\quad , \quad x\in \ll \:,\:  v \in E^u_j(x) .
\ee

We will also use the norm  $|u| = \max\{ \|\ui\| : 1 \leq i \leq \tk\}$. Clearly, 
$$\|u\| = \| u^{(1)} + \ldots + u^{(\tk)}\| \leq \sum_{i=1}^{\tk} \|\ui\| \leq \tk \, |u| .$$
Taking the regularity function $\Gamma(x)$ appropriately,  
we have $|u|\leq \Gamma(x) \|u\|$, so
\be
\frac{1}{\tk}\, \|u\| \leq |u| \leq \Gamma (x) \|u\| \quad , \quad x\in \ll\;,\; u\in E^u(x) .
\ee


Next, Taylor's formula (see also section 3 in \cite{PS}) implies that there exists 
a Lyapunov $\hep$-regularity function  $D = D_{\hep} : \ll \longrightarrow  [1,\infty)$ such that 
for any $i = \pm 1$ we have
\be
\quad\|\hf_x^i(v) - \hf_x^i(u) - d\hf^i_x(u) \cdot (v-u)\| \leq D(x) \, \|v-u\|^{1+\beta} \:\:,\: 
x\in \ll\:, \: u,v \in E^u (x ; r(x)) ,
\ee
and
\be
\|d\hf^i_x(u) - d\hf^i_x(0)\| \leq D (x)\, \|u\|^\beta \quad , \quad x\in \ll \:, \: u \in E^u (x; r(x)) .
\ee

Finally, we state here a Lemma from \cite{St3} which will be used several times later.

\bs

\noindent
{\bf Lemma 3.1.} (Lemma 3.3 in \cite{St3}) {\it There exist a Lyapunov $6\hep$-regularity function\\
$L = L_{6\hep}: \ll \longrightarrow [1,\infty)$ and a Lyapunov $7\hep/\beta$-regularity 
function $r = r_{7\hep/\beta}: \ll \longrightarrow (0,1)$
such that for any $x\in \ll$, any integer $p \geq 1$ and any $v\in E^u(z, r (z))$ with 
$\|\hf_z^p(v)\|\leq r(x)$,  where $z = f^{-p}(x)$, we have
$$\|\wo_p - \vo_p\| \leq L(x) |v_p|^{1+\beta} ,$$
where $v_p = \hf^p_z(v) \in E^u (x)$ and $w_p = d\hf_z^p(0)\cdot v \in E^u (x)$.
Moreover, if $|v_p| = \|\vo_p\| \neq 0$, then}  $1/2 \leq \|\wo_p\|/\|\vo_p\| \leq 2$.

\bs

\noindent
{\bf Remark.} Notice that if $v\in E^u_1(z, r(z))$ in the above lemma, then
$v_p, w_p \in E^u_1(x)$, so $\|w_p - v_p \| \leq L(x)\, \|v_p\|^{1+\beta}$.

\ms


Let $\Phi : M \longrightarrow \R$ be a H\"older continuous functions as in Sect. 3.1 and let $\m$ be the 
{\it Gibbs measure} determined by $\Phi$.
Let $\rr = \{ R_i\}_{i=1}^{k_0}$  be a pseudo-Markov family for $\phi_t$ as in Sect. 2, and let
$\tau: R = \cup_{i=1}^{k_0} R_i \longrightarrow [0,1/2]$ and $\pp : R\longrightarrow R$ be the 
corresponding first return map and the Poincar\'e map. As before fix constants $0 < \tau_0 < \htau_0 \leq 1/2$ 
so that $\tau_0 \leq \tau(x) \leq \htau_0$  for all $x\in R$. The Gibbs measure $\m$ induces a {\it Gibbs measure} 
$\mu$ on $R$ (with respect to the Poincar\'e map $\pp$) for the function
$$F(x) = \int_0^{\tau(x)} \Phi(\phi_s(x))\, ds \quad, \quad x\in R .$$
The function $F$ is H\"older and, using Sinai's Lemma, it is cohomologous to a {\it H\"older function 
$f: R \longrightarrow \R$  which is constant on stable leaves in rectangles} $R_i$ in $R$. Thus, $\mu$ 
coincides with the Gibbs measure determined by $f$.
For every continuous function $H$ on $M$ we then have (see e.g. \cite{PP})
\be
\int_M H\, d\m = \frac{\int_R \left( \int_0^{\tau(x)} H(\phi_s(x))\, ds\right) d\mu(x)}{\int_R \tau\, d\mu } .
\ee


Given a Lyapunov regularity function $R_\ep$ with (3.3) and (3.4), any set of the form
$$Q_{p}(\ep) = \{ x\in \ll : R_{\ep}(x) \leq e^{p} \} $$
is called a  {\it Pesin set}. Given $p > 0$, $\ep > 0$, $\delta  > 0$ and an integer $n \geq 1$ 
set
\be
\Xi_n = \Xi_n(p,\ep,\delta) = \left\{ x\in \ll \cap R: \:\: \sharp \, 
\left\{ j : 0 \leq j \leq n-1 \: \mbox{\rm and } \:  \pp^j(x) \notin Q_{p}(\ep) \right\} 
\geq \delta\, n \: \: \right\}.
\ee

\ms

\noindent
{\bf Definition.} (\cite{GS})
Consider a log-integrable linear cocycle $M$ above
a transformation $(T,\mu)$, with Lyapunov exponents $\lambda_1 \geq \ldots
\geq \lambda_d$. We say that $M$ has {\it exponential large deviations for all exponents} if, for
any $i\leq d$ and any $\epsilon>0$, there exists $C>0$ such that, for all $n\geq 0$,
\begin{equation}
\label{eq:exp_dev_all_exp}
  \mu\{  x : |\log \|\Lambda^i M^n(x)\|
     - n (\lambda_1+\ldots+\lambda_i)| \geq n \epsilon\} \leq C e^{- n/C}.
\end{equation}

The following theorem, which is a special case of Theorem 1.7 in \cite{GS},  
shows that if $df$ has exponential large deviations 
for all exponents, then most points in $\ll$ return exponentially often to some Pesin set. 

\bs

\noindent
{\bf Theorem 3.2.} (\cite{GS})
{\it  Assume that $df$ has exponential large deviations for all exponents with respect to $\mu$.
Let $\hep_0 > 0$ and $\hd_0 > 0$. Then there exist $p_0 > 0$,  $C > 0$ and $c  > 0$ such that
$$\m \left( \left\{ x\in \ll : \:\: \sharp \, 
\left\{ j : 0 \leq j \leq n-1 \: \mbox{\rm and } \:  f^j(x) \notin Q_{p_0}(\hep_0) \right\} \geq \hd_0 n \right\}\right)
\leq C e^{- c n} ,$$
for all $n\geq 1$. Thus, there exist constants $p > 0$,  $C' > 0$ and $c' > 0$ such that
\begin{equation}\label{eq:tailestimate}
\mu (\Xi_n(p_0,\hep_0,\hd_0))\leq C' e^{- c' n}
\end{equation}
for all $n \geq 1$. }

\bs

Clearly, if (3.21) holds for $p_0$, then it will hold with $p$ replaced by any $p \geq p_0$.

As established in \cite{GS} (see Theorem 1.5 there), for a transitive subshift of finite type
$T$ on a space $\Sigma$, if $\mu$ is a Gibbs measure for a H\"older-continuous
potential and $M$ is a continuous linear cocycle on a vector bundle $E$
above $T$, each of the following conditions is sufficient for $M$ to have exponential large deviations 
for all exponents: (i) if all its Lyapunov exponents coincide;
(ii) if there us a continuous decomposition of $E$ as a direct sum of
    subbundles $E=E_1 \oplus \ldots \oplus E_k$ which is invariant under
    $M$, such that the restriction of $M$ to each $E_i$ has exponential
    large deviations for all exponents;
(iii) more generally, if there is an invariant continuous flag
    decomposition $\{0\}=F_0 \subseteq F_1 \subseteq \ldots \subseteq F_k =
    E$, such that the cocycle induced by $M$ on each $F_i/F_{i+1}$ has
    exponential large deviations for all exponents;
(iv) if the cocycle $M$ is locally constant in some trivialization of the
    bundle $E$ (this is equivalent to the existence of invariant continuous
    holonomies which are commuting);
(v) if the cocycle $M$ admits invariant continuous holonomies, and if it
    is pinching and twisting in the sense of Avila-Viana \cite{AV};
(vi) if the cocycle $M$ admits invariant continuous holonomies, and the
    bundle is $2$-dimensional.

It follows from the above and Theorem 9.18 in \cite{V} that generic linear cocycles have exponential large 
deviations for all exponents. Moreover,  amongst fiber bunched cocycles\footnote{Which are the most frequently 
met H\"older cocycles.} generic cocycles in the H\"older topology also have exponential large deviations for all exponents.

\bs

Finally we will make an important remark concerning Lyapunov regularity functions on $\ll$ (see Sect. 3.1).
A Lyapunov $\ep$-regularity function $H(x)$, $x\in \ll$, will be called a {\it large } (resp. {\it small}) 
{\it canonical $\ep$-regularity function} if there exist constants $p > 0$ and $H_0 > 0$ such that 
$$1 \leq H(x) \leq H_0\, (R_\ep(x))^p \quad \mbox{\rm (resp. $1 \geq H(x) \geq \frac{1}{H_0\, (R_\ep(x))^{p}}$)}$$
 for all $x\in \ll$, where $R_\ep(x)$ is given by (3.18). It follows from the above that
in the case of a large (reps. small) canonical $\ep$-regularity function $H$ there exists a constant $H'_0$ so that
$$1 \leq H(x) \leq H'_0\, (\hD_{\ep/4}(x))^p  \quad \mbox{\rm (resp. 
$1 \geq H(x) \geq \frac{1}{H_0\, (\hD_{\ep/4}(x))^{p}}$) .}$$

\ms

\noindent
{\bf Proposition 3.3.} {\it The function $r_{\hep}$ in Sect. 3.1 can be chosen to be a small canonical 
$\hep$-regularity function, while $\Gamma(x)$ and $D(x)$ can be chosen to be large canonical $\hep$-regularity 
functions. Moreover, we can always choose $r_\ep(x)$ so that 
$r_\ep(x) \leq \frac{1}{(R_\ep(x))^6}$
for all $x\in \ll \cap R$.}

\bs
 
This follows e.g. from the arguments in Sect. 3 in \cite{PS}. Moreover, following the arguments in \cite{St3},
the function $L(x)$ in Lemma 3.1 is a large canonical regularity function, and similarly the arguments
in Sect. 9 show that all regularity functions constructed there are canonical regularity functions.

\def\chR{\check{R}}
\def\llm{\ll^{(m)}}
\def\llk{\ll^{(k)}}

\section{Non-integrability of Anosov flows}
\setcounter{equation}{0}

\subsection{Choice of constants, sets of Lyapunov regular points}


In what follows we assume that $\tRR = \{ \tR_i\}_{i=1}^{k_0}$  is a fixed Markov partition for 
$\phi_t$ on $M$ of size $< 1/2$ and $\rr = \{ R_i\}_{i=1}^{k_0}$  is the related pseudo-Markov 
family as in Sect. 2. We will use the notation associated with these from Sect. 2, and
we will assume that for any $i = 1, \ldots,k_0$, $z_i$ is chosen so that $z_i \in \Intu(W^u_{R_i}(z_i))$.
For any $x\in R$, any $y \in \tR$ and $\delta > 0$ set 
$$B^u(x,\delta) = \{ y\in W^u_{R_i}(x) : d(x,y) < \delta\}  \quad , \quad
\tB^u(y,\delta) = \{ z\in W^u_{\tR_i}(z) : d(z,y) < \delta\} .$$
In a similar way define $B^s(x,\delta)$. For brevity sometimes we will use the notation 
$U_i(z) = W^u_{R_i}(z)$
for $z\in  R_i .$ 

Fix constants $0 < \tau_0 < \htau_0 < 1$ so that $\tau_0 \leq \tau(x) \leq \htau_0$ 
for all $x\in R$ and $\tau_0 \leq \ttau(x) \leq \htau_0$ for all $x\in \tR$.



Let $\alpha_1 > 0$ be as in Sect. 2, and let $f$ be an essentially 
$\alpha_1$-H\"older continuous potential
on $R$.  Set $g = f - P_f \tau$, where $P_f \in \R$ is chosen so that
the topological pressure of $g$ with respect to the Poincar\'e map 
$\pp : R = \cup_{i=1}^{k_0} R_i \longrightarrow R$  is $0$.
Let $\mu = \mu_g$ be the Gibbs measure on $R$ determined by $g$; then $\mu(R) = 1$.
We will assume that $f$ (and therefore $g$) depends on forward coordinates only\footnote{If the initial 
potential $F$ on $R$ is $\alpha^2$-H\"older, applying Sinai's Lemma 
(see e.g. \cite{PP}) produces an $\alpha$-H\"older 
potential $f$ depending on forward coordinates only.}  i.e. it is constant on stable leaves of 
$R_i$ for each $i$. 

 
Since $g$ is constant on stable leaves, it generates a {\it Gibbs measure} $\nu^u$ on $U$. 
Let $g^{(s)}$ be a function on $R$ which is homological to $g$ and constant on unstable 
leaves in $R$; then $g^{(s)}$ can be regarded as a function on $S$
and determines a {\it Gibbs measure} $\nu^s$ on $S$. 

A sequence $i_p, i_{p+1}, \ldots, i_{q}$ of elements of $\{ 1, \ldots, k_0\}$ 
for some integers $p \leq q$, will be called {\it admissible} if 
$\tpp(\Int(\tR_{i_j})) \cap \Int(\tR_{i_{j+1}}) \neq \e$ for all $j = p, p+1, \ldots, q-1$.
Given such a sequence, consider the {\it cylinder}
$$\cc_R[i_p, i_{p+1}, \ldots, i_{q}] = \{ x\in R_{i_p} : \pp^j(x) \in R_{i_{p+j}} \:, \: 1 \leq j \leq q-p\} $$
in $R$. When $p = 0$, we can define similarly the usual `unstable' cylinders in $U_{i_0}$:
$$\cc^u[i_0, \ldots, i_q] = \{ x\in U_{i_0} : \pp^j(x) \in R_{i_{j}} \:, \: 1 \leq j \leq q\} .$$
and for $q = 0$ we define a `stable' cylinder in $S_{i_0}$:
$$\cc^s[i_p, \ldots, i_0]  = \{ x\in S_{i_0} : \pp^j(x) \in R_{i_{j}} \:, \: p \leq j \leq 0\} .$$
Then (see Proposition A2.2 in \cite{P2} or Sect. 2.3 in \cite{Ch3}) there exist constants $0 < A_1 < A_2$ 
such that for every cylinder $\cc_R$ as above with $p \leq 0 \leq q$ we have
$$A_1 \leq \frac{\mu(\cc_R[i_p, i_{p+1}, \ldots, i_{q}])}{\nu^s(\cc^s[i_p, \ldots, i_0] ) 
\, \nu^u(\cc^u[i_0, \ldots, i_q])} \leq A_2 .$$
Moreover we have
$\nu^u(\cc^u[i_0, \ldots, i_q]) = \mu (\cc_R[i_0, \ldots, i_q])$, and
$\nu^s(\cc^s[i_p, \ldots, i_0]) = \mu(\cc_R[i_p, \ldots, i_0])$.

It follows from the above that $\mu$ is almost the direct product of $\nu^u$ and $\nu^s$. More precisely,
let $\hmu$ be {\it the probability measure} on $R$ such that $\hmu = \nu^u\times \nu^s$ 
on each $R_i$, where we use the 
natural (Borel measurable) isomorphism $R_i = [U_i,S_i] \approx U_i \times S_i$. 
It then follows from the above that  for every bounded Borel measurable function 
$H$ on $R$ and every $i = 1, \ldots, k_0$ we have
\be
A_1\, \int_{U_i} \int_{S_i} H([x,y])\, d\nu^u(x)\, d\nu^s(y) \leq \int_R H\, d\mu 
\leq A_2 \, \int_{U_i} \int_{S_i} H([x,y])\, d\nu^u(x)\, d\nu^s(y) .
\ee
For later convenience, for every $i$ and every $z\in R_i$ we will denote by 
$\nu^u_z$  {\it the measure on} $W^u_{R_i}(z)$
determined by $\nu^u$ and the  projection $\pi_{z} : U_i \longrightarrow W^u_{R_i}(z)$ 
along stable manifolds in $R_i$, i.e.
$\nu^u_z(\pi_z(A)) = \nu^u(A)$ for every Borel measurable subset $A$ of $U_i$. In a similar way we 
define\footnote{In general $\nu^u_z$ and $\nu^s_z$ are not the conditional measures 
determined by $\mu$.} $\nu^s_z$.

Given an unstable leaf $W = W^u_{\tR_i}(z)$ in some rectangle $\tR_i$ 
and an admissible sequence $\ii = i_0, \ldots,i_m$ of integers $i_j \in \{ 1,\ldots, k_0\}$, the set
$$C_W[\ii] = \{ x\in W : \tpp^j(x) \in \tR_{i_j} \:, \: j = 0,1,\ldots,m\}$$
will be called a {\it cylinder of length} $m$ in $W$ (or an {\it unstable cylinder} in $\tR$ in general). 
When $W = U_i$ we will simply write $C[\ii]$. In a similar way one defines cylinders $C_V[\ii]$, 
where $V = W^u_{R_i}(z)$ is an unstable leaf in some rectangle $R_i$. 

Let
$\pr_D : \cup_{i=1}^{k_0} \phi_{[-\ep,\ep]}(D_i) \longrightarrow \cup_{i=1}^{k_0} D_i$
be the {\it projection along the flow}, i.e. for all $i = 1,\ldots,k_0$ and all 
$x\in \phi_{[-\ep,\ep]}(D_i) $ we have $\pr_D(x) = \pr_{D_i}(x)$ (see Sect. 2). For any $z\in R$ denote by 
$\chU(z)$ the {\it part of the unstable manifold} $W^u_{\ep_0}(z)$ such that $\pr_D(\chU(z)) = W^u_{\tR_i}(z)$.
The shift along the flow determines a bi-H\"older continuous bijections 
$\T_z :  W^u_{R}(z) \longrightarrow \chU(z)$ 
and $\tPsi : W^u_{R}(z) \longrightarrow W^u_{\tR}(z)$ for all $i$. 
These define bi-H\"older continuous bijections   
$\Psi: R \longrightarrow \check{R} = \cup_{i=1}^{k_0} \check{R}_i$, 
where $\check{R}_i = \cup_{z\in S_i} \chU(z)$ and
$\Psi_{| W^u_R(z)} = (\T_z)_{|W^u_R(z)}$ for $z\in S_i$, and $  \tPsi : R \longrightarrow \tR$. 
Notice that there exists a global constant $C > 1$ such that 
$\frac{1}{C} d(x,y) \leq d(\T_z(x), \T_z(y)) \leq C\, d(x,y)$ for any $z\in \tR$ and any 
$x,y \in W^u_{\tR}(z)$.


\def\Xim{\Xi^{(m)}}
\def\Omm{\Omega^{(m)}}
\def\Xin{\Xi^{(n)}}
\def\Omn{\Omega^{(n)}}


\subsection{ Regular distortion of cylinders}

In \cite{St4} we established some nice properties concerning diameters of cylinders for 
Axiom A flows on basic sets satisfying a pinching condition which we called {\it regular distortion 
along unstable manifolds}. In \cite{St3} something similar was established for Anosov flows with 
Lipschitz local stable holonomy maps. It seems unlikely that any Anosov flow will 
have such properties, however it turns out that for general Anosov flows something similar holds 
for cylinders in $R$ that intersect `at both ends' a compact set of Lyapunov regular points with bounded
from below regularity functions $r(x)$ (i.e. a Pesin set).  More precisely we have the following.

\bs

\noindent
{\bf Lemma 4.1.} 
(a) {\it   There exists a constant $0 < \rho_1  < 1$ 
such that for any unstable leaf $W$ in $R$, any cylinder $C_W[\ii] = C_W[i_0, \ldots,i_m]$  in $W$ 
and any sub-cylinder $C_W[\ii'] = C_W[i_0,i_{1}, \ldots, i_{m+1}]$ of $C_W[\ii]$ of co-length $1$ 
such that there exists $z\in C_W[\ii']$ with $\pp^m(z) \in P_0$  we have}
$$\rho_1 \; \diam ( \tPsi(C_W [\ii]) ) \leq  \diam ( \tPsi(C_W [\ii'] )) .$$

\ms

(b) {\it For any constant $\rho' \in (0,1)$ there exists an integer $q' \geq 1$ such that
for  any unstable leaf $W$ in $R$, any cylinder $C_W[\ii] = C_W[i_0, \ldots,i_m]$ of length $m$ in $W$ 
and any sub-cylinder $C_W[\ii'] = C_W [i_0,i_{1}, \ldots, i_{m+1}, \ldots, i_{m+q'}]$ of $C_W[\ii]$ of co-length $q'$ 
such that there exists $z\in C_W[\ii'] $ with $\pp^{m+q'}(z) \in P_0$ we have}
$$\diam(\tPsi(C_W[\ii']) ) \leq \rho' \, \diam (\tPsi(C_W[\ii])) .$$

(c)  {\it There exist an integer $q_0 \geq 1$ and a constant $\rho_1 \in (0,1)$ 
such that for any unstable leaf $W$ in $R$  and any cylinder $C_W[\ii] = C_W[i_0, \ldots,i_m]$ in $W$ 
such that there exists $z\in C_W[\ii']$ with $\pp^m(z) \in P_0$
there exist points $z, x\in C_W[\ii]$ such that if $C_W[\ii'] = C_W[i_0,i_{1}, \ldots, i_{m+1}, \ldots, i_{m+q_0}]$ 
is the sub-cylinder of  $C_W[\ii]$ of co-length $q_0$
containing $x$ then $d(z,y) \geq \rho_1 \, \diam (\tPsi(C_W[\ii]))$ for all $y \in C_W[\ii']$. }

\bs

This Lemma will be used essentially in the proof of the main result in Sects. 5-7 below. 
Its proof is given in Sect. 9.

\subsection{Non-integrability}

Throughout we assume that $\phi_t$ is a $C^2$ contact Anosov flow on $M$ with a $C^2$ invariant  contact
form $\omega$. Then the two-form $d\omega$ is $C^1$, so there exists a constant $C_3 > 0$ such that
\be
|d\omega_x(u,v)| \leq C_0 \|u\|\, \|v\| \quad , \quad u, v \in T_xM\:, \: x\in M .
\ee
Moreover, there exists a constant $\theta_0 > 0$ such that for any $x\in M$ and any $u\in E^u(x)$ with $\|u\| = 1$
there exists $v\in E^s(x)$ with $\|v\| = 1$ such that $|d\omega_x(u,v)| \geq \theta_0$.

The main ingredient in this section is the following lemma of Liverani (Lemma B.7 in \cite{L1})
which significantly strengthens a lemma of Katok and Burns (\cite{KB}).

\bs

\noindent
{\bf Lemma 4.2.} (\cite{L1}) {\it Let $\phi_t$ be a $C^2$ contact flow on $M$ with a 
$C^2$ contact form $\omega$.
Then there exist constants $C_0  > 0$, $\vartheta > 0$ and $\hep_0 > 0$ such that
for any $z\in M$, any $x\in  W^u_{\hep_0}(z)$ and any  $y\in W^s_{\hep_0}(z)$ we have
\be\label{eq:liverbound}
|\Delta(x,y)  - d\omega_z(u,v)| \leq C_0\,\left[ \|u\|^2\, \|v\|^\vartheta + \|u\|^\vartheta \|v\|^2 \right]\;,
\ee
where $u \in E^u(z)$ and $v \in E^s(z)$ are such that $\exp^u_z(u) = x$ and $\exp^s_z(v) = y$.}

 \bs

\noindent
{\bf Note.}
Actually Lemma B.7 in \cite{L1} is more precise with a particular choice of the constant $\vartheta$
determined by the (uniform) H\"older exponents of the stable/unstable foliations and the corresponding
local holonomy maps. However in this paper we do not need this extra information.

\ms

{\bf From now on we will assume that $C_0 > 0$, $\vartheta > 0$ and $\hep_0 \in (0,\ep_0/4)$ satisfy (4.2) and (4.3),}
where $\ep_0$ is as in Theorem 3.2.
Assume also that the constant $\hd_0 $ is fixed so that  $ \lambda_k^{\hd_0} < e^{\hep}$
 (some conditions will be listed later) with.   
As in Sect. 3, set
$$Q_p(\ep) = \{ x\in \ll  : R_{\ep}(x) \leq e^p \}$$
for all $\ep \in (0,\ep_0)$ and $p > 0$. Then $Q_0(\hep_0) \subset Q_1(\hep_)) \subset Q_n(\hep_0) \subset \ldots $
and $\cup_{p=0}^\infty Q_p(\hep_0) = \ll$.
{\bf Fix an integer} $p_0 \geq 1$ so large that $\mu(Q_{p_0}(\hep_))) > 1- \delta$ for some small appropriately 
chosen $\delta > 0$ (to be determined later). Set 
$$\ll_0 = \cup_{p=0}^\infty Q_p(\hep_0) .$$ 
Then $\mu(\ll_0 \cap R) = 1$. 
Set
\be \label{eq:pesinset}
P_0 = Q_{p_0}(\hep_0) \quad , \quad \tP_0 = \phi_{[-1,1]}(P_0) .
\ee
Then the Lyapunov regularity function $R_{\hep}(x)$ is bounded by some constant on $P_0$ (and therefore on
$\tP_0$ as well), and according to Proposition 3.3, we may assume
\be\label{eq:constbounds}
R_{\hep}(x) \leq R_0 \quad, \quad r(x) \geq r_0 \quad, \quad \Gamma(x) \leq \Gamma_0 \quad , \quad
L(x) \leq L_0 \quad, \quad D(x) \leq D_0 \quad ,\quad x\in \tP_0 ,
\ee
for some positive constants $R_0, \Gamma_0, L_0, D_0 \geq 1$ and $r_0 > 0$. 

It follows easily from the properties of Markov 
families\footnote{Easy proof by contradiction.} that there exists a constant 
$r_1 > 0$  such that for every $i$ and every $x\in \partial R_i$  
there exists $y \in R_i$ such that $\dist(y, \partial R_i) \geq r_1$ and $d(x,y) < r_0/2$.
{\bf Fix a constant} $r_1 < \frac{r_0}{2R_0}$ with this property.

Take a large constant $L > 1$ (to be determined later) and set
\be\label{eq:Xin}
\Xi_n = \Xi_n(p_0,\hep_0, \hd_0) \quad, \quad
\Xi^{(n)}_L = \ll\setminus \left(\cup_{n/L \leq \ell \leq L n} \Xi_\ell \right)  .
\ee
It follows from  Theorem 3.2 that, choosing the  constants $C_0, c_0 > 0$ appropriately (depending on $L$), we have
\be\label{eq:XinB}
\mu(R\setminus \Xin_L) \leq C_0 \, e^{-c_0 n/L}
\ee
 for all  $n \geq 1$.
 

We will show below that for Lyapunov regular points $x\in \ll_0$ the estimate (4.3) can be 
improved what concerns the involvement of $u$ for certain choices of $u$ and $v$. More precisely, 
we will show that choosing $v$ in a special way, $\Delta(x,y)$ becomes a $C^1$ function of 
$x = \exp^u_z(u)$ with a non-zero uniformly bounded derivative in a certain direction. 


We will now state two Main Lemmas. Their proofs, both using Liverani's Lemma 4.2, are given in Sect. 8.

We will assume that $L$ is a {\bf fixed constant } with $L > 3/\tau_0$.


\bs


\noindent
{\bf Lemma 4.3.}  {\it There exist constants $C_1 > 0$ and $\beta_1 \in (0,1)$ with the 
following properties:}

(a)  {\it For any unstable cylinder $\cc$ in $R$ of length $m$ with $\cc  \cap P_0\cap \Xi_m \neq \e$ 
and any $z\in \cc$ we have
\be\label{eq:cylbound}
\frac{1}{C_1 \lambda^p_1 e^{2\hep p}} \leq \diam(\tPsi(\cc)) \leq \frac{C_1 e^{2\hep p}}{\lambda_1^p} ,
\ee
where $p = [\tau_m(z)]$.}

\ms 

(b) {\it For any unstable cylinder $\cc$ of length $m$ in $R$ with 
$\cc \cap P_0\cap \Xim_L \neq \e$,
any $\hx_0,\hz_0 \in \cc$ and any $\hy_0, \hb_0 \in W^s_{R}(\hz_0)$ we have
$$|\Delta(\hx_0 ,\hy_0)  - \Delta(\hx_0, \hb_0)| \leq C_1 \, \diam (\tPsi(\cc)) \, (d(\hy_0,\hb_0))^{\beta_1} .$$
In particular,}
$$|\Delta(\hx_0 , \hy_0)| \leq C_1 \, \diam (\tPsi(\cc)) \, (d(\hy_0,\hz_0))^{\beta_1}  \leq C_1 \, \diam (\tPsi(\cc)) .$$



\bs

{\bf Fix a constant} $C_1 > 0$ with properties in Lemma 4.3. We take $C_1 \geq C_0$.
Set  $\di \beta_0 = \frac{1}{\sqrt{1+ \theta_0^2/(64 C_1^2)}}$.
Next, {\bf fix an integer} $\ell_0 = \ell_0(\delta) \geq 1$ so large that we can find
unit vectors $\eta_1, \eta_2 , \ldots,\eta_{\ell_0}$ in $\R^{n_1}$ such that for any unit vector $\xi \in \R^{n_1}$ 
there exists $j$ with $\la \xi , \eta_j\ra \geq \beta_0$. Then {\bf fix measurable families}
$\eta_1(x), \eta_2(x), \ldots, \eta_{\ell_0}(x)$ ($x\in \ll_0$) of unit vectors in $E^u_1(x)$ such that for any
$x\in \ll_0$ and any $\xi\in E^u_1(x)$ with $\|\xi\| = 1$ there exists $j$ with  $\la \xi, \eta_j(x) \ra \geq \beta_0$.
Recall the projections $\T_z : W^u_R(z) \longrightarrow \chU(z) \subset W^u_{\ep_0}(z)$ for $z\in R$.

\bs

\noindent
{\bf Lemma 4.4.} {\it Let $\phi_t$ be a $C^2$ contact Anosov flow on $M$. Let  
$\eta_1(x), \eta_2(x), \ldots, \eta_{\ell_0}(x)$ ($x\in \ll_0$) be families of unit vectors in $E^u_1(x)$ as above, 
and let $\kappa \in (0,1)$ be a constant.  Then there exist constants  $\ep'' > 0$, 
$0 < \delta'' < \delta'$ (depending on $\kappa$ in general), $\delta_0 \in (0,1)$,  
 with the following properties:}

\ms

(a) {\it  For any integer $m \geq 1$ 
and any $Z \in P_0 \cap \Xim_L$ there exist families of points 
$y_j(Z) \in  B^s (Z,\delta')$ ($ j = 1, \ldots,\ell_0$) such that 
if $\cc$ is a cylinder of length $m$ in $R$ with $Z \in \cc$, then for any $x_0  \in \T_Z(\cc)$,  
$z_0  \in \T_Z(\cc \cap P_0)$ of the form $x_0 = \Phi^u_{Z}(u_0)$,  $z_0 = \Phi^u_{Z}(w_0)$ such that
\be
d(x_0,z_0) \geq \kappa \, \diam (\cc') , 
\ee 
where $\cc' = \T_z(\cc)$, and
\be
\left\langle \frac{w_0 - u_0}{\|w_0 - u_0\|} , \eta_j(Z) \right\rangle \geq \frac{\beta_0}{2R_0} 
\ee
for some $j = 1, \ldots, \ell_0$, then we have
\be \label{eq:deltabound}
\frac{\beta_0 \delta_0 \kappa}{16R^2_0} \, \diam(\cc') 
\leq |\Delta(x_0,\pi_{d_1}(z_0))  - \Delta(x_0, \pi_{d_2}(z_0))| 
\ee
for any $d_1 \in B^s(y_j (Z), \delta'')$ and  $d_2 \in B^s(Z, \delta'')$.}

\ms

(b)  {\it There exists an integer $N_0 \geq 1$ such that for any integer $N \geq N_0$, 
any integer $m \geq 1$ 
and any $Z \in P_0 \cap \Xim_L$ there exist families of points 
$$\yjo(Z), \yjt(Z) \in \pp^N(B^u(Z;\ep'')) \cap B^s (Z,\delta')
 \quad , \quad j = 1, \ldots,\ell_0 ,$$
such that if $\cc$ is a cylinder of length $m$ in $R$ with $Z \in \cc$, $x_0  \in \T_Z(\cc)$,  
$z_0  \in \T_Z(\cc \cap P_0)$ have the form $x_0 = \Phi^u_{Z}(u_0)$,  $z_0 = \Phi^u_{Z}(w_0)$ and
(4.10) and (4.11) hold for some $j = 1, \ldots, \ell_0$, then (4.12) holds
for any $d_1 \in B^s(\yjo (Z), \delta'')$ and  $d_2 \in B^s(\yjt(Z), \delta'')$.}



\def\tomega{\tilde{\omega}}
\def\vlim{v^{(\ell)}_{i,m}}
\def\tvlim{\tilde{v}^{(\ell)}_{i,m}}
\def\thetlim{\theta^{(\ell)}_{i,m}}
\def\lamblim{\lambda^{(\ell)}_{i,m}}
\def\vl{v^{(\ell)}}
\def\vlik{v^{(\ell)}_{i,k}}
\def\vlom{v^{(\ell)}_{1,m}}
\def\vltm{v^{(\ell)}_{2,m}}
\def\tvlik{\tilde{v}^{(\ell)}_{i,k}}
\def\tvlom{\tilde{v}^{(\ell)}_{1,m}}
\def\tvltm{\tilde{v}^{(\ell)}_{2,m}}
\def\Cum{\cc^{(u_m)}}
\def\Cu{\cc^{(u)}}
\def\gao{\gamma^{(1)}}
\def\gat{\gamma^{(2)}}
\def\diamtef{{\footnotesize\mbox{\rm  diam}_\theta}}
\def\tcc{\tilde{\cc}}
\def\halpha{\hat{\alpha}}
\def\hgamma{\hat{\gamma}}
\def\tXijl{\widetilde{X}^{(\ell)}_{i.j}}
\def\tdd{\widetilde{\dd}}
\def\hcc{\widehat{\cc}}
\def\lambdam{\lambda^{(m)}}
\def\thetam{\theta^{(m)}}
\def\tomijl{\tilde{\omega}^{(\ell)}_{i,j}}
\def\jj{{\cal J}}
\def\oloij{\omega^{(\ell,0)}_{i,j}}
\def\olooj{\omega^{(\ell,0)}_{1,j}}
\def\olotj{\omega^{(\ell,0)}_{2,j}}
\def\olkij{\omega^{(\ell,k)}_{i,j}}
\def\olkoj{\omega^{(\ell,k)}_{1,j}}
\def\olktj{\omega^{(\ell,k)}_{2,j}}
\def\oloim{\omega^{(\ell,0)}_{i,m}}
\def\oloom{\omega^{(\ell,0)}_{1,m}}
\def\olotm{\omega^{(\ell,0)}_{2,m}}
\def\olkim{\omega^{(\ell,k)}_{i,m}}
\def\olkom{\omega^{(\ell,k)}_{1,m}}
\def\olktm{\omega^{(\ell,k)}_{2,m}}
\def\vlo{v^{(\ell)}_{1}}
\def\vlt{v^{(\ell)}_2}
\def\vli{v^{(\ell)}_{i}}
\def\vloij{v^{(\ell,0)}_{i,j}}
\def\vlooj{v^{(\ell,0)}_{1,j}}
\def\vlotj{v^{(\ell,0)}_{2,j}}
\def\vlkij{v^{(\ell,k)}_{i,j}}
\def\vlkoj{v^{(\ell,k)}_{1,j}}
\def\vlktj{v^{(\ell,k)}_{2,j}}

\def\vloim{v^{(\ell,0)}_{i,m}}
\def\vloom{v^{(\ell,0)}_{1,m}}
\def\vlotm{v^{(\ell,0)}_{2,m}}
\def\vlkim{v^{(\ell,k)}_{i,m}}
\def\vlkom{v^{(\ell,k)}_{1,m}}
\def\vlktm{v^{(\ell,k)}_{2,m}}

\def\vlrim{v^{(\ell,r)}_{i,m}}
\def\vlrom{v^{(\ell,r)}_{1,m}}
\def\vlrtm{v^{(\ell,r)}_{2,m}}

\def\Xloij{X^{(\ell,0)}_{i,j}}
\def\Xlooj{X^{(\ell,0)}_{1,j}}
\def\Xlotj{X^{(\ell,0)}_{2,j}}
\def\Xlkij{X^{(\ell,k)}_{i,j}}
\def\Xlkoj{X^{(\ell,k)}_{1,j}}
\def\Xlktj{X^{(\ell,k)}_{2,j}}

\def\Eekj{\Ee^{(k)}_{j}}
\def\Eerj{\Ee^{(r)}_{j}}

\def\Wr{W^{(r)}}
\def\Wk{W^{(k)}}

\def\tvlo{\tilde{v}^{(\ell)}_{1}}
\def\tvlt{\tilde{v}^{(\ell)}_{2}}
\def\tvli{\tilde{v}^{(\ell)}_{i}}

\def\omegak{\omega^{(k)}}
\def\nnk{\nn^{(k)}}
\def\omegao{\omega^{(o)}}
\def\nno{\nn^{(o)}}


\section{Construction of a `contraction set' $K_0$}
\setcounter{equation}{0}

\subsection{Normalized Ruelle operators and the metric $\dte$}


Let the constants $C_0 > 0$, $c_0 > 0$, $1 < \gamma < \gamma_1$  be as in Sects. 2 and 4, so that
(2.1) and (4.7) hold. 
{\bf Fix a constant} $\theta$ such that
\be
\frac{1}{\gamma^{\alpha_1 }} 
= \hat{\theta} \leq \theta < 1,
\ee
where $\alpha_1 > 0$ is the constant chosen at the end of Sect. 2.
 
Recall the metric $\dte$ on $\hU$  and the space $\ff_\theta(\hU)$ from Sect. 1.1.
In the same way we define the distance $\dte(x,y)$ for $x,y$ in 
$W \cap \hR$.
Lemma 5.2 below shows that $\tau\in \ff_\theta(\hU)$.
For a non-empty subset  $A$ of $U$ (or some $W^u_R(x)$) let $\diamte(A)$ be the {\it diameter} 
of $A$ with respect to $\dte$.
 


Let $f \in \ff_\theta (\hU)$ be a fixed real-valued function  and let 
$g = f - P_f\, \tau$, where $P_f\in \R$ is such that $\Pr_{\sigma} (g) = 0$. 
Since $f$ is a H\"older continuous function on $\hU$,
it can be extended to a H\"older continuos function on $R$ which is constant on stable leaves. 

Set $\Fa = f - (P_f + a) \tau$.
By Ruelle-Perron-Frobenius' Theorem (see e.g. chapter 2 in \cite{PP}) for any real number $a$  
with $|a|$ sufficiently small, as an operator on $\ff_\theta(\hU)$, $L_{\Fa}$ has a 
{\it largest eigenvalue} 
$\lambda_{a}$ and there exists a (unique) regular probability measure $\hnu_a$ on $\hU$ with 
$L_{\Fa}^*\hnu_a = \lambda_a\, \hnu_a$, i.e.  
$\int L_{\Fa} H \, d\hnu_a = \lambda_a\, \int H\ d\hnu_a$ for any $H \in \ff_\theta(\hU)$.
Fix a corresponding (positive) eigenfunction $h_{a} \in \ff_\theta(\hU)$ such that $\int h_{a} \, d\hnu_a = 1$. 
Then $d\nu = h_0\, d\hnu_0$ defines a {\it $\sigma$-invariant  probability measure} $\nu = \nu^u$ on $U$, which is in fact
the Gibbs measure $\nu^u$ determined by $G$ on $U$ (see Sect. 4.1).
Since $\Pr_\sigma (f- P_f\tau) = 0$, it follows from the main properties of pressure 
(cf. e.g. chapter 3 in \cite{PP}) 
that $|\Pr_\sigma(\Fa)| \leq \|\tau \|_0 \, |a|$.  Moreover, for small $|a|$ the maximal eigenvalue
$\lambda_{a}$ and the eigenfunction $h_{a}$ are Lipschitz in  $a$, so
there exist constants $a'_0 > 0$ and $C' > 0$ such that $|h_{a} - h_0| \leq C'\, |a|$ on $\hU$ and
$|\lambda_a - 1| \leq C' |a|$ for  $|a| \leq a'_0$.

For $|a|\leq a'_0$,  as in \cite{D}, consider the function
$$\fa(u) = f (u) - (P_ f+ a) \tau(u) + \ln h_{a}(u) -  \ln h_{a}(\sigma(u)) - \ln \lambda_{a}$$
and the operators 
$$ \lab = L_{\fa - \i\,b\tau} : \ff_\theta(\hU) \longrightarrow \ff_\theta (\hU)\:\:\: , 
\:\:\: \ma = L_{\fa} : \ff_\theta(\hU) \longrightarrow \ff_\theta(\hU) .$$
One checks that $\ma \; 1 = 1$ and $\di |(\lab^m h)(u)| \leq (\ma^m |h|)(u)$ for all $u\in \hU$, 
$h\in \ff_\theta (\hU)$ and $m \geq 0$. It is also easy to check that $L_{f^{(0)}}^*\nu = \nu$, 
i.e.  $\int L_{f^{(0)}} H \, d\nu = \int H\, d\nu$ for any $H \in \ff_\theta (\hU)$ 
(in fact, for any bounded continuous function $H$ on $\hU$).

Since $g$ has zero topological pressure with respect to the shift
map $\sigma : U \longrightarrow U$, there exist constants $0 < c_1 \leq c_2$  
such that for any cylinder $\cc = \cc^u[i_0, \ldots, i_m]$  of length $m$ in $U$ we have
\be\label{eq:Gibbs}
c_1 \leq  \frac{\nu(\cc)}{e^{g_m(y) }} \leq c_2 \quad , \quad y \in \cc ,
\ee
(see e.g. \cite{PP} or \cite{P2}).

We now state some basic properties of the metric $D_\theta$ that will be needed later.

\bs

\noindent
{\bf Lemma 5.1.}  
(a) {\it For any cylinder $\cc$ in $U$ the characteristic function $\chi_{\hcc}$ of
$\hcc$ on $\hU$ is Lipschitz with respect to $\dte$ and $\Lip_\theta(\chi_\cc) \leq 1/\diamte(\cc)$.}

\ms

(b) {\it There exists a constant $C_2 > 0$ such that if $x,y\in \hU_i$ for some $i$,  then 
$$|\tau(x) - \tau(y) |  \leq  C_2\, \dte(x,y) .$$
That is, $\tau\in \ff_\theta(\hU)$. Moreover, we can choose $C_2 > 0$ so that
$$|\tau_m(x) - \tau_m(y)| \leq C_2\,  \dte(\sigma^m(x),\sigma^m(y))$$
whenever $x,y\in \hU_i$ belong to the same cylinder $X$ of length $m$.}

\ms

(c) {\it There exist  constants $C_2 > 0$ and $\alpha_2 > 0$  such that for any $z\in R$, any cylinder 
$\cc$ in $W^u_R(z)$ and any $x,y\in \cc$ we have
$d(\tPsi(x),\tPsi(y))  \leq  C_2\, \dte(x,y)$ and
$\diamte(\cc) \leq C_2 (\diam(\tPsi(\cc)))^{\alpha_2}$. Moreover, we can take $\alpha_2 > 0$ so that
$1/(\gamma_1)^{\alpha_2} = \hat{\theta}$.}

\bs

\noindent
{\it Proof.} (a) Let $\cc$ be a cylinder in $U$ and let $x,y\in \hU$. 
If $x,y\in \cc$ or $x\notin \cc$ and $y\notin \cc$,
then $\chi_\cc(x) - \chi_\cc(y) = 0$. Assume that $x\in \cc$ and $y \notin \cc$.
Let $\dte(x,y) = \theta^{N+1}$ and let $\cc'$ be a cylinder of length $N$ containing both $x$
and $y$. Since $x\in \cc$, as well, and $x$ is an interior point of $\cc$, we must have
$\cc\subset \cc'$. Thus, $\diamte(\cc) \leq \dte(x,y)$. This gives
$$|\chi_\cc(x) - \chi_\cc(y)| = 1 = \frac{\diamte (\cc)}{\diamte(\cc)} \leq
\frac{1}{\diamte(\cc)}\, \dte(x,y) ,$$
which proves the assertion.

\ms

(b), (c)  Assume $x\neq y$ and let $\cc$ be the cylinder of largest length $m$
containing both $x$ and $y$. Set $\tx = \tPsi(x), \ty = \tPsi(y) \in \tR$.
Then $\dte(x,y) = \theta^{m+1}$. On the other hand, (2.1) and (5.1)  imply
\begin{eqnarray*}
|\tau(x) - \tau(y)| 
& \leq & |\tau|_{\alpha_1}\, (d(\tx,\ty))^{\alpha_1}  
\leq \frac{\Con}{(\gamma^{\alpha_1})^{m}} \leq  \Con \theta^m \leq C_2\, \dte(x,y)
\end{eqnarray*}
for some global constant $C_2 > 0$. The above also shows that 
$ d(\tx ,\ty ) \leq   \Con \theta^m \leq C_2\, \dte(x,y)$, which proves half of  part (c).
The second part of (c) follows by using a similar estimate and the other half of (2.1).

Next, assume that $x,y$ belong to the same cylinder $\cc$ of length $m$.
Let $\pp_j(x), \pp^j(y) \in R_{i_j}$ for all $j = 0,1, \ldots, m$.  Assume that 
$\dte(x',y') = \theta^{p+1}$, where
$x' = \sigma^m(x)$ and $y' = \sigma^m(y)$. Then $\dte (x,y) = \theta^{m+p+1}$
and moreover $\dte(\sigma^j(x), \sigma^j(y)) = \theta^{m-j+p+1}$ for all $j = 0,1, \ldots,m-1$. 
Then (2.1) and (5.1)  imply
\begin{eqnarray*}
|\tau(\sigma^j(x)) - \tau( \sigma^j(y))| 
& \leq & |\tau|_{\alpha_1} (d (\sigma^j(\tx), \sigma^j(\ty))^{\alpha_1} 
\leq  \Con\,  (d (\Psi(\sigma^j(x)), \Psi(\sigma^j(y)))^{\alpha_1} \\
& \leq & |\tau|_{\alpha_1} \left(\frac{1}{c_0 \gamma^{m-j+p}} 
d(\tpp^{m+p-j}(\Psi(\sigma^j x)), \tpp^{m+p-j}(\Psi(\sigma^j y)))  \right)^{\alpha_1 }\\
& \leq & \frac{\Con}{(\gamma^{\alpha_1 })^{(m-j+p)}} \leq \Con\, \theta^{m-j+p} 
\leq \Con \theta^{m-j+1} \dte(x',y') .
\end{eqnarray*}
So
\begin{eqnarray*}
|\tau_m(x) - \tau_m(y)| 
& \leq & \sum_{j=0}^{m-1}|\tau(\sigma^j(x)) - \tau( \sigma^j(y))|  
\leq \Con \dte(x',y')\, \sum_{j=0}^{m-1} \theta^{m-j+1}  \leq  \Con \, \dte(x',y') ,
\end{eqnarray*}
which proves the statement.
\endofproof

\bs


It follows from Lemma 5.1 that $\tau \in \ff_\theta(\hU)$, so assuming $f \in \ff_\theta(\hU)$, we have
$h_a\in \ff_\theta(\hU)$ for all $|a|\leq a'_0$. Then $\fa \in \ff_\theta(\hU)$ for all such $a$.
Moreover, using the analytical dependence of $h_a$ and $\lambda_a$ on $a$ and assuming that the 
constant $a'_0 > 0$ is sufficiently small, there exists $T = T(a_0')$ such that
\begin{equation}
T \geq \max \{ \, \|\fa \|_0 \, , \, |\fa_{|\hU}|_\theta \, , |\tau_{|\hU}|_\theta \, \}\;
\end{equation}
for all $|a| \leq a'_0$. Fix $a_0' > 0$ and $T > 0$ and with these properties.
Taking the constant $T > 0$ sufficiently large, we have  $|\fa - f^{(0)}| \leq T\, |a|$ 
on $\hU$ for $|a| \leq a'_0$.

The following Lasota-Yorke type inequality is  similar to that in \cite{D}, 
and in fact the same as the corresponding one
in \cite{St2} (although we now use a different metric) and its proof is also the same.

\bs

\noindent
{\bf Lemma 5.2.}  {\it There exists a constant $A_0 > 0$ 
such that for all $a\in \R$ with $|a|\leq a'_0$ the following holds:
If the functions $h$ and  $H$ on $\hU$, the constant $B > 0$  
and the integer $m \geq 1$ are such that $H > 0$ on $\hU$ and 
$|h(v) - h(v')| \leq B\, H(v')\, \dte (v,v')$ for any $i$ and any $v,v'\in \hU_i$, 
then for any $b\in \R$ with  $|b|\geq 1$ we have 
$$ |\lab^m h(u) - \lab^m h(u')| \leq  A_0\,\left[ B\,\theta^m \, (\ma^m H)(u') 
+ |b|\, (\ma^m |h| )(u')\right]\, \dte (u,u')$$
whenever $u,u'\in \hU_i$ for some $i = 1, \ldots,k_0$.} \endofproof

\bs

\noindent
{\bf Remark.} From the proof of this lemma (see e.g. the Appendix in \cite{St2}) that the constant
$A_0$ depends only on $\|f\|_\theta$ and some global constants, e.g. $c_0$ and $\gamma$ in (2.1).

\subsection{First step -- fixing $N$, a few compact sets of positive measure}

Let the constants $c_1$ and $c_2$ be as in (5.2).
{\bf Fix constants} $\rho_1 \in (0,1)$ and $q_0 \geq 1$ such that Lemma 4.1(a), (b), (c)
hold with $\rho' = \rho_1/8$ and $q' = q_0$. 


{\bf In what follows we will use the entire set-up and notation from Sect. 4}, 
e.g. the subsets $P_0$ and $\tP_0$ of $\ll_0\cap R$,
the numbers $r_0 \geq r_1 > 0$, $R_0 > 1$, etc.,  satisfying  (4.2), (4.3), etc.
Let $\eta_1(x), \eta_2(x), \ldots, \eta_{\ell_0}(x)$ ($x\in M$) 
be families of unit vectors in $E^u_1(x)$ as in the text just  before Lemma 4.4, and let $\ep'' \in (0,\ep')$, 
$0 < \delta'' < \delta' $, $\delta_0 > 0$ (depending on the  choice of $\kappa$), $\beta_1 \in (0,1)$,  $C_1 > 0$
 be constants with the properties described in Lemmas 4.4 and 4.3.  
Fix an integer  $N_0 \geq 1$ with the property described in Lemma 4.4(b). and then
{\bf fix an integer} $N \geq N_0$. A few additional conditions on $N$ will be imposed later.
Let $\ylo(Z) \in B^s (Z,\delta') \cap \pp^N(B^u(Z,\ep''))$, 
$\ylt(Z) \in B^s (Z,\delta') \cap \pp^N(B^u(Z,\ep''))$, ($Z\in P_o$; $\ell = 1, \ldots,\ell_0$) 
be families of points satisfying the requirements of Lemma  4.4. 

Assume the integer $n_0 \geq 1$ is chosen so large that for any $z\in R$ and any
unstable cylinder $\cc$ of length $\geq n_0$ in $R$ we have $\diam(\Psi(\cc)) \leq \ep''$
and $\diam(\T_z(\cc)) \leq \ep''$ for any $z\in \cc$.
Set 
\be \label{eq:hddefn}
\di \hd = \frac{\beta_0 \delta_0 \rho_1}{512 R^2_0} .
\ee

Let $E > 1$ be a constant -- we will see later how large it should be,  let
$\ep_1 > 0$ be a constant with
\be\label{eq:epone}
0 < \ep_1 \leq  \min\left\{ \;\frac{1}{32 C_0 }\;, \;  \frac{1}{4E} \, \right\} .
\ee
We will assume  $N > N_0$ is chosen so that
\be\label{eq:Ncond}
\gamma^N \geq \frac{1}{\delta''}   \quad , \quad 
\theta^N < \frac{\rho^2_1 \beta_0 \delta_0\ep_1}{256  E} 
\quad , \quad \theta_2^N < \frac{\hd \hrho \ep_1}{64 E} ,
\ee
where $\theta_2 = \max\{ \theta , 1/ \gamma^{\alpha_1 \beta_1}\}$, 
$\beta_1 > 0$ being the constant from Lemma 4.3. 

\bs

\noindent
{\bf Lemma 5.3.} {\it Let $\cc$ be an unstable cylinder in $R$ of length 
$m \geq 1$  with $\pp^m(\cc) \cap P_0 \neq \e$. }

\ms

(a) {\it There exist sub-cylinders $\dd$ and $\dd'$ of $\cc$ of co-length $q_0$ such that 
$d(\tPsi(y),\tPsi(x)) \geq \frac{\rho_1}{2}  \, \diam(\tPsi(\cc))$ for all $y \in \dd'$ and $x\in \dd$.
Moreover, we can take one of the sub-cylinders, e.g. $\dd$, so that it contains $z$.} 

\ms

(b)   {\it  There exists an integer $q_1 \geq q_0$ such that for any  sub-cylinder 
$\cc'$ of $\cc$ of co-length $q_1$ with $\pp^m(\cc')\cap P_0 \neq \e$ 
we have} $\diam(\tPsi(\cc')) \leq \min\left\{ \frac{\rho_1}{8} , \frac{\hd }{8 C_1} \right\}\, \diam (\tPsi(\cc)) .$

\bs

\noindent
{\it Proofs.} (a) Take $z,x  \in \cc$ as in Lemma 4.1(c), and let $\dd$ and $\dd'$ be
the sub-cylinders of $\cc$ of co-length $q_0$ containing $z$ and $x$, respectively. By Lemma 4.1
and the choice of $q_0$ it follows that $\diam(\tPsi(\dd)) \leq \frac{\rho_1}{8} \diam(\tPsi(\cc))$. 

Next, by the choice of $z,x$ in Lemma 4.1(c), for any $y \in \dd'$ we have 
$d(\tPsi(y),\tPsi(z)) \geq \rho_1\, \diam(\tPsi(\cc))$. 
Then
$$d(\tPsi(x),\tPsi(y)) \geq d(\tPsi(y), \tPsi(z)) - d(\tPsi(x), \tPsi(z))  
\geq \rho_1 \diam(\tPsi(\cc)) - \frac{\rho_1}{8} \diam (\tPsi(\cc))  > \frac{\rho_1}{2}\,  \diam (\tPsi(\cc))  $$
for any $y\in \dd'$ and any $x\in \dd$.

\ms

(b) This follows from Lemma 4.1(b): take $q_1 = q_0^r$ for some sufficiently large integer $r \geq 1$.
\endofproof


\subsection{Main consequence of Lemma 4.3}

We will use Lemma 4.4 with $\kappa = \hrho/2$, where $\hrho = \frac{\rho_1}{8C_0^2}$, $\rho_1$ being
the constant from Lemma 4.1.
As before set $\di \hd = \frac{\beta_0 \delta_0 \rho_1}{512 R^2_0}$,
where  $\beta_0 > 0$ and $\delta_0 > 0$  are fixed constants with the properties described in Lemma 4.4.
We will also use the integers $N_0  \geq 1$ and the constants $\ep'' > 0$ and $\delta' > \delta'' > 0$ from Lemma 4.4.
Assume $L > 2$ is a large constant (to be determined later) and for any integer $m \geq L$, let the set
$\Xim$ be defined as in Sect. 4.3.

\bs

\noindent
{\bf Lemma 5.4.} 
{\it For any  $m \geq L$, any point $Z \in P_0\cap \Xim$, any integer $N \geq N_0$, any 
$\ell = 1, \ldots, \ell_0$ and any $i = 1,2$ 
there exists a (H\"older) continuous map $B^u(Z,\ep'') \ni x \mapsto v^{(\ell)}_i (Z, x) \in U$ 
such that $\sigma^N(\vl_i(Z, x)) = x$  for all $x\in B^u(Z,\ep'')$
and the following property holds:

For  any cylinder $\cc$ in $B^u(Z,\ep'')$ of length $m$ with $Z\in \cc$ 
there exist sub-cylinders $\dd$ and $\dd'$ of $\cc$ of co-length $q_1$ 
and  $\ell = 1, \ldots, \ell_0$ such that $Z \in \dd$ and for any points $x\in \dd$ 
and $z\in \dd'$, setting $x' = \piU(x)$, $z' = \piU(z)$, 
we have $d(\T_Z(x), \T_Z(z)) \geq \frac{\hrho}{2} \, \diam(\T_Z(\cc))$ and 
$$I_{N, \ell}(x',z') =  \left|  \varphi_{\ell,N}(Z,x') - \varphi_{\ell,N}(Z,z') \right|  
\geq \hd \, \diam(\tPsi(\cc)) ,$$
where 
$\varphi_{\ell.N}(Z, x) = \tau_{N}(\vl_1(Z,x)) - \tau_{N}(\vl_2(Z,x)) .$
Moreover, $I_{N,\ell}(x',z') \leq C_1\, \diam(\tPsi(\cc))$ for any $x,z \in \cc$, 
where $C_1 > 0$ is the constant from Lemma} 4.3.

\bs

\noindent
{\it Proof.}  Fix for a moment  $m \geq L$, $Z \in P_0 \cap \Xim$, $N \geq N_0$ and $\ell = 1, \ldots, \ell_0$. 
Using  Lemma 4.4, there exist points  $\ylo = \ylo(Z), \ylt = \ylt(Z) \in W^s_{\delta'} (Z)$ such that 
the property (b) in Lemma 4.4 holds.

Thus, given $i = 1,2$, there exists  a cylinder $L^{(\ell)}_i = L^{(\ell)}_i(Z)$ of length $N$ in $B^u(Z,\ep'')$ 
so that 
$$\pp^{N} : L^{(\ell)}_i \cap R \longrightarrow W^u_{R_{i_0}}(\tyli)\cap R$$
is a bijection; then it is a bi-H\"older homeomorphism. Consider its inverse and its H\"older continuous
extension  $\pp^{-N} : W^u_{R_1}(\yli)   \longrightarrow L^{(\ell)}_i$.

Set $M^{(\ell)}_i = M^{(\ell)}_i(Z) = \piU(L^{(\ell)}_i(Z)) \subset U$; this is then a cylinder 
of length $N$ in $U$. Define the maps
$$\tvli (Z, \cdot) : U_{i_0} \longrightarrow L^{(\ell)}_i \subset B^u(Z,\ep'')  
\quad , \quad \vli (Z, \cdot)   : U_{i_0} \longrightarrow M^{(\ell)}_i  \subset U$$
by 
$$\tvli(Z, y) = \pp^{- N} (\pi_{\yli}(y)) \quad, \quad \vli(Z, y) = \piU(\tvli(Z, y)) .$$
Then  
\be
\pp^{N}(\tvli(Z, y)) = \pi_{\yli}(y) = W^s_{\ep_0}(y) \cap W^u_{R_{i_0}}(\yli) ,
\ee
and
\be
\pp^{N}(\vli(Z, y)) = W^s_{\ep_0}(y) \cap \pp^{N}(M^{(\ell)}_i) = \pi_{\dli}(y) ,
\ee
where $\dli = \dli(Z) \in W^s_{R}(Z)$ is such that $\pp^{N}(M^{(\ell)}_i) = W^u_{R}(\dli)$.
Next, there exist $x' \in M^{(\ell)}_i$ and $y\in L^{\ell}_i$ with $\pp^{N} (x') = \dli$ and
$\pp^{N}(y') = \yli$. Since stable leaves shrink exponentially fast, using (2.1) and (\ref{eq:Ncond}) we get
\be
d(\dli , \yli) \leq \frac{1}{c_0 \gamma^{N}} d(x',y') \leq  \frac{1}{\gamma^{N}} < \delta'' .
\ee
{\bf Thus, $\dlo, \dlt$ satisfy the assumptions and therefore the conclusions of Lemma 4.4(b).}

Let $\cc$ be a cylinder in $B^u(Z,\ep'')$ of length $m$ with $Z \in \cc \cap P_0 \cap \Xim$. 
Set $\cc' = \T_Z(\cc)$ and $\tcc = \tPsi(\cc)$. By the choice of
the constant $C_0$, we have $\frac{1}{C_0} \diam(\tcc) \leq \diam(\cc') \leq C_0\, \diam (\tcc)$.
Let $\dd$ be the sub-cylinder of $\cc$ of co-length $q_1$. 

Next, by Lemma 5.3(a), there exist $\tx \in \cc$ such that if $\dd'$ is the sub-cylinder of $\cc$ of co-length 
$q_1$ containing $\tx$, then $d(\tPsi(y), \tPsi(x)) \geq \frac{\rho_1}{2} \diam(\tcc)$ for all
$y \in \dd$ and $x\in \dd'$. Thus, 
$d(\T_Z(y), \T_z(x)) \geq \frac{\rho_1}{2C_0^2} \diam(\cc') \geq \hrho \, \diam(\cc')$ for all
$y \in \dd$ and $x\in \dd'$.

Set $x_0 = \T_Z(\tx)$ and let $x_0 = \Phi^u_{\tZ}(u_0)$, where $u_0 \in E^u(Z)$.
By the choice of the constant $\beta_0$ and the family of unit vectors 
$\{\eta_\ell(Z)\}_{\ell=1}^{m_0}$, there exists some $\ell = 1, \ldots, \ell_0$ such that 
$ \left\langle \frac{u_0}{\|u_0\|}, \eta_\ell(Z) \right\rangle \geq \beta_0 .$
Moreover, $d(x_0, Z) \geq \hrho\, \diam(\T_Z(\cc))$.
It then follows from Lemma 4.4(b) with $\kappa = \hrho/2$ and (5.9) that 
\be
\frac{\beta_0 \delta_0 \hrho}{32 R^2_0 } \, \diam(\T_z(\cc)) \leq \left| \Delta(x_0, \dlo) - \Delta(x_0, \dlt)\right| .
\ee
(In the present situation, since $\dlo , \dlt \in W^s_R(Z)$, we have $\pi_{\dli}(Z) = \dli$
for $i = 1,2$.)

Consider the projections of $\tx, Z$ to $U$ along stable leaves:
$x' = \piU(\tx) \in U_i$, $ Z' = \piU(Z) \in U_i$, where $R_i$ is the rectangle containing
$Z$ (and therefore $\cc$).  We have
\begin{eqnarray*}
I_{N, \ell}(x',Z')
& = & \left| \; [\tau_{N}(\vl_1(Z,x')) - \tau_{N}(\vl_2(Z,x'))]
 - [\tau_{N}(\vl_1(Z,Z')) - \tau_{N}(\vl_2(Z,Z'))] \; \right|\\
& = & \left| \; [\tau_{N}(\vl_1(Z, x')) - \tau_{N}(\vl_1(Z,Z'))] 
- [\tau_{N}(\vl_2(Z,x')) - \tau_{N}(\vl_2(Z,Z'))] \; \right|\\
& = & \left| \, \Delta(\pp^{N} (\vl_1(Z, x')), \pp^{N}(\vl_1(Z, Z'))) 
-  \Delta(\pp^{N} (\vl_2(Z, x')), \pp^{N}(\vl_2(Z, Z'))) \,\right|\\
& = & \left| \, \Delta(\pi_{\dlo}(x'), \pi_{\dlo}(Z') ) -  \Delta(\pi_{\dlt}(x') , \pi_{\dlt} (Z') ) \,\right|\\
& = & \left| \, \Delta( x' , \dlo ) -  \Delta( x', \dlt ) \,\right| 
= \left| \, \Delta( \tx , \dlo ) -  \Delta( \tx  , \dlt ) \,\right|.
\end{eqnarray*}
We claim that the latter is the same as the right-hand-side of (5.10). Indeed, let
$\Delta(\tx, \dlo) = s_1$ and $\Delta(\tx, \dlt) = s_2$. 
Then $\phi_{s_1}([\tx, \dlo]) \in W^u_{\ep_0}(\dlo)$ and $\phi_{s_2}([\tx, \dlt]) \in W^u_{\ep_0}(\dlt)$. 
Let $\phi_s(x_0) = \tx$. It is then straightforward to see that
$\Delta(x_0, \dlo) = s + s_1$ and $\Delta(x_0, \dlt) = s+ s_2$. Thus,
$$\left| \Delta(x_0, \dlo) - \Delta(x_0, \dlt)\right| = |(s+s_1) - (s+s_2)| = |s_1-s_2|
= \left| \Delta(\tx, \dlo) - \Delta(\tx, \dlt)\right| .$$
Combining this with (5.10) gives
$I_{N,\ell}(x' ,Z') \geq \frac{\beta_0 \delta_0 \hrho}{8R_0} \, \diam(\T_z(\cc)) \geq 2\hd \, \diam(\tPsi(\cc)) .$

For arbitrary $x,z\in \cc$, setting $x' = \piU(x)$, $z' = \piU(z)$,
the above calculation and Lemma 4.3 give
$ I_{N,\ell}(x' ,z') = \left| \Delta(x, \pi_{\dlo}(z)) 
- \Delta(x, \pi_{\dlt}(z))\right| \leq C_1\, \diam(\tPsi(\cc)) .$
The same argument shows that for any $z\in \dd$, using Lemma 5.3(b) the fact that $Z \in \dd$,  we have
\begin{eqnarray*}
I_{N,\ell}(z',Z') =  \left| \Delta(z, \dlo) 
- \Delta(z, \dlt)\right| \leq C_1\, \diam(\tPsi(\dd))  \leq \frac{\hd}{8} \,\diam (\tPsi(\cc)) .
\end{eqnarray*}
Similarly, for any $z\in \dd'$ and $z_0 = \T_Z(z)$ we have
$$I_{N,\ell}(z',Z') = \left| \Delta(z_0, \dlo) - \Delta(z_0, \dlt)\right| \leq C_1\, \diam(\tPsi(\dd')) 
\leq \frac{\hd}{8} \,\diam (\tPsi(\cc)) .$$
Since 
$\Delta(x,\pi_{\dli}(z)) = \Delta(x, \dli) - \Delta(z, \dli) ,$
it follows that
\begin{eqnarray*}
I_{N,\ell}(x' ,z') \geq   I_{N,\ell}(x' , Z') -  I_{N,\ell}(z' , Z') 
\geq  2\hd \, \diam(\tPsi(\cc)) - \hd \,\diam (\tPsi(\cc)) = \hd\, \diam(\tPsi(\cc)) .
\end{eqnarray*}
This completes the proof of the lemma.
\endofproof


\def\hB{\widehat{B}}


\section{Contraction operators}
\setcounter{equation}{0}

We use the notation and the set-up from Sect. 5. 

\subsection{Choice of cylinders, definition of the contraction operators}

{\bf Fix a constant} $A > 0$ so large that for any $i =1, \ldots, k_0$ and any $x,x\rq{}\in R_i$ we have 
$\diam(\tPsi(W^u_{R_i}(x))) \leq A\, \diam(\tPsi(W^u_{R_i}(x')))$.

{\bf Fix integers} $d \geq 1$ and $t_0$ such that
\be\label{eq:t0condition}
c_0 \gamma^d > \frac{1}{\hrho} \quad, \quad 
t_0 \geq \frac{1}{\beta_1 \log \gamma} \left|\log  \frac{32C_4 r_0^2}{\beta_0 \delta_0 \hrho^{q_1}}\right| + \left|\frac{\log c_0}{\log \gamma} \right| ,
\ee
where $q_1 \geq 1$ is the constant from Lemma 5.3(b), while 
$\gamma > 1$ is the constant from the end of Sect. 2. {\bf Fix an integer} $N \geq N_0$ as in Sect. 5.
Assume that the constant $\hd_0 > 0$ that appears in $\Xi_n(p_0,\hep_0,\hd_0)$ (see Sect. 3) is chosen so small that
\be\label{eq:dcond}
1-\hd_0 < \frac{d+1}{d+2} \quad, \quad 
\delta = (d+3) \hd_0 < \frac{1}{2}.
\ee
Set
\be\label{eq:mucond}
\mu_0 = \mu_0(N) = \min \left\{\;
\frac{\theta^{2N + 2d} }{6 e^{T/(1-\theta)}}\; , \; 
\frac{1}{8\,e^{2 T N}} \, \sin^2\left(\frac{\hd \hrho \ep_1}{16}\right) \; , \; \frac{1}{40}\; \right\} ,
\ee
and
\be\label{eq:b0cond}
b_0 = b_0(N) = \max\left\{ \theta^{-N}\,  , \, \left(\frac{2 C_0 \gamma_1^{d}}{c_0 \hd} \right)^{1/\alpha_1} 
\, , \, \left(\frac{3C_2 T e^{T/(1-\theta)}}{(1-\theta) } \right)^{1/\alpha_2}\right\} ,
\ee
where $\alpha_2 > 0$ is as in Lemma 5.1(c).

It follows easily from  Lemmas 4.1 and  5.3 that for any maximal unstable cylinder $\cc$ in $R$ with
$\diam(\tPsi(\cc)) \leq \frac{\ep_1}{|b|} $ we have
$$\hrho\, \frac{\ep_1}{|b|} \leq \diam(\tPsi(\cc)) \leq \frac{\ep_1}{|b|}  ,$$
and by Lemma 5.1, if $\ell$ is the length of $\cc$, then
$$\frac{-\log C_2 - \beta' \log (\hrho \ep_1)}{|\log \theta|} + \frac{\beta'}{|\log\theta|}\, \log |b| \leq
\ell \leq \frac{\log C_2 - \log \ep_1}{|\log \theta|} + \frac{1}{|\log\theta|}\, \log |b| .$$
Thus, there exists a global constant $B > 1$ (independent of $b$) such that if $|b| \geq b_0$, then
\be\label{eq:Bcond}
\frac{1}{B} \log |b| \leq \ell \leq B\, \log |b|  .
\ee
{\bf Fix a constant } $B > 1$ with this property. Later we may have to impose some further
requirements on $B$.  Assume that the constant $L > 3/\tau_0$ from Sect. 4.3 is such that
$$L > \hB = \frac{B}{\hd_0} .$$
Some further conditions on $L$ will be imposed later.

Throughout the rest of Sect. 6, $b$ will be a {\bf fixed real number}  with $|b| \geq b_0$. 
Set 
\be\label{eq:hbdef}
\hb = \lceil \log |b| \rceil  .
\ee

For every $z\in \ll$ denote by $\cc(z)$ {\it the maximal cylinder} in $W^u_{R}(z)$ with
$\diam (\tPsi(\cc(z))\leq \ep_1/|b|$. 
If $\pp^k(\cc(z)) \cap P_0 \neq \e$, where $k$ is the length of $\cc(z)$,
then the maximality and Lemma 4.1(a) imply $\diam (\tPsi(\cc(z))\geq \hrho \,\ep_1/|b|$.

We will now define a subset $P_1 = P_1 (b)$ of  $\ll$ as follows. Given $x\in P_0$, let 
$y(x) = \pp^{-m_x}(x)$ for some integer $m_x \geq 0$ be such that the cylinder
$\cc_{m_x}(y(x))$ of length $m_x$ satisfies
$$\diam(\tPsi(\cc_{m_x}(y(x)))) \leq \frac{\ep_1}{|b|}$$
and $m_x$ is the smallest number with this property. Thus,
$\diam(\tPsi(\pp(\cc_{m_x}(y(x))))) > \frac{\ep_1}{|b|} .$
Set
$$P_1 = \{ y(x) : x \in P_0 \} .$$
and
\be\label{eq:Lambda}
\Lambda_N (b) = \left\{ x\in \ll \cap R : \:\: \sharp \, \left\{ j : 0 \leq j <  L \, N\, \hb \; , \;  
\pp^j(x) \notin P_1 \cap P_0 \cap \Xib_{B}\right\}  \geq \frac{\delta}{N}\,  L\, N \, \hb \: \: \right\}.
\ee

For later use we have to fix some constants.
Having fixed $\theta \in  [\hat{\theta}, 1)$ earlier, now {\bf fix a constant}  $\theta_1 \in (0,\theta)$ with 
$$0 < \theta_1 \leq \frac{1}{\gamma_1} \quad , \quad (\theta_1)^{\alpha_2/2} \leq \theta .$$
Set  
$\theta_2 = \max \{ \theta, 1/ \gamma^{\alpha_1 \beta_1}\} ,$
where $\beta_1 \in (0,1)$ is the constant from Lemma 4.2.
There exists $\beta_2 \in (0,1)$ with $\theta = \theta_1^{\beta_2}$. Set $s = \frac{2}{\beta_2} > 0$.
We will {\bf assume that} the constant $L > 0$ is chosen so large that
\be
c_3 = c_0 L/B > 16s = \frac{32}{\beta_2} ,
\ee
where $c_0 > 0$ and $C_0 > 0$ are the constant from (\ref{eq:XinB}).

\bs

\noindent
{\bf Lemma 6.1.} {\it Assuming $b_0 > 0$ is chosen sufficiently large, we have
$\mu (\Lambda_N (b)) \leq \frac{2C_0}{|b|^{ c_3}}$
for all $b$ with $|b| \geq b_0$.}

\bs

\noindent
{\it Proof.} 
Set $m = \hb$.
We claim that $\Lambda_N (b) \subset \Xi_{L m} \cup Y$, where
$Y = \cup_{j=0}^{L m} \pp^{-j}(R\setminus \Xib_{B}) .$

 Assume that there exists
$x \in \Lambda_N(b) \setminus (\Xi_{L m} \cup Y)$. Then $x \notin \Xi_{Lm}$, so  $\pp^j(x) \in P_0$ for more than 
$(1- \hd_0) L m $ values of $j = 0,1, \ldots, L m -1$. 
Notice that by (\ref{eq:dcond}), $\frac{\delta}{N}\,  L\, N \, \hb = (d+3) \hd_0 L m$. 
Now $x \in \Lambda_N (b)$ and (\ref{eq:Lambda}) imply $\pp^j(x) \notin P_1 \cap \Xib_{B}$ for at least $(d+2) \hd_0 L m$ values of 
$j = 0,1, \ldots, L m -1$.
However $x \notin Y$ shows that $\pp^j(x) \notin  \Xib_{B}$ for all $j = 0,1, \ldots, L m - 1$, so we must have $x \in Z$, where
$$Z = \{ x \in \ll \cap R : \pp^j(x) \notin P_1 \: \mbox{\rm for at least 
$(d+2) \hd_0 L m$ values of}\:  j = 0,1, \ldots, L m -1 \}.$$
We will show that $Z\setminus \Xi_{Lm} = \e$; this will lead to a contradiction.  

Assume that $x \in Z \setminus \Xi_{Lm}$.
Given $y = \pp^j(x)$ for some $j \leq (1-\hd_0) L m$, let $k_j$ be the length of the maximal cylinder $\cc(y)$ in $W^u_{R}(y)$ with
$\diam (\tPsi(\cc(y))\leq \ep_1/|b|$. Then by (\ref{eq:Bcond}), $k_j \leq B \hb = Bm$. Moreover, if $y \notin P_1$, 
then $\pp^{k_j}(y) \notin P_0$, i.e. $\pp^{k_j + j} (x) \notin P_0$. On the other hand, by the choice of the constant $d$ we have 
$c_0\gamma^d \hrho > 1$. Since  $j \leq  (1-\hd_0) L m$, we have $\pp^i(y) \in P_0$ for some $i > 0$, and now Lemma 4.1(a) gives
$\diam (\tPsi(\cc(y))\geq \hrho \ep_1/|b|$, while (2.1) implies 
$\diam (\tPsi(\pp^d\cc(y))) \geq c_0 \gamma^d \hrho \ep_1/|b| > \ep_1/|b|$.
Thus, $\pp^d(\cc(\pp^j(x)) )$ contains  $\cc (\pp^{j+d}(x))$ as a proper subset, so the length of 
$\pp^d(\cc(\pp^j(x)) )$ is strictly less than
the length of $\cc(\pp^{j+d}(x))$, i.e. $k_j - d <  k_{j+d}$. This yields $k_j + j < k_{j+d} + (j+d)$ for all $j \leq (d+1) B m$ with
$\pp^j(x) \notin P_1$. Notice that for $j \leq (1-\hd_0) L m$ we have 
$k_j + j \leq Bm + (1-\hd_0)  L  m  < \delta_0 L m + (1-\hd_0) L m  = L m$.

There are at least $(d+2) \hd_0 L m$ values of $j = 0,1, \ldots, L m -1$ with $\pp^j(x) \notin P_1$, so
at least $(d+1) \hd_0 L m $ of them satisfy $j \leq (1-\hd_0) L m$. For such $j$, the sequence $\{ k_j + j\}$ contains a strictly increasing
subsequence with at least $(d+1) \hd_0 L m/ d $  members. Since $(d+1) \hd_0 L m/ d > \hd_0 L m$, we have
$\pp^{k_j+j}(x) \notin P_0$ for at least  $\hd_0 L m$ values of $j \leq (1-\hd_0) L m$, which is a contradiction with
$x\notin \Xi_{Lm}$.

This proves that $\Lambda_N (b) \subset \Xi_{L m} \cup Y$, and now (\ref{eq:tailestimate}) and (\ref{eq:XinB}) give
$$\mu(\Lambda_N (b)) \leq C' e^{-c' L \hb} + C_0 L\, \hb \, e^{-c_0 L\hb/B}  
\leq 2C_0 e^{-(c_0L/B) \log |b|} = \frac{2C_0}{|b|^{c_0L/B}} ,$$
assuming $|b| \geq b_0$ and $b_0$ is chosen sufficiently large. This proves the lemma.
\endofproof

\bs

\noindent
{\bf Definitions 6.2} (Choice of cylinders) Here we define an important family of cylinders in 
$R$ and $U$ and some sub-cylinders of
theirs that will play an important role throughout Sections 6 and 7.

Set 
$$K_0 = \piU(P_1 \cap P_0 \cap \Xib_B) .$$ 
For any $u \in K_0$ amongst the cylinders 
$\cc(z)$ with $z\in P_1\cap P_0 \cap \Xib_B$ and $\piU(z) = u$, there 
is {\bf one of maximal length}. Choose one of these -- it has the form $C(Z(u))$ for some 
$Z(u) \in P_1 \cap P_0\cap \Xib_B$ with $\piU(Z(u)) = u$. Set $\cc'(u) = \piU(C(Z(u)))$. 
It follows from this choice that for any $z\in R$ with $\piU(z) = u$ we have $\cc'(u) \subset \piU(\cc(z))$.

Since the lengths of the cylinders $C'(u)$ are bounded above and 
$K_0 \subset \cup_{u\in K_0} \cc'(u)$, there exists 
finitely many different cylinders $\cc'_m = \cc'(u_m)$ for some $m =1, \ldots, m_0$ such that
$$K_0 \subset \cup_{m=1}^{m_0} \cc'_m .$$
Different cylinders have no common interior points, so
$\cc'_m \cap \cc'_{m'} \cap \hU = \e$ for $m \neq m'$. For each $m$, set $\cc_m = \cc(Z(u_m))$; then
$Z_m = Z(u_m) \in \cc_m \cap P_1 \cap P_0 \cap \Xib_B$ is so that $\piU(Z_m) = u_m$. 
According to the definitions of the cylinders $\cc(z)$, $\cc_m = \cc(Z_m)$ is a 
{\it maximal  closed cylinder in $W^u_R(Z_m)$ with} $\diam(\tPsi(\cc_m)) \leq \ep_1/|b|$.
Let $\dd_1, \ldots, \dd_{\tilde{j}}$ be the list of
{\it all closed unstable cylinders  in $R$ which are sub-cylinders of co-length} $q_1$ of some $\cc_m$. 
Here $q_1 \geq 1$ is the constant from Lemma 5.3(b). Set $\dd_j' = \piU(\dd_j) \subset U$.
Re-numbering the cylinders  $\dd_j$ if necessary, we may assume there exists $j_0 \leq \tilde{j}$ such that
and $\dd_1, \ldots, \dd_{j_0}$ is {\it the list of all sub-cylinders
$\dd_j$ } such that $\dd_j \cap P_1 \cap P_0 \cap \Xib_B \neq \e$.

\bs

From the choice of the cylinders $\cc_m$  and Lemmas 4.1 and  5.3 we get:
\be
\hrho\, \frac{\ep_1}{|b|} \leq \diam(\tPsi(\cc_m)) \leq \frac{\ep_1}{|b|}  \quad , \quad 1 \leq m \leq m_0 .
\ee
If $\ell_m$ is {\it the length of the cylinder} $\cc_m$, it follows from  (6\ref{eq:Bcond}) that
\be
\frac{1}{B} \log |b| \leq \ell_m \leq B\, \log |b| \quad , \quad m = 1, \ldots, m_0 .
\ee

Set 
\be
V_b = \cup_{j=1}^{j_0} \dd'_j  \subset U .
\ee
It follows from the construction that
$K_0 \subset V_b  .$

\def\tnn{\widetilde{\nn}}

We are now ready to define an important family of {\it contraction operators}.
For any $\ell = 1, \ldots, \ell_0$, $i = 1,2$ and $j = 1, \ldots, j_0$, consider the unique 
$m = 1, \ldots, m_0$ with $\dd_j \subset \cc_m$, and set 
$$\vli = \vli(Z_m, \cdot) \quad , \quad \xijl = \vli(\dd'_j) \subset \hU .$$
where $\vli(Z_m, \cdot)$ is the map from Lemma 5.4 for the integer $N$. We will consider this map only on $\cc'_m$.
By Lemma 5.1(a), the {\it characteristic function} 
$\omijl = \chi_{\xijl} : \hU \longrightarrow [0,1]$
of $\xijl$ belongs to $\ff_\theta(\hU)$ and  $\Lip_\theta(\omijl) \leq 1/\diamte(\xijl)$. 



\bs

A subset $J$ of the set
$\Pi(b)  = \{\; (i, j , \ell) \; :  \; 1\leq i \leq 2\; ,\;  1\leq j\leq j_0 \;, 
\; 1 \leq \ell \leq \ell_0 \;\} $
will be called {\it representative} if for every $j = 1, \ldots, j_0$ there exists 
at most one pair $(i,\ell)$ such that
$(i,j,\ell) \in J$, and for any $m = 1, \ldots, m_0$ there exists $(i,j, \ell) \in J$ 
such that $\dd_j \subset \cc_m$.
Let $\jj(b)$ be the {\it family of all representative subsets} $J$  of $\Pi(b)$.

Given $J  \in \jj(b)$, define the function   
$\omega = \omega_{J}(b)  : \hU \longrightarrow [0,1]$ by
$$\di \omega  = 1- \mu_0 \,\sum_{(i, j,\ell) \in J} \omijl .$$
Clearly $\omega \in \ff_\theta(\hU)$ and $\frac{1}{2} \leq 1-\mu_0 \leq \omega (u) \leq 1$ for any $u \in \hU$. 
Define the {\it contraction operator}: 
$$\nn = \nn_J(a,b) : \ff_\theta (\hU) \longrightarrow \ff_\theta (\hU) \quad 
\mbox{\rm by } \quad \nn h = \ma^{N} (\omega_J \cdot h) .$$

\subsection{Main properties of the contraction operators}

First, we  derive an important consequence of the construction of the cylinders $\cc_m$ and $\dd_j$.

\bs

\noindent
{\bf Lemma 6.3.} {\it If $\sigma^p(\dd'_j) \subset \cc'_{k}$  for some $p\geq 0$, $j \leq j_0$ and  
$k \leq m_0$, then $p \leq t_0$, where $t_0$ is the integer given by (6.1). Moreover, there exists
a constant $s_0 \geq 0$ independent of $b$  such that the co-length of $\sigma^p(\dd'_j)$ in $\cc'_k$ 
does not exceed $s_0$.}

\bs

\noindent
{\it Proof.} Assume  $\dd = \sigma^p(\dd'_j) \subset \cc'_{k}$  for some $p > 0$, $j \leq j_0$ and  
$k \leq m_0$. 
From the assumptions we get $\pi_{Z_k}(\dd) \subset \cc_{k} \subset W^u_R(Z_k)$.

According to the choice of the sub-cylinders $\dd_j$, there exists $Z \in \dd_j \cap P_1 \cap P_0\cap \Xib_B$.
Then using Lemma 4.4(a) with $\kappa = 1/2$ and an appropriately chosen 
$X \in \chU(Z) \subset W^u_{\ep}(Z)$ with 
\be
d(X, Z) \geq \frac{1}{2} \diam(\T_z(\dd_j)) \geq \frac{\hrho^{q_1}\ep_1}{2|b|} ,
\ee
we can find points $d_1,d_2 \in W^s_R(Z)$ such that
$$|\Delta(X, \pi_{d_1}(Z)) - \Delta (X, \pi_{d_2}(Z))|
\geq \frac{\beta_0 \delta_0}{32R_0^2}\, \diam (\T_Z(\dd_j))$$
that is
\be 
|\Delta(X, d_1) - \Delta (X, d_2)|
\geq \frac{\beta_0 \delta_0}{32R_0^2}\, \frac{\hrho^{q_1}\ep_1}{|b|} .
\ee

Let $T = \tau_p(Z)$; then $z = \pp^p(Z) = \phi_T(Z)$. 
Next, consider the points $d'_i = \phi_T(d_i) \in W^s_R(z)$ ($i = 1,2$) and
$x = \phi_T(X) \in W^u_{\ep}(z) \subset \phi_T(W^u_{\ep}(X))$. 
It follows from (6.12) and the properties of temporal distance that
\be
|\Delta(x, d'_1) - \Delta (x, d'_2)|
\geq \frac{\beta_0 \delta_0}{32R_0^2}\, \frac{\hrho^{q_1}\ep_1}{|b|} ,
\ee
while (2.1) yields
\be
d(d'_1, d'_2) = d(\pp^p(d_1), \pp^p(d_2)) \leq \frac{1}{c_0 \gamma^p}\, d(d_1,d_2)
\leq \frac{1}{c_0 \gamma^p} .
\ee
Let $y \in W^s_R(z)$ be such that $\cc_m \subset W^u_R(y)$. Since $\pp^p(\dd'_j) \subset \cc'_m$,
we have $y \in \cc_m$. Using this again, for the point $x' = \pi_y(x) \in W^{sc}_{\ep}(y)$ we have
$\phi_t(x') \in \cc_m$ for some $t\in \R$, so $x'' = \T_y(\phi_t(x')) \in \T_y(\cc_m)$.
Moreover it is easy to see, using just the definition of the temporal distance function
and the fact that $x'' = \phi_s(x')$ for some $s\in \R$, that
$$|\Delta(x, d'_1) - \Delta(x, d'_2)| = |\Delta(x', d'_1) - \Delta(x', d'_2)|
= |\Delta(x'', d'_1) - \Delta(x'', d'_2)| .$$
This and (6.13) give
$$|\Delta(x'', d'_1) - \Delta(x'', d'_2)| 
\geq \frac{\beta_0 \delta_0 \hrho^{q_1}}{32R_0^2}\, \frac{\ep_1}{|b|} .$$
Combining the latter with $\diam(\T_y(\cc_m)) \leq \ep_1/|b|$, $y, x'' \in \T_y(\cc_m)$,
$d'_1, d'_2 \in W^s_R(y)$,
since $\cc_m \cap P_1 \cap \Xib_B \neq \e$, Lemma 4.3(b) implies
$$\frac{\beta_0 \delta_0 \hrho^{q_1}}{32R_0^2}\, \frac{\ep_1}{|b|} \leq 
|\Delta(x'', d'_1) - \Delta(x'', d'_2)| \leq C_1 \diam(\T_y(\cc_m)) \, (d(d'_1, d'_2))^{\beta_1}
\leq C_1 \frac{\ep_1}{|b|} \, (d(d'_1, d'_2))^{\beta_1} .$$
This and (6.14) yield
$\left(\frac{\beta_0 \delta_0 \hrho^{q_1}}{32C_1 R_0^2}\right)^{1/\beta_1} \leq 
\frac{1}{c_0 \gamma^p} ,$
so $p \leq t_0$, where $t_0 > 0$ is the integer from (6.1).

Next, let $s$ be the co-length of $\sigma^p(\dd'_j)$ in $\cc'_m$. Denote by $Q$ the cylinder
in $W^u_R(z)$ such that $Q \parallel \cc_m$. Then $\pp^p(\dd_j)$ is a sub-cylinder of $Q$
of co-length $s$, so $\dd_j$ is a sub-cylinder of co-length $s$ of $Q' = \pp^{-p}(Q)$. 
Since $Z \in \dd_j \cap P_1 \cap P_0 \cap \Xib_B$, it follows from Lemma 4.4(a) with $\kappa = 1/2$
that there exist $x_0 \in \T_Z(Q')$ and $y_1,y_2 \in W^s_R(Z)$ such that
$$\frac{\beta_0\delta_0}{32 R_0^2} \, \diam(\T_Z(Q')) \leq |\Delta(x_0, y_1) - \Delta(x_0,y_2)| .$$
Setting $x'_0 = \phi_T(x_0) \in \T_z(Q)$ and $y'_i = \phi_T(y_i) \in W^s_R(z)$, $i = 1,2$,
we have $|\Delta(x_0, y_1) - \Delta(x_0,y_2)| = |\Delta(x'_0, y'_1) - \Delta(x'_0,y'_2)| $, so
$$\frac{\beta_0\delta_0}{32 R_0^2} \, \diam(\T_Z(Q')) \leq |\Delta(x'_0, y'_1) - \Delta(x'_0,y'_2)| .$$
As above, denoting by $x''_0 \in \T_{Z_k}(\cc_m)$ the unique point such that $x''_0 \in W^{sc}(x'_0)$,
and using Lemma 4.3(b), we get
$$\frac{\beta_0\delta_0}{32 R_0^2} \, \diam(\T_Z(Q')) \leq |\Delta(x'_0, y'_1) - \Delta(x'_0,y'_2)| 
= |\Delta(x''_0, y'_1) - \Delta(x''_0,y'_2)| \leq C_1 \diam(\T_{Z_k}(\cc_m)) ,$$
so
$$\frac{\beta_0\delta_0}{32 R_0^2} \, \diam(\T_Z(Q')) \leq  \frac{C_1 \ep_1}{|b|} 
\leq  \frac{C_1}{\hrho^{q_1}} \frac{\hrho^{q_1} \ep_1}{|b|}
\leq \frac{C_1}{\hrho^{q_1}} \, \diam(\T_Z(\dd_j)) .$$
On the other hand, it follows from Lemma 4.1(b) with $\rho' = \rho_1$  that
$\diam(\T_Z(\dd_j)) \leq \rho_1^s\, \diam(\T_Z(Q')) .$
Thus,
$$\frac{\beta_0\delta_0}{32 R_0^2} \, \diam(\T_Z(Q')) 
\leq \frac{C_1}{\hrho^{q_1}} \,\rho_1^s\,  \diam(\T_Z(Q')) ,$$
so $\frac{\beta_0\delta_0 \hrho^{q_1}}{32 C_1 R_0^2} \leq \rho_1^s$, which implies
$s \leq \frac{1}{|\log \rho_1|}\, \left|\log \frac{\beta_0\delta_0 \hrho^{q_1}}{32 C_1 R_0^2} \right|$. 
\endofproof

\bs

Given $u, u'\in \hU$, we will denote by $\ell(u,u') \geq 0$ {\it the length of the 
smallest cylinder $Y(u,u')$ in $\hU$ containing} $u$ and $u'$.

Define the {\bf distance} $\dd(u,u')$ for $u,u' \in \hU$ by\footnote{Clearly $\dd$ depends on 
the cylinders $\cc_m$ and therefore
on the parameter $b$ as well.}: (i) $\dd(u,u') = 0$ if $u = u'$; (ii)  Let $u \neq u'$, and assume there exists 
$p \geq 0$ with $\sigma^p(Y(u,u')) \subset \cc'_m$, $\ell(u,u') \geq p$, for some $m = 1, \ldots, m_0$. 
Take the maximal $p$ with this property and the corresponding $m$ and set 
$\dd(u,u') = \frac{\dte(u,u')}{\diamtef (\cc_m)}$;  (iii) Assume $u \neq u'$,
however there is no $p \geq 0$ with the property described in (ii). Then set $\dd(u,u') = 1$.

Notice that $\dd(u,u') \leq 1$ always. 
Some other properties of $\dd$ are contained in the following, part (b) 
of which needs Lemma 6.3.

\bs


\noindent
{\bf Lemma 6.4.} {\it Assume that $u,u' \in \hU$, $u \neq u'$, and $\sigma^N(v) = u$, $\sigma^N(v') = u'$
for some $v,v'\in \hU$ with $\ell(v,v') \geq N$. Assume also that there exists 
$p \geq 0$ with $\sigma^p(Y(u,u')) \subset \cc'_m$, $\ell(u,u') \geq p$, for some $m = 1, \ldots, m_0$. }

\ms

(a) {\it We have $\dd(v,v') = \theta^N\, \dd (u,u')$.}

\ms

(b) {\it  Assume in addition that $\omega_J(v) < 1$ and $\omega_J(v') = 1$ for some $J \in \J(b)$. 
Then $p \leq t_0$ and}
$$|\omega_J (v) - \omega_J (v')| \leq \frac{\mu_0}{\theta^{t_0 + s_0 }}\, \dd(u,u') .$$

\ms

\noindent
{\it Proof.} 
(a) Let $p$ be the maximal integer with the given property and let $m \leq m_0$ correspond to $p$. 
Then $\sigma^{p+N}(Y(v,v')) \subset \cc'_m$, $\ell(v,v') \geq p+N$, and $p+N$ is the maximal 
integer with this property. Thus,
$\dd(v,v') = \frac{\dte(v,v')}{\diamtef (\cc_m)} = \theta^N \, \frac{\dte(u,u')}{\diamtef (\cc_m)} 
= \theta^N \, \dd(u,u') .$

\ms

(b)  $\omega_J(v) < 1$ means that $v\in \Xijl$ for some $(i,j,\ell) \in J$, and so 
$v = \vl_i(u)$ for some $u \in \dd'_j$. Then  $u = \sigma^N(v)$.  
If $u' \in \dd'_j$, then $v'' = \vl_i(u') \in \Xijl$ and $\sigma^N (v'') = u'$, so we must have
$v'' = v'$, which implies $\omega_J(v') = \omega_J(\vl_i(u')) = 1$, a contradiction.
This shows that $u'\notin \dd'_j$, and therefore
$\dte(u,u') \geq \diamte(\dd'_j) .$

Since $u\in \dd'_j$, $u' \notin \dd'_j $ and $\ell(u,u') \geq p$, it follows that 
$\sigma^p(u) \in \sigma^p(\dd'_j)$ and $\sigma^p(u') \notin \sigma^p(\dd'_j) $. 
On the hand, by assumption, $\sigma^p(u), \sigma^p(u') \in \cc'_m$.  Thus, the cylinder
$\sigma^p(\dd'_j)$ must be contained in $\cc'_m$. Now
Lemma 6.3 gives $p \leq t_0$ and the co-length $s$ of $\sigma^p(\dd'_j)$ in $\cc'_m$ is
$s \leq s_0$. If $\ell_m = \length(\cc_m)$, and $\ell = \length(\dd_j)$ we have 
$\ell-p -s = \length (\sigma^p(\dd'_j)) -s = \length (\cc'_m) = \ell_m$. Hence
$\ell = \ell_m + p + s \leq \ell_m + t_0 + s_0$, and using
$\dd(u,u') = \frac{\dte(u,u')}{\diamtef (\cc_m)}$, we get
\begin{eqnarray*}
|\omega (v) - \omega (v')| 
& = &  \mu_0 = \mu_0 \, \frac{\dte(u,u')}{\dte(u,u')} 
\leq \mu_0 \frac{\dte(u,u')}{\diamtef(\dd'_j) } 
=  \mu_0 \frac{\dte(u,u')}{  \theta^{\ell} }\\
& \leq & \mu_0 \frac{\dte(u,u')}{\theta^{\ell_m + t_0 + s_0}}  
= \mu_0 \frac{\dte(u,u')}{\theta^{t_0 + s_0}\diamtef (\cc_m) } 
= \frac{\mu_0}{ \theta^{t_0 + s_0}} \, \dd(u,u').
\end{eqnarray*}
This proves the lemma. 
\endofproof

\bs

Given $E > 0$ as in Sect. 5.2, let $\kk_E$ be {\bf the set of all functions $H \in \ff_\theta(\hU)$ 
such that } $H > 0$ on $\hU$ and $\frac{|H(u) - H(u')|}{H(u')} \leq E\, \dd(u,u')$ for all 
$u,u'\in \hU$ for which there exists an integer $p \geq 0$ with $\sigma^p(Y(u, u')) \subset \cc_m$ 
for some $m \leq m_0$ and $\ell(u,u') \geq p$.

Using Lemma 6.4 we will now prove the main lemma in this section.

\bs

\noindent
{\bf Lemma 6.5.} {\it For any $J \in \J(b)$ we have 
$\nn_J(\kk_E) \subset \kk_{E}$.} 

\bs

\noindent
{\it Proof.} 
Let  $u, u' \in \hU$  be such that  there exists
an integer $p \geq 0$ with $\sigma^p(Y(u, u')) \subset \cc_m$ for some $m = 1, \ldots,  m_0$ and $\ell(u,u') \geq p$.

Given $v \in \hU$ with $\sigma^{N}(v) = u$, let $C[\ii] = C[i_0, \ldots,i_{N}]$ 
be the  cylinder of length $N$ containing $v$.  
Set $\hC[\ii] = C[\ii]\cap \hU$. Then $\sigma^{N}(\hC[\ii]) = \hU_i$. Moreover, 
$\sigma^N : \hC[\ii] \longrightarrow \hU_i $ is a homeomorphism, so
there  exists a unique $v' = v'(v)\in \hC[\ii]$ such that $\sigma^N(v') = u'$.  Then 
$\dte (\sigma^j(v),\sigma^j(v'(v))) = \theta^{N-j}\, \dte (u,u')$ for all $j = 0,1, \ldots, N-1$. 
Also $\dte(v,v'(v))= \theta^N \dte (u,u')$ and $\dd(v,v'(v)) = \theta^N \dd(u,u')$. Using (5.3), we get
\begin{eqnarray}
|\fa_N(v) - \fa_N(v')| 
& \leq & \sum_{j=0}^{N-1} |\fa(\sigma^j(v)) - \fa (\sigma^j(v'))| \leq
\sum_{j=0}^{N-1} |\fa|_\theta \,\theta^{N-j}\, \dte (u,u')\nonumber\\
& \leq & \frac{T}{1-\theta}\, \dte (u,u'). 
\end{eqnarray}

 Let $J \in \J(b)$ and let $H \in \kk_E$. Set $\nn = \nn_J$.
We will show that $\nn H \in \kk_{E }$. 

Using the above  and the definition of $\nn = \nn_J$, and setting  $v' = v'(v)$ for brevity, we get
\begin{eqnarray*}
&        &\frac{|(\nn H)(u) - (\nn H)(u')|}{\nn H (u')} 
 =      \frac{\di \left| \sum_{\sigma^N v = u} e^{\fa_N(v)}\, \omega(v) H(v) -  
\sum_{\sigma^N v = u} e^{\fa_N(v'(v))}\, \omega(v'(v)) H(v'(v)) \right|}{\nn H (u')} \\
& \leq & \frac{\di \left| \sum_{\sigma^N v = u} 
e^{\fa_N(v)}\, [\omega(v) H(v) -  \omega(v') H(v')]\right|}{\nn H (u')}   +
\frac{\di \sum_{\sigma^N v = u}  
\left|e^{\fa_N(v)}-  e^{\fa_N(v')}\right| \, \omega(v') H(v')}{\nn H (u')} \\
& \leq & \frac{\di \sum_{\sigma^N v = u} e^{\fa_N(v) - \fa_N(v')} 
e^{\fa_N(v')}\, |\omega(v) - \omega(v')| H(v')}{\nn H (u')}   +
\frac{\di \sum_{\sigma^N v = u} e^{\fa_N(v)}\, \omega(v) | H(v) -   H(v')|}{\nn H (u')} \\
&       & + \frac{\di \sum_{\sigma^N v = u}  \left| e^{\fa_N(v) - \fa_N(v')} - 1\right| 
\,  e^{\fa_N(v')} \omega(v') H(v')}{\nn H (u')} .
\end{eqnarray*}

By the definition of $\omega$, either $\omega(v) = \omega(v')$ or at least one of these numbers is  $ < 1$.
Using Lemma 6.4 we then get 
$|\omega(v) - \omega(v')| \leq \frac{\mu_0}{\theta^{t_0 + s_0 }} \dd(u,u')$. Apart from that
$H \in \kk_E$ implies $|H(v) - H(v')| \leq E H(v') \dd(v,v') = E H(v') \theta^N \dd(u,u')$, while 
$\left| e^{\fa_N(v) - \fa_N(v')} - 1\right| \leq e^{T/(1-\theta)} \frac{T}{1-\theta}\, \dte (u,u')$. Thus,
\begin{eqnarray*}
&        & \frac{|(\nn H)(u) - (\nn H)(u')|}{\nn H (u')} 
\leq e^{T/(1-\theta)} \frac{\mu_0}{\theta^{t_0 + s_0 }} \, 
\frac{\di \sum_{\sigma^N v = u} e^{\fa_N(v')}\, \dd(u,u')\, H(v')}{\nn H (u')} \\
&        &  + \frac{\di \sum_{\sigma^N v = u} 
e^{\fa_N(v)- \fa_N(v')}\,e^{\fa_N(v')}\,  2\omega(v')\, E H(v') \theta^N \dd(u,u')}{\nn H (u')} 
    + e^{T/(1-\theta)} \frac{T}{1-\theta}\, \dte (u,u')\\
& \leq & 2 e^{T/(1-\theta)} \frac{\mu_0}{\theta^{t_) +s_0 + q_1}} \, \dd(u,u') 
+ 2 e^{T/(1-\theta)} E \theta^N \, \dd(u,u') +  e^{T/(1-\theta)} \diamte(\cc_m) \frac{T}{(1-\theta) }\, \dd (u,u')\\
& \leq & E \, \dd(u,u') , 
\end{eqnarray*}
using (\ref{eq:mucond}) and Lemma 5.1, and assuming $2 e^{T/(1-\theta)} \theta^N \leq 1/3$ and 
$e^{T/(1-\theta)} C_2 (\ep_1/|b|)^{\alpha_2} \frac{T}{(1-\theta) } \leq \frac{1 }{3}\leq \frac{E}{3}$;
the latter follows from $|b| \geq b_0$ and (\ref{eq:b0cond}). Hence $\nn H \in \kk_E$.
\endofproof

\def\tkk{\widetilde{\kk}}

\subsection{Main properties of the operators $\lab^N$}

Recall the numbers $\theta_1, \theta_2 \in (0,1)$ defined in Sect. 6.1.
Then using the proof of Lemma 5.1(c) and taking $C_2 > 0$ sufficiently large we have\footnote{Notice that
for (6.16) choosing $\theta_1$ with $\theta_1 ^{\beta'} \leq \theta$ would be enough. However in the beginning
of Sect. 6.1 we imposed a stronger condition on $\theta_1$ which will be used later on (see the end of the proof
of Theorem 1.1 in Sect. 7).}
\be
\diam_{\theta_1}(\cc) \leq C_2\, \diam(\Psi(\cc))
\ee
for any cylinder $\cc$ in $U$.

Throughout the rest of this section {\bf we assume that} $f \in \ff_{\theta_1}(\hU)$.

\medskip

Given points $u,u' \in U$ we will denote $\tu = \Psi(\pi_{\hz_0}(u))$ and $\tu' = \Psi(\pi_{\hz_0}(u'))$; 
these are then points on the true unstable manifold $W^u_{\ep_0}(\hz_0)$. In this section we will 
frequently work under the following assumption for points $u,u' \in \hU$ contained in some cylinder 
$\cc_m$ ($1\leq m \leq m_0$), an integer $p \geq 0$ and points $v,v'\in \hU$:
\be\label{eq:u-cond}
u,u'\in \cc_m \:  , \: \sigma^p(v) = \vli(u)\;, \; \sigma^p(v')  = \vli(u')\;,\; \ell(v,v') \geq N  ,
\ee
for some $i = 1,2$.
From (\ref{eq:u-cond}) we get $\ell(v,v') \geq N + p$ and $\sigma^{N+p}(v) = u$, $\sigma^{N+p}(v') = u'$. 
We will use the notation $\tcc_m = \tPsi(\cc_m) \subset \tR$.

The following estimate plays a central role in this section.

\bs

\noindent
{\bf Lemma 6.6.} {\it There exists a global constant $C_3 > 0$ independent of $b$ and $N$ 
such that if the points  $u,u' \in \hU$, the  cylinder $\cc_m$, the integer 
$p \geq 0$ and the points $v,v'\in \hU$ satisfy {\rm (\ref{eq:u-cond})}  for some 
$i = 1,2$ and $\ell = 1, \ldots, \ell_0$, and $w,w'\in \hU$ are such that $\sigma^N w = v$, 
$\sigma^N w' = v'$ and $\ell(w,w') \geq N$, then}
$$|\tau_N(w) - \tau_N(w')| \leq C_3 \, \theta_2^{p+N} \, \diam(\tcc_m) .$$

\ms

\noindent
{\it Proof.} Assume that the points $u,u',v,v',w,w'$ and the cylinder $\cc$
satisfy the assumptions in the lemma. Clearly, $\ell(w,w') \geq p+2N$ and
\be
\tau_N(w) - \tau_N(w') = [\tau_{p+2N}(w) - \tau_{p+2N}(w')] - [\tau_{p+N}(v) - \tau_{p+N}(v')] .
\ee

Recall the construction of the map $\vl_i$ from the proof of Lemma 5.4. In particular by (5.7),
$\pp^N(\vl(u)) = \pi_{\dl_i} (u)$, where we set $\dl_i = \dl_i(\hz_0) \in W^s_{R_1}(\hz_0)$ for brevity. 
Since $\sigma^{p}(v) = \vl_i(u)$ and $\sigma^p(v') = \vl_i(u')$, we have $\sigma^{p+N}(v) = u$ and 
$\sigma^{p+N}(v') = u'$, so $\pp^{p+N}(v), \pp^{p+N}(v') \in W^u_{R_1}(d')$ for some $d' \in W^s_{R_1}(\hz_0)$.
Moreover, $\pp^p(v) \in W^s_{R_1}(\vl_i(u))$ and the choice of $N$ imply (as in the proof of Lemma 5.4) that 
$d(\dl_i, d') < \delta''$, the constant from Lemma 4.4.
Similarly, $\pp^{p+2N}(w), \pp^{p+2N}(w') \in W^u_{R_1}(d'')$ for some $d'' \in W^s_{R_1}(\hz_0)$ with
$d(\dl_i, d'') < \delta''$.  Moreover, since the local stable/unstable holonomy maps are uniformly 
$\alpha_1$-H\"older (by the choice of $\alpha_1$), there exists a global constant $C'_3 > 0$ such that 
$d(d', d'') \leq C'_3 (d(\pp^{p+N}(v), \pp^{p+2N}(w)))^{\alpha_1}$.
Using this and (2.1) for points on local stable manifolds, i.e. going backwards along the flow, we get
$$d(d', d'') \leq C'_3 (d(\pp^{p+N}(v), \pp^{p+2N}(w)))^{\alpha_1} 
\leq C'_3 \left( \frac{d(v, \pp^N(w))}{c_0 \gamma^{p+N}} \right)^{\alpha_1} 
\leq \frac{C'_3}{c_0^{\alpha_1} \gamma^{\alpha_1 (p+N)}}  .$$
Hence
$$(d(d', d''))^{\beta_1} \leq  (C'_3/c_0^{\alpha_1})^{\beta_1} (1/\gamma^{p+N})^{\alpha_1 \beta_1} 
\leq  C''_3 \theta_2^{p+N} .$$

We are preparing to use Lemma 4.3. Set $\hu = \pi_{\hz_0}(u)$ and $\hu' = \pi_{\hz_0}(u')$.
Then for $\Psi(\hu) = \tu$ and
$\Psi(\hu') = \tu'$ we have $\tu = \phi_{t(u)} (\hu)$ and  $\tu' = \phi_{t(u')} (\hu')$ 
for some $t(u), t(u') \in \R$. So
$$\tau_{p+N}(v) - \tau_{p+N}(v') = \Delta(\pp^{p+N}(v), \pp^{p+N}(v')) = \Delta(\hu, \pi_{d'}(\hu'))
= \Delta(\tu, \pi_{d'}(\tu')) + t(u) - t(u') ,$$
and similarly
$$\tau_{p+2N}(w) - \tau_{p+2N}(w') = \Delta(\pp^{p+2N}(w), \pp^{p+2N}(w')) = \Delta(\hu, \pi_{d''}(\hu'))
= \Delta(\tu, \pi_{d''}(\tu')) + t(u) - t(u') .$$
This, (6.19), Lemma 4.3  and the above estimate yield
\begin{eqnarray*}
|\tau_N(w) - \tau_N(w')|
& =     & |[ \Delta(\tu, \pi_{d'}(\tu')) - t(u) + t(u')] - [ \Delta(\tu, \pi_{d''}(\tu')) - t(u) + t(u')]|\\
& =     & |  \Delta(\tu, \pi_{d'}(\tu')) -  \Delta(\tu, \pi_{d''}(\tu'))| 
\leq C_1 \diam (\tcc) \, (d(d',d''))^{\beta_1}\\
& \leq & C_1 C''_3 \theta_2^{p+N} \, \diam(\tcc) .
\end{eqnarray*}
This proves the lemma.
\endofproof

\bs

Let $M_0 > 0$ be a fixed constant (it is enough to take $M_0 = M + a_0$) and let 
$$E_1 = 2C_4 e^{C_4} \quad \mbox{\rm where } \quad
C_4 = \frac{T_0 C_2}{1-\theta} + M_0 C_3 ,$$ 
and $C_3 > 0$ is the constant from Lemma 6.6. Assume $N$ is so large that
$\theta_2^N e^{C_7} \leq 1/2$.

Denote by $\kk_0$
{\it the set of all $h \in \ff_{\theta}(\hU)$ such that $h \geq 0$ on $\hU$ and for any 
$u,u' \in \hU$ contained in some cylinder $\cc_m$ ($1\leq m \leq m_0$), any integer 
$p \geq 0$ and any points $v,v'\in \hU$ satisfying {\rm (\ref{eq:u-cond})} for some $i = 1,2$ and 
$\ell = 1, \ldots, \ell_0$ we have}
\be
|h(v) - h(v') | \leq E_1 \, \theta_2^{p+N} \, h(v')\, \diam(\tcc_m) .
\ee

We are going to show that the eigenfunctions $h_a \in \kk_0$ for $|a| \leq a'_0$ (see Sect. 5.1). 
This will be derived from the following.

\bs

\noindent
{\bf Lemma 6.7.} {\it For any real constant $s$ with $|s| \leq M_0$ we have 
$L^q_{f - s \tau}(\kk_0) \subset \kk_0$ for all integers $q \geq N$.} 

\bs

\noindent
{\it Proof.} We will use Lemma 6.6 and a standard argument.

Assume that $u,u' \in \hU$, the  cylinder $\cc$ in $U$, the integer $p \geq 0$ 
and the points $v,v'\in \hU$ satisfy (\ref{eq:u-cond}) for some $i = 1,2$ and $\ell = 1, \ldots, \ell_0$, and 
$w,w'\in \hU$ are such that $\sigma^N w = v$, $\sigma^Nw' = v'$ and $\ell(w,w') \geq N$; then $w' = w'(w)$ is
uniquely determined by $w$. 

Using $f \in \ff_{\theta_1}(\hU)$, the choice of $\theta_1$ and (6.15), we get 
\begin{eqnarray*}
|f_N (w) - f_N(w')| 
& \leq & \frac{T_0}{1-\theta_1} D_{\theta_1}(v,v') = \frac{T_0}{1-\theta_1} \theta_1^{p+N}\, D_{\theta_1}(u,u')\\
& \leq &  \frac{T_0}{1-\theta_1} \theta_1^{p+N}\, \diam_{\theta_1}(\cc) \leq C'_4 \theta_2^{p+N} \, \diam(\tcc) ,
\end{eqnarray*}
where $T_0 = |f|_\theta$ and $C'_4 = C_2 T_0/ (1-\theta_1)$. This and Lemma 6.6 imply
$$|(f-s\tau)_N(w) - (f- s\tau)_N(w')| \leq C''_4 \theta_2^{p+N} \, \diam(\tcc) $$
for all $s\in \R$ with $|s| \leq M_0$, where $C''_4 > 0$ is as above.

Thus, given $s$ with $|s| \leq M_0$ and $h \in \kk_0$ we have:
\begin{eqnarray*}
&        &|(L^N_{f-s\tau}h)(v) - (L^N_{f- s\tau} h)(v')|
 =      \left| \sum_{\sigma^N w = v} e^{(f-s \tau)_N(w)}\, h(w) 
 -  \sum_{\sigma^N w = v} e^{(f - s \tau)_N(w'(w))}\, h(w'(w)) \right| \\
& \leq & \left| \sum_{\sigma^N w = v} e^{(f- s\tau)_N(w)}\, [h(w) -  h(w')]\right|  +
\sum_{\sigma^N w = v}  \left|e^{(f- s\tau)_N(w)} -  e^{(f- s\tau)_N(w')}\right| \, h(w') \\
& \leq & \sum_{\sigma^N w = v} e^{(f- s\tau)_N(w) - (f- s\tau)_N(w')} 
e^{(f- s \tau)_N(w')}\,E_1 \theta_2^{p+2N} \diam (\tcc)\,  h(w') \\  
&       & + \sum_{\sigma^N w = v}  \left| e^{(f- s \tau)_N(w) - (f- s \tau)_N(w')} - 1\right| 
\,  e^{(f- s \tau)_N(w')} h(w') \\
& \leq & E_1 \theta_2^{p+2N} \diam (\tcc)\,e^{C''_4} \, (L^N_{f-s \tau} h)(v') 
+  e^{C_4} \,C_4 \theta_2^{p+N} \, \diam(\tcc)\, (L^N_{f-s \tau} h)(v')\\
& \leq & E_1 \,\theta_2^{p+N}\,   \diam(\tcc)\, (L^N_{f-s \tau} h)(v') , 
\end{eqnarray*}
since $e^{C_4} \,C_4 \leq E_1/2$ and $\theta_2^N \,e^{C_4} \leq 1/2$ by (\ref{eq:Ncond}).
Hence $L^N_{f-s \tau} h \in \kk_0$.
\endofproof

\bs

\noindent
{\bf Corollary 6.8.} {\it For any real constant $a$ with $|a| \leq a_0$ we have $h_a \in \kk_0$.} 

\bs

\noindent
{\it Proof.} Let $|a| \leq a_0$. Since the constant  function $h = 1 \in \kk_0$, it follows from 
Lemma 6.7 that $L^{mN}_{f-(P+a)\tau} 1 \in \kk_0$ for all $m \geq 0$. Now the Ruelle-Perron-Frobenius 
Theorem (see e.g. \cite{PP}) and the fact that $\kk_0$ is closed in $\ff_\theta(\hU)$ imply 
$h_a\in \kk_0$. 
\endofproof

\subsection{The main estimate for $\lab^N$}

We will now define a class of pairs of functions similar to $\kk_0$ however involving the parameter $b$. 
We continue to {\bf assume that} $f \in \ff_{\theta_1}(\hU)$.

Denote by $\kk_b$ 
{\it the set of all pairs $(h, H)$ such that $h \in \ff_{\theta}(\hU)$, 
$H \in \kk_E$ and the following two properties hold:}

(i) $|h| \leq H$ on $\hU$,

(ii) for any $u,u' \in \hU$ contained in a cylinder $\cc_m$ for some $m = 1, \ldots, m_0$, 
any integer  $p \geq 0$ and any points $v,v'\in \hU_1$
satisfying (\ref{eq:u-cond})  for some $i = 1,2$ and $\ell = 1, \ldots, \ell_0$ we have
\be
|h(v) - h(v') | \leq E \, |b|\, \theta_2^{p+N} \, H(v') \, \diam(\tcc_m) .
\ee
Recall that here $\tcc_m = \tPsi(\cc_m)$.

Our aim in this section is to prove the following.

\ms

\noindent
{\bf Lemma 6.9.} {\it Choosing $E > 1$ and $\mu_0$ as in Sect. {\rm 5.2} and assuming $N$ is sufficiently large, 
for any $|a| \leq a'_0$, any $|b| \geq b_0$ and any $(h,H) \in \kk_b$ 
there exists $J \in \J(b)$ such that $(\lab^{N} h, \nn_J H ) \in \kk_b$.}

\bs

To prove this  we need the following lemma, whose proof is essentially the same as that of 
Lemma 14 in \cite{D}.  For completeness we prove it in the Appendix.

\bs

\noindent
{\bf Lemma 6.10.}  {\it Let $(h,H) \in \kk_b$. Then for any $m \leq m_0$,  any $j = 1, \ldots, j_0$
with $\dd_j \subset \cc_m$, any $i = 1,2$ and  $\ell = 1,\ldots,  \ell_0$ we have:}

(a) {\it $\di\frac{1}{2} \leq \frac{H(\vli (u'))}{H(\vli (u''))} \leq 2$ for all} $u', u'' \in \dd'_j$;

(b) {\it Either for all $u\in \dd'_j$ we have
$|h(\vli (u))|\leq \frac{3}{4}H(\vli (u))$, or $|h(\vli (u))|\geq \frac{1}{4}H(\vli (u))$ 
for all $u\in \dd'_j$.} 


\bs

\noindent
{\it Proof of Lemma} 6.9. 
The constant $E_1 > 1$ from Sect. 6.4 depends only on $C_4$, and we take $N$ so large that
$E_1 \theta^N \leq 1/4$; then  $C_4 \theta^N \leq 1/2$ holds, too.

Let $|a| \leq a'_0$, $|b| \geq b_0$ and $(h, H) \in \kk_b$. We will construct a representative set 
$J \in \J(b)$ such that  $(\lab^N h , \nn_J H) \in \kk_b$. 

Consider for a moment an arbitrary (at this stage) representative set $J$. 
We will first show that  $(\lab^N h , \nn_J H) $ has property (ii). 

Assume that the points $u,u'$, the cylinder $\cc_m$ in $U$, the integer $p \geq 0$ and the points 
$v,v'\in \hU$ satisfy (\ref{eq:u-cond})  for some $i = 1,2$ and $\ell = 1, \ldots, \ell_0$. 


From the definition of $\fa$, for any $w, w'$ with $\sigma^N w = v$, $\sigma^N(w') = v'$ 
and $\ell(w,w') \geq N$ we have 
\begin{eqnarray*}
\fa_N(w) 
& =  & f_N(w) - (P+a)\tau_N(w) + (\ln h_a - \ln h_a\circ \sigma)_N(w) - N \lambda_a\\
& =  & f_N(w) - (P+a)\tau_N(w) + \ln h_a (w)  - \ln h_a (v) - N \lambda_a .
\end{eqnarray*}
Since $h_a \in \kk_0$ by Corollary 6.8,
$$|\ln h_a(w) - \ln h_a(w') | \leq \frac{| h_a(w) - h_a(w')|}{\min \{ |h_a(w)|, |h_a(w')|\}} 
\leq E_1 \,\theta_2^{p+2N} \, \diam(\tcc_m) , $$
and similarly, $|\ln h_a(v) - \ln h_a(v') | \leq E_1 \,\theta_2^{p+N} \, \diam(\tcc_m)$.
Using this and Lemma 6.6, as in the proof of Lemma 6.7 we get
\begin{eqnarray}
|\fa_N(w) - \fa_N(w')| 
& \leq & C_4 \theta_2^{p+2N} \, \diam(\tcc_m) + 2E_1 \, \theta_2^{p+N}\, \diam(\tcc_m)\nonumber\\
& \leq & (C_4 + 2 E_1)  \, \theta_2^{p+N}\, \diam(\tcc_m)  \leq 1 ,
\end{eqnarray}
by the choice of $N$.

Hence for any $a$ and $b$ with $|a| \leq a'_0$ and $|b| \geq b_0$   we have:
\begin{eqnarray*}
&        &|(\lab^N h)(v) - (\lab^N  h)(v')|
 =      \left| \sum_{\sigma^N w = v} e^{(\fa_N - \i b \tau_N) (w)}\, h(w) 
 -  \sum_{\sigma^N w = v} e^{(\fa_N - \i b \tau_N)(w'(w))}\, h(w'(w)) \right| \\
& \leq & \left| \sum_{\sigma^N w = v} e^{(\fa_N - \i b \tau_N)(w)}\, [h(w) -  h(w')]\right|  +
\sum_{\sigma^N w = v}  \left|e^{(\fa_N - \i b \tau_N)(w)} -  e^{(\fa_N - \i b \tau_N) (w')}\right| \, |h(w')| \\
& \leq & \sum_{\sigma^N w = v} e^{(\fa_N (w) - \fa_N(w')} e^{\fa_N (w')}\,E |b| 
\theta_2^{p+2N} \diam (\tcc_m)\,  H(w') \\  
&       & + \sum_{\sigma^N w = v}  \left| e^{(\fa_N - \i b \tau_N)(w) - 
(\fa_N - \i b  \tau_N)(w')} - 1\right| \,  e^{\fa_N(w')} H(w') \\
& \leq & e\, E |b|\, \theta_2^{p+2N} \diam (\tcc_m) \, (\ma^N H)(v') 
+  e \,(C_4 + 2E_1 +  C_3 |b| ) \theta_2^{p+N} \, \diam(\tcc_m)\, (\ma^N H)(v')\\
& \leq & [2 e\, E \theta_2^N  + 2 e (C_4 + 2E_1 + C_3) ]  \,|b|\, \theta_2^{p+N}\,   
\diam(\tcc_m)\, (\nn_J H)(v') \leq  E |b|\, \theta_2^{p+N}\,   \diam(\tcc_m)\, (\nn_J H)(v') , 
\end{eqnarray*}
assuming $2 e \theta^N \leq 1/2$ and $2 e (C_4 + 2E_1 + C_3)  \leq E/2$. 
Thus, $(\lab^N h , \nn_J H)$ has property (ii).

So far the choice of $J$ was not important. We will now construct a representative set $J$ so that 
$(\lab^N h , \nn_J H)$ has property (i), namely 
\be\label{eq:Lab}
|\lab^N h|(u) \leq (\nn_J H) (u)
\ee
for all $u \in \hU$. 

Define the functions $\psi_\ell, \gao_\ell, \gat_\ell : \hU  \longrightarrow \C$  by
$$\di \psi_\ell(u) = e^{(\fa_{N}+\i b\tau_{N})(\vl_1(u))} h(\vl_1(u)) 
+ e^{(\fa_{N}+\i b\tau_{N})(\vl_2(u))} h(\vl_2(u)) ,$$
$$\di \gao_\ell(u) = (1-\mu_0)\, e^{\fa_{N} (\vl_1(u))} H(\vl_1(u)) + e^{\fa_{N}(\vl_2(u))} H(\vl_2(u)) ,$$
while $\gat_\ell(u)$ is defined similarly with a coefficient $(1-\mu_0)$ 
in front of the second term. 

Recall the functions $\varphi_\ell(u) = \varphi_\ell(\hz_0, u)$, $u \in U$,  from Sect. 5.3. 

Notice that (\ref{eq:Lab}) is trivially satisfied for $u \notin V_b$ for any choice of $J$.

Consider an arbitrary  $m = 1, \ldots, m_0$. We will construct  $j \leq j_0$ with $\dd_{j} \subset \cc_m$,  
and a pair  $(i, \ell)$ for which $(i,j,\ell)$ will be included in $J$.

\ms

\noindent
{\bf Case 1.}
There exist $j\leq j_0$ with $\dd_j \subset \cc_m$, $i = 1,2$ and $\ell \leq \ell_0$   such that the 
first alternative in Lemma 6.10(b) holds for $\hdd_j$,  $i$ and $\ell$. 
For such $j$, choose $i= i_j$ and $\ell = \ell_j$ with this property and include $(i,j,\ell)$ in $J$. 
Then $\mu_0 \leq 1/4$ implies $|\psi_\ell(u)| \leq \gamma^{(i)}_\ell(u)$ for all 
$u\in \hdd'_j$, and regardless how the rest
of $J$ is defined, (\ref{eq:Lab}) holds for all $u\in \hdd'_j$, since
\begin{eqnarray}
 \left| (\lab^N h)(u)\right| 
& \leq & \left| \sum_{\sigma^N v = u, \;v\neq \vl_1(u),\vl_2(u)} e^{(\fa_N+\i b\tau_N)(v)} h(v) \right| 
+ |\psi_\ell(u)|\nonumber \\
& \leq & \sum_{\sigma^N v = u, \;v\neq \vl_1(u),\vl_2(u)} e^{\fa_N(v)} |h(v)| 
+ \gamma^{(i)}_\ell(u)\nonumber\\
& \leq & \sum_{\sigma^N v = u, \;v\neq \vl_1(u),\vl_2(u)} e^{\fa_N(v)} \omega(v) H(v)\nonumber\\
&      & + \left[e^{\fa_N(\vl_1(u))} \omega_J(\vl_1(u)) H(\vl_1(u)) 
+  e^{\fa_N(\vl_2(u))} \omega_J(\vl_2(u)) H(\vl_2(u))\right]  \nonumber\\
& \leq &  (\nn_J H) (u) .
\end{eqnarray}

\noindent
{\bf Case 2.}
For all $j\leq j_0$ with $\dd_j \subset \cc_m$, $i = 1,2$ and $\ell \leq \ell_0$ the 
second alternative in Lemma 6.10(b) holds for $\hdd_j$,  $i$ and $\ell$, i.e.
\be
|h(\vl_i(u))|\geq \frac{1}{4}\, H(\vl_i(u)) > 0
\ee
for any  $u \in \hcc'_m$.


Let $u,u' \in \hcc'_m$, and let $i = 1,2.$ 
Using (6.19) and the assumption that $(h,H) \in \kk_b$, and in particular property 
(ii) with $p = 0$, $v = \vl_i(u)$ and
$v' = \vl_i(u')$, and assuming e.g.\\   $\min\{ |h(\vl_i(u))| , |h(\vl_i(u'))|\}  = |h(\vl_i(u'))|$, we get
\begin{eqnarray*}
\frac{|h(\vl_i(u)) - h(\vl_i(u'))|}{\min\{ |h(\vl_i(u))| , |h(\vl_i(u'))| \}}
 \leq  \frac{E|b|\,\theta_2^N H(\vl_i(u'))}{|h(\vl_i(u'))| } \diam (\Psi(\cc_m)) 
 \leq  4 E|b|\,   \theta_2^{N} \frac{\ep_1}{|b|} = 4 E \theta_2^N \ep_1 .
\end{eqnarray*}
So, the angle between the complex numbers 
$h(\vl_i(u))$ and $h(\vl_i(u'))$ (regarded as vectors in $\R^2$)  is  
$< 8 E \theta_2^N \ep_1 < \pi/6$ by (\ref{eq:Ncond}).  
In particular,  for each $i = 1,2$ we can choose a real continuous
function $\thetam_i(u)$, $u \in  \hcc'_m$, with values in $[0,\pi/6]$  and a constant $\lambdam_i$ such that
$$\di h(\vl_i(u)) = e^{\i(\lambdam_i + \thetam_i(u))}|h(\vl_i(u))| \quad , \quad u\in \hcc'_m .$$
Fix an arbitrary $u_0\in \hcc'_m$ and set $\lambdam = |b| \varphi_\ell(u_0)$. 
Replacing e.g $\lambdam_2$ by $\lambdam_2 +  2r \pi$ for some integer $r$, 
we may assume that $|\lambdam_2 - \lambdam_1 + \lambdam | \leq \pi$.

Using the above, $\theta \leq 2 \sin \theta$ for $\theta \in [0,\pi/3]$,  
and some elementary geometry  yields\\
$|\thetam_i(u) - \thetam_i(u')|\leq 2 \sin |\thetam_i(u) - \thetam_i(u')| < 16 E \theta_2^N \ep_1$ for all
$u,u' \in \hcc'_m$.

The difference between the arguments of the complex numbers
$e^{\i \,b\,\tau_N(\vl_1(u))} h(\vl_1(u))$ and $e^{\i \,b\, \tau_N(\vl_2(u))} h(\vl_2(u))$
is given by the function
\begin{eqnarray*}
\Gl(u) 
& = & [b\,\tau_N(\vl_2(u)) + \thetam_2(u) + \lambdam_2] -  
[b\, \tau_N(\vl_1(u)) + \thetam_1(u) + \lambdam_1]\\
& = & (\lambdam_2-\lambdam_1) + |b| \varphi_\ell(u) + (\thetam_2(u) - \thetam_1(u)) .
\end{eqnarray*}

Using Lemma 5.4, we can now choose $j  =1, \ldots, j_0$ and $j = 1, \ldots, \tj$ with $j \neq j'$
such that  $\dd_j, \dd_{j'} \subset \cc_m$ (and $\dd_j \cap P_1 \cap P_0 \cap \Xib_B \neq \e$ by the choice of $j$)
and  $\ell = 1, \ldots, \ell_0$ such that
for all $u \in \hdd'_j $ and $u'\in \hdd'_{j'}$  we have 
\be
\frac{\hd \hrho \ep_1}{|b|} \leq \hd_0 \, \diam(\Psi(\cc_m)) \leq |\varphi_\ell(u) - \varphi_\ell(u')| \leq 
C_1 \, \diam(\Psi(\cc_m)) \leq C_1 \frac{\ep_1}{|b|} .
\ee
Fix $\ell_m = \ell$ with this property.
Then for $u \in \hdd'_j$ and $u'\in \hdd'_{j'}$ we have
\begin{eqnarray*}
|\Gl(u)- \Gl(u')|
& \geq & |b|\, |\varphi_\ell(u) - \varphi_\ell(u')| - |\thetam_1(u) - \thetam_1(u')|
- |\thetam_2(u) - \thetam_2(u')|\\
& \geq & \hd \hrho \epsilon_{1} - 32 E \theta_2^N \ep_1 >  2\ep_3 ,
\end{eqnarray*}
since $32 E \theta_2^N  < \hd_0 \hrho/2$ by (5.6), where 
$\di \ep_3 = \frac{ \hd_0 \hrho \epsilon_{1}}{4} .$

Thus,  $|\Gl(u)- \Gl(u')|\geq 2\epsilon_{3}$ for all
$u\in \hdd'_j$ and $u'\in \hdd'_{j'}$. Hence either 
$|\Gl(u)| \geq \ep_3$ for all $u\in \hdd'_j$ or
$|\Gl(u')| \geq \ep_3$ for all $u'\in \hdd'_{j'}$.

Assume for example that $|\Gl(u)| \geq \ep_3$ for all $u\in \hdd'_j$. 
On the other hand, (6.26) and the choice of $\ep_1$ imply that
for any $u \in \hcc'_m$ we have
$$|\Gamma_\ell(u)| \leq |\lambdam_2-\lambdam_1 +\lambdam | + 
|b|\, |\varphi_\ell(u) -\varphi_\ell(u_0)| + |\thetam_2(u) - \thetam_1(u)| 
\leq \pi + C_1 \ep_1 + 16 E \theta_2^N \ep_1 < \frac{3\pi}{2} .$$
Thus, $\ep_3 \leq |\Gl(u)| <  \frac{3\pi}{2}$ for all $u \in \hdd'_j$.

Hence, we see that for $u\in \hdd'_j$ the difference $\Gl(u) $ between the 
arguments of the complex numbers
$e^{\i \,b\,\tau_N(\vl_1(u))} h(\vl_1(u))$ and $e^{\i \,b\, \tau_N(\vl_2(u))} h(\vl_2(u))$, 
defined as a number in the interval $[0, 2\pi)$, satisfies
$\Gl(u) \geq \ep_3$ for all $u\in \hdd'_j .$

\def\halpha{\hat{\alpha}}
\def\hgamma{\hat{\gamma}}

As in \cite{D} it follows from Lemma 6.10 that either $H(\vl_1(u)) \geq H(\vl_2(u))/4$ for all
$u \in \dd'_j$ or $H(\vl_2(u)) \geq H(\vl_1(u))/4$ for all $u \in \dd'_j$. Indeed, fix
an arbitrary $u'\in \dd'_j$ and assume e.g. $H(\vl_1(u')) \geq H(\vl_2(u'))$. Then for any $u \in \dd'_j$
using Lemma 6.10(a) twice we get 
$H(\vl_1(u)) \geq H(\vl_1(u'))/2 \geq H(\vl_2(u'))/2 \geq H(\vl_2(u))/4 .$
Similarly, if $H(\vl_2(u')) \geq H(\vl_1(u'))$, then $H(\vl_2(u)) \geq H(\vl_1(u))/4$ for all $u \in \dd'_j$.

Now assume e.g. that $H(\vl_1(u)) \leq H(\vl_2(u))/4$ for all $u \in \dd'_j$.
As in \cite{D} (and \cite{St1}) we derive that 
$|\psi_\ell(u)| \leq \gamma^{(1)}_\ell (u)$ for all $u \in \hdd'_j$.

This completes the construction of the set $J = \{ (i_m, j_m, \ell_m) : m = 1, \ldots, m_0\} \in \J(b)$ 
and also the proof  of (\ref{eq:Lab}) for all $u \in V_b$. As we mentioned in the beginning of the proof, (\ref{eq:Lab}) 
always holds for $u \in \hU\setminus V_b$.
\endofproof


\section{Proofs of the Main Results}
\setcounter{equation}{0}

Here we prove Theorems 1.3 and 1.1 and Corollary 1.4. The main step is to obtain
$L^1$-contraction estimates for large powers of the contraction operators.
using the properties of these operators on $K_0$ and the strong mixing
properties of the shift map $\pp : R \longrightarrow R$.

For any $J \in \J(b)$ set
$W_J = \cup\{ \hdd'_j : (i,j,\ell) \in J \; \mbox{\rm for some }\; i, \ell \} \subset V_b .$
Using Lemma 6.3 and the class of functions $\kk_E$ we will now prove the following 
important estimates\footnote{This should be regarded as the analogue of 
Lemma 12 in \cite{D} (and Lemma 5.8 in \cite{St2}).}.

\bs

\noindent
{\bf Lemma 7.1.}  {\it Let $f \in \ff_{\theta_1}(\hU)$. }

\medskip

(a) {\it There exists a global constant $C''_5 > 0$ such that for any $H\in \kk_E$  
and any $J\in \J(b)$ we have}
\be
\int_{V_b} H^2 \, d\nu \leq C''_5\, \int_{W_J} H^2 \, d\nu .
\ee

\ms

(b) {\it For any $H\in \kk_E$   and any $J\in \J(b)$ we have
\be
\int_{V_b} (\nn_J H)^2 \, d\nu \leq \rho_3\, \int_{V_b} L^N_{\f0} (H^2) \, d\nu ,
\ee
where $\di \rho_3 = \rho_3(N) = \frac{e^{a_0NT}}{1+ \frac{\mu_0 e^{-NT}}{C''_5}} < 1$,
assuming that $a_0 > 0$ is sufficiently small.}

\bs

\noindent
{\it Proofs.} (a) Let $H \in \kk_E$ and $J \in \J(b)$. Consider an arbitrary 
$m = 1, \ldots,m_0$. There exists $(i_m,j_m,\ell_m) \in J$ such that $\dd_{j_m} \subset \cc_m$.
It follows from (\ref{eq:Gibbs}) that there exists a global constant $\omega_0 \in (0,1)$ such that
 $\frac{\nu(\dd'_{j_m})}{\nu(\cc'_m)} \geq 1- \omega_0$.
Since $H \in \kk_E$, for any $u,u'\in \cc'_m$ we have $\frac{|H(u) - H(u')|}{H(u')} \leq E\dd(u,u') \leq E$, so
$H(u)/H(u') \leq 1+ E \leq 2E$. Thus, if $M_1 = \max_{\cc'_m} H$ and  $M_2 = \min_{\cc'_m}H$ 
we have $M_1/M_2 \leq 2E$. This gives
$$\int_{\cc'_m} H^2 \, d\nu \leq M_1^2 \nu(\cc'_m) 
\leq  \frac{4E^2}{1-\omega_0}\, \int_{\dd'_{j_m}} H^2\, d\nu . $$
Hence
$$\int_{V_b} H^2 \, d\nu \leq \sum_{m=1}^{m_0} \int_{\cc'_m} H^2 \, d\nu  
\leq  \frac{4E^2}{1-\omega_0}\, \sum_{m=1}^{m_0}  \int_{\dd'_{j_m}} H^2\, d\nu
\leq C''_5\, \int_{W_J} H^2 \, d\nu ,$$
with $C''_5 =  \frac{4E^2}{1-\omega_0}$, since  $\cup_{m=1}^{m_0} \hdd'_{j_m} = W_J$. 
This proves (7.1).

\ms

(b) Let again $H \in \kk_E$ and $J \in \J(b)$. By Lemma 6.4, $\nn_J H \in \kk_E$, while
the Cauchy-Schwartz  inequality implies
$$(\nn_J H)^2 = (\ma^N \omega H)^2 \leq (\ma^N \omega_J^2)\, (\ma^N H^2) \leq 
(\ma^N \omega_J)\, (\ma^N H^2) \leq \ma^N H^2 .$$
Notice that if $u \notin W_J$, then $\omega_{J}(u) = 1$. Let $u \in W_J$; then $u \in \hdd'_j$ 
for some (unique) $j\leq j_0$, and there exists a unique $(i(j), j, \ell(j)) \in J$. 
Set $ i = i(j)$, $\ell = \ell(j)$ for brevity.
Then $\vl_i(u) \in \hXijl$, so $\omijl(\vl_i(u)) = 1$, and therefore  
$\omega(v^{(\ell)}_i(u)) \leq 1-\mu_0\, \omijl(\vl_i(u)) = 1-\mu_0$. 
In fact,  if $\sigma^N(v) = u$ and $\omega(v) < 1$, then
$\omega^{(\ell')}_{i',j'}(v) = 1$ for some $(i',j',\ell') \in J$, so $v \in X^{(\ell')}_{i',j'}$. Then
$u = \sigma^N(v)  \in \sigma^N(X^{(\ell')}_{i',j'}) = \hdd'_{j'}$. 
Thus, we must have $j' = j$, and since for a given $j$, there is 
only one element $(i(j) , j , \ell(j) )$ in $J$, we must have also 
$i' = i(j)$ and $\ell' = \ell(j)$. Assuming e.g. that
$i = 1$,  this implies $v = \vl_1(u)$. Thus, 
\begin{eqnarray*}
(\ma^N \omega_{J})(u)
& = &  \sum_{\sigma^N v = u, \;v\neq \vl_1(u)} e^{\fa_N(v)}   + e^{\fa_N(\vl_1(u))} \omega_{J} (\vl_1(u))  \\
& =  & \sum_{\sigma^N v = u, \;v\neq \vl_1(u)} e^{\fa_N(v)}   + (1-\mu_0)  e^{\fa_N(\vl_1(u))} \\ 
&  =   & \sum_{\sigma^{N}v = u} e^{\fa_{N}(v)} - \mu_0 \,e^{\fa_{N}(\vl_1(u))}  
\leq  (\ma^{N}\; 1)(u) - \mu_0 \, e^{-N T} = 1 - \mu_0 \, e^{-NT} .
\end{eqnarray*}
This holds for all $u \in W_J$, so  $(\nn_J H)^2 \leq (1- \mu_0 e^{-NT})\, (\ma^N H^2)$ on $W_J$.
Using this and part (a) we get:
\begin{eqnarray*}
\int_{V_b} (\nn_J H)^2\, d\nu
& =     & \int_{V_b \setminus W_J} (\nn_J H)^2 \, d\nu +  \int_{W_J} (\nn_J H)^2 \, d\nu\\
& \leq &  \int_{V_b \setminus W_J} (\ma^N H)^2 \, d\nu 
+ (1-\mu_0 e^{-NT})\,  \int_{W_J} (\ma^N H)^2 \, d\nu\\
& =     & \int_{V_b} (\ma^N H)^2 \, d\nu  -\mu_0 e^{-NT}\,  \int_{W_J} (\ma^N H)^2 \, d\nu\\
& \leq & \int_{V_b} (\ma^N H)^2 \, d\nu  - \mu_0 e^{-NT}\,  \int_{W_J} (\nn_J H)^2 \, d\nu \\
& \leq & \int_{V_b} (\ma^N H)^2 \, d\nu  - \frac{\mu_0 e^{-NT}}{C''_5} \,  \int_{V_b} (\nn_J H)^2 \, d\nu .
\end{eqnarray*}
From this and 
$$(\ma^N H)^2 \leq (\ma^N 1)^2 (\ma^NH^2) \leq \ma^N H^2 
= L^N_{\f0} (e^{\fa_N - \f0_N} H^2)\leq e^{a_0 NT} (L^N_{\f0} H^2) ,$$ 
we get
$$(1+\mu_0 e^{-NT}/C''_5)\, \int_{V_b} (\nn_J H)^2\, d\nu \leq  \int_{V_b} (\ma^N H)^2 \, d\nu
\leq e^{a_0NT}\, \int_{V_b} L^N_{\f0} H^2 \, d\nu .$$
Thus (7.2) holds with $\di \rho_3 = \frac{e^{a_0NT}}{1+ \frac{\mu_0 e^{-NT}}{C''_5}} > 0$. 
Taking $a_0 = a_0(N) > 0$ sufficiently small, we have $\rho_3 < 1$.
\endofproof

\bs

We can now prove that iterating sufficiently many contraction operators provides an $L^1$-contraction on $U$.

Recall the set $\Lambda_N(b)$ defined  by (\ref{eq:Lambda}) and Lemma 6.1.
Set 
$$ \rho_3 = \frac{e^{2a_0NT}}{1+ \frac{\mu_0 e^{-NT}}{C''_5}} < 1 \quad , 
\quad R = e^{2a_0 N T} > 1 \quad, \quad h = \rho_3\, \chi_{V_b} + R\, \chi_{U\setminus V_b} ,$$
and notice that $\rho_3$ is as in Lemma 7.1.
We will assume $a_0 = a_0(N) > 0$ is chosen so small that
$$\frac{e^{4a_0NT}}{1+ \frac{\mu_0 e^{-NT}}{C''_5}} <
\rho_4 = \frac{1}{1+ \frac{\mu_0 e^{-NT}}{2C''_5}}  < 1 \quad , \quad 8a_0 N^2T < c' ,$$
where $c' > 0$ is the constant from (3.21). Notice that the latter and (6.7) imply
$8 a_0NTL < c_4$. Moreover, $\rho_3 R = \frac{e^{4a_0NT}}{1+ \frac{\mu_0 e^{-NT}}{C''_5}} < \rho_4$.
Assuming $L $ is taken sufficiently large, we will assume that
\be
L |\log \rho_4| = L \log \left( 1+ \frac{\mu_0 e^{-NT}}{2C''_5} \right) > \frac{16}{\beta_2} = 8s .
\ee

\ms

\noindent
{\bf Lemma 7.2.}  {\it Let $f \in \ff_{\theta_1}(\hU)$. Given the number $N$, there exist  constants $C_5 \geq 1$,   
$a_0 = a_0(N) > 0$  and $b_0 = b_0(N) \geq 1$
such that for any $|a| \leq a_0$ and $|b| \geq b_0$ and any sequence $J_1, J_2, \ldots, J_r \ldots $ of 
representative subsets of $\J(b)$, setting
$H^{(0)} = 1$ and $H^{(r+1)} = \nn_{J_r} (H^{(r)})$ ($r \geq 0$) we have}
\be
\di\int_{U} (H^{(L \hb)})^2 \, d\nu \leq \frac{C_5}{|b|^{8s}} .
\ee

\ms

\noindent
{\it Proof of Lemma} 7.2. Set  $\omega_r = \omega_{J_r}$, $W_r = W_{J_r}$ and $\nn_r = \nn_{J_r}$.
Since $H^{(0)} = 1 \in \kk_E$, it follows from Lemma 6.5 that $H^{(r)} \in \kk_E$ for all $r \geq 1$. 

Using $L^N_{\fa}((h\circ \sigma^N) \, H) = h\, (L_{\f0}^N H)$ and Lemma 7.1(b) we get
\begin{eqnarray*}
\int_{U} (H^{(m)})^2 \, d\nu
& =      & \int_{V_b} (H^{(m)})^2 \, d\nu + \int_{U\setminus V_b} (H^{(m)})^2 \, d\nu\\
& \leq  & \rho_3\, \int_{V_b} L^{N}_{\f0} (H^{(m-1)})^2\, d\nu 
+ e^{a_0 NT}\, \int_{U\setminus V_b} L^N_{\f0} (H^{(m-1)})^2\, d\nu\\
& =     & \int_{U} h\,  (L_{\f0}^N (H^{(m-1)})^2)\, d\nu 
= \int_{U}   L_{\f0}^N ( (h \circ \sigma^N)\, (H^{(m-1)})^2)\, d\nu\\
& =     &  \int_{U} (h\circ \sigma^N)\,  (H^{(m-1)})^2\, d\nu .
\end{eqnarray*}
Similarly,
\begin{eqnarray*}
\int_{U} (h\circ \sigma^N)\,  (H_{m-1})^2 \, d\nu
& \leq  &  \int_{U}  (h\circ \sigma^{2N})\, (h\circ \sigma^N)\,  (H^{(m-2)})^2\, d\nu .
\end{eqnarray*}
Continuing by induction and using $H^{(0)} = 1$, we get
\be
\int_{U} (H^{(m)})^2 \, d\nu \leq 
 \int_{U}  (h\circ \sigma^{mN})\, (h\circ \sigma^{(m-1)N}) 
 \ldots  (h\circ \sigma^{2N})\, (h\circ \sigma^N)\, d\nu .
\ee


Let $m = L \hb$ and let  $\delta = (d+3) \hd_0 \in (0,1)$, where $\hd_0$  
is as in Sects. 3 and 5. Set 
$$W = \{ x \in U : x\in \sigma^{-jN}(U\setminus V_b) \: \mbox{\rm for at least
$\delta m$ values of }\: j = 0,1, \ldots,m-1 \} .$$
Since $K_0 \subset V_b$, for such $j$ we have  $x\in \sigma^{-jN}(U\setminus K_0)$.
Notice that 
\be
(\piU)^{-1}(W) \subset \Lambda_{N} (b) ,
\ee
the set defined by (\ref{eq:Lambda}). Indeed, if $x\in W$ and $y \in W^s_R(x)$, then for any $j = 0,1, \ldots, m-1$ with
$\sigma^{jN}(x) \notin K_0 = \piU(P_1\cap P_0 \cap \Omb)$, since $\piU(\pp^{jN}(y)) = \sigma^{jN}(x)$,
we have $\pp^{jN}(y) \notin P_1 \cap P_0 \cap \Omb$. Thus, the latter holds for at least $\delta m$ values of
$j = 0, 1, \ldots, m-1$, so $\pp^i(y) \notin P_1 \cap P_0 \cap \Omb$ for at least $\delta m = \frac{\delta}{N} N L \hb$ values of
$i = 0,1, \ldots, NL \hb-1$. It follows from (\ref{eq:Lambda}) that $y \in \Lambda_{N}(b)$. This proves (7.6), and now Lemma 6.1 implies
\be
\nu(W) \leq \frac{2C_0}{|b|^{c_3}} .
\ee

Notice that if $x\in U\setminus W$, then $x\in \sigma^{-j N}(V_b)$ for at least $(1-\delta)m$ values of
$j = 0,1, \ldots, m-1$, so $(h\circ \sigma^{jN})(x) = \rho_3$ for that many $j's$.  Thus, (7.5) gives
\begin{eqnarray*}
\int_U (H^{(m)})^2\, d\nu 
& \leq & \int_{U\setminus W} \prod_{j=1}^{m} (h\circ \sigma^{jN}) \, d\nu  
+ \int_{W} \prod_{j=1}^{m} (h\circ \sigma^{jN}) \, d\nu \\
& \leq & \rho_3^{(1-\delta)m}R^{\delta m} + R^{m} \nu(W) 
\leq \rho_4^m + \frac{2C_0 \, R^m}{|b|^{c_3}} \\
& \leq & e^{(L \log |b|) |\log \rho_4|} + \frac{2C_0 \, e^{2a_0 NT L \log |b|}}{|b|^{c_3}}
 =  \frac{1}{|b|^{L |\log \rho_4|}} + \frac{2C_0}{|b|^{c_3/2}} \leq \frac{C_5}{|b|^{4s}}
\end{eqnarray*}
for some global constant $C_5 > 0$, where as in Sect. 6.1, $s = 2/\beta_2$,  $c_3 = 16 s$, and by (7.3),
$L |\log \rho_4| > 8s$.
This completes the proof of the lemma.
\endofproof

\bs

\noindent
{\it Proof of Theorem } 1.3. We will again assume
that $f \in \ff_{\theta_1}(\hU)$; the general case $f \in \ff_\theta(\hU)$
will be done later using an approximation procedure.

Let $\hat{\theta} \leq \theta < 1$, where $\hat{\theta}$ is as in (5.1), and
let $N \geq N_0$. Let $L$ be as in Sect. 6.1. Choose $a_0 = a_0(N)$, $b_0 = b_0(N)$, $\rho_4 = \rho_4(N) \in (0,1)$, 
$C_5 > 0$ and $c_5 > 0$ as in Lemmas 7.1 and 7.2. Take $\theta_1 = \theta_1(\theta) \in (0,\theta]$ and  
$\theta_2 = \theta_2(\theta)\in [\theta,1)$  as in Sect. 6.3. Recall the set $\kk_b$ of pairs $(h,H)$ 
from Sect. 6.4. For any $h\in \ff_\theta(\hU)$ we set $\|h\|_\theta = \|h\|_0 + |h|_\theta$.

Let $|a|\leq a_0$ and  $|b|\geq b_0$, and let $h\in \ff_{\theta_1}(\hU)$ be such that 
$\| h\|_{\theta_1,b} \leq 1$.  Then  $|h(u)| \leq 1$ for all $u\in \hU$ and  $|h|_{\theta_1} \leq |b|$.  

Define $\hb$ by (6.5).
Assume that the points $u,u'$, the cylinder $\cc$ in $U$, the integer $p \geq 0$ and the points $v,v'\in \hU_1$
satisfy (6.17)  for some $i = 1,2$. Then, using (6.16)  and $|h|_{\theta_1} \leq |b|$ we get
\begin{eqnarray*}
|h(v) - h(v')| 
& \leq & |b|\, D_{\theta_1} (v,v') = |b|\, \theta_1^{p+N}\, D_{\theta_1} (u,u')
\leq  |b|\, \theta_1^{p+N}\, \diam_{\theta_1}(\cc) \\
& \leq & |b|\, \theta_1^{p+N}\, C_2 \, \diam(\Psi(\cc)) \leq E  |b|\, \theta_2^{p+N} \, \diam(\Psi(\cc)) ,
\end{eqnarray*}
since $C_2 \leq E$. Thus, $(h, 1) \in \kk_b$.
Set $h^{(m)} = \lab^{m N} h$ for $m \geq 0$. Define the sequence of functions $\{ H^{(m)}\}$ recursively by
$H^{(0)} = 1$ and $H^{(m+1)} = \nn_{J_m} H^{(m)}$, where $J_m \in \J(b)$ is chosen by induction as follows.
Since $(h^{(0)}, H^{(0)}) \in \kk_b$,  using Lemma 6.9  we find $J_0 \in \J(b)$ such that
for $h^{(1)} = \lab^{N} h^{(0)}$ and $H^{(1)} = \nn_{J_0} H^{(0)}$ we have
$(h^{(1)} , H^{(1)}) \in \kk_b$. Continuing in this way we construct by induction
an infinite sequence of functions $\{ H^{(m)}\}$ with
$H^{(0)} = 1$, $H^{(m+1)} = \nn_{J_m} H^{(m)}$ for all $m \geq 0$, such that $(h^{(m)} , H^{(m)}) \in \kk_b$.

Now set $m = L \hb$. Then Lemma 7.2 implies
$\int_U (H^{(m)})^2 d\nu \leq \frac{C_5}{|b|^{8s}} .$
Hence
$$\int_{U} |\lab^{m N} h|^2 \; d\nu  = \int_{U} |h^{(m)}|^2\; d\nu \leq \int_{U}  (H^{(m)})^2\; d\nu  
\leq \frac{C_5}{|b|^{8s}}.$$
From this it follows that for any $h \in \ff_{\theta_1}(\hU)$ we have 
$\int_{U} |\lab^{m N} h|^2 \; d\nu \leq  \frac{C_5}{|b|^{8s}} \, \|h\|^2_{\theta_1, b}$, and so
\be
\int_{U} |\lab^{m N} h| \; d\nu \leq  \frac{C_6}{|b|^{4s}} \,\|h\|_{\theta_1,b} .
\ee

\def\E{\mathcal E}
\def\tp{\tilde{p}}
\def\tq{\tilde{q}}

Next, we apply an approximation procedure to deal with functions 
$h \in \ff_\theta(\hU)$. Fix an arbitrary $\ep > 0$.
Assume that $\|h\|_{\theta,b} \leq 1$; then 
$\|h\|_0 \leq 1$ and $|h|_\theta \leq |b|$. So, using Lemma 5.2 with $H = 1$, 
we get 
\be
|\lab^{r} h|_{\theta} \leq A_0 [ |b| \theta^{r} + |b| ]\leq 2A_0 |b|
\ee
for any integer $r \geq 0$.

Recall from Sect. 6.1 that $\beta_2 > 0$ satisfies $\theta = \theta_1^{\beta_2}$.
Take the smallest integer $p$ so that $\theta^p \leq 1/|b|^{2}$. It is known 
(see e.g. the end of Ch. 1 in \cite{PP}) that there exists 
$h' \in \ff_{\theta_1}(\hU)$ which is constant on cylinders of 
length $p$ so that $\|h - h'\|_0 \leq |h|_\theta \, \theta^p$. Then 
$\|h-h'\|_0  \leq \frac{1}{|b|}$ and so $\|h'\|_0 \leq 2$, and it follows easily from this that 
$|h'|_{\theta_1} \leq \frac{4}{\theta_1^{p}} \leq \frac{4}{\theta^{p/\beta_2}} \leq C'_7 |b|^s ,$
where $s = \frac{2}{\beta_2} > 0$ and $C'_7 = 4/ \theta^{1/\beta_2}$. This and (7.8) imply
\be
\di \int_{U} |\lab^{m N} h' | \; d\nu \leq  \frac{2C'_7}{|b|^{3s}}  .
\ee
Moreover, as in (7.9) for $h$, we get 
\be
|\lab^{r} h'|_{\theta_1} \leq A_0 [ C'_7 |b|^s \theta_1^{r} + 2|b| ]\leq E |b|^s
\ee
for any integer $r \geq 0$, assuming $A_0( C'_7 + 2) \leq E$.

Next, recall from the Perron-Ruelle-Frobenius Theorem (see e.g. \cite{PP}) 
that there exist global constants $C_7 \geq 1$ and $\rho_5 \in (0,1)$ such that 
\be
\| L_{\f0}^n w - h_0 \, \int_{U} w \, d\nu \| \leq C_7\, \rho_5^n \, \|w\|_\theta
\ee
for all $w \in \ff_\theta(\hU)$ and all integers $n \geq 0$. 
The same estimate holds (we will assume with the same constants $C_7$
and $\rho_5$) for $h \in \ff_{\theta_1}$ replacing $h_0$ with the 
corresponding eigenfunction $h'_0$ and $\|w\|_\theta$ by $\|w\|_{\theta_1}$.

Write $\rho_5 = e^{-\beta_3}$ for some $\beta_3 > 0$, and assume (for later use) that
$LN\beta_3 > 4s = \frac{8}{\beta_2} .$
We have
\begin{eqnarray*}
|\lab^{2mN} h'|& =    & |\lab^{mN} \left(|\lab^{mN} h' |\right) \leq \ma^{mn} |\lab^{mN} h'| 
= L_{\f0}^{mN} \left(e^{\fa_{mN} - \f0_{mN}}\, |\lab^{mN} h'|\right)\\
& \leq & \left(L_{\f0}^{mN}\left(e^{\fa_{mN} -\f0_{mN}}\right)^{2}\right)^{1/2}\,
\left( L_{\f0}^{mN} |\lab^{mN} h'|^{2}|\right)^{1/2} .
\end{eqnarray*}
For the first term in this product (5.3) implies 
$$\left(L_{\f0}^{mN} \left(e^{\fa_{mN} - \f0_{mN}})\right)^{2}\right)^{1/2} \leq e^{Ta_0 mN} 
\leq  e^{T a_0 NL \log |b|} = |b|^{Ta_0 NL} < |b|^s ,$$
assuming 
$Ta_0 LN < s = \frac{2}{\beta_2} .$

For the second term, using (7.12) with $w = |\lab^{mN} h'|$, we get
$$L_{\f0}^{mN} |\lab^{mN} h'|^2 \leq L_{\f0}^{mN} |\lab^{mN} h'| \leq \|h'_0\|\, \int_U |\lab^{mN} h'|\, d\nu 
+ C_7\, \rho_5^{mN} \, \|\lab^{mN} h'\|_{\theta_1}.$$
By (7.11), $\|\lab^{mN} h'\|_{\theta_1} \leq E|b|^s$, so by (7.8),
$$L_{\f0}^{mN} |\lab^{mN} h'|^2 \leq \frac{4C'_7}{|b|^{3s}} + 2 E C_7 |b|^s \rho_5^{mN} .$$
Now 
$\rho_5^{mN} \leq e^{-\beta_3 N L \log |b|} = \frac{1}{|b|^{NL\beta_3}} < \frac{1}{|b|^{4s}} ,$
so we get
$L_{\f0}^{mN} |\lab^{mN} h'|^2 \leq \frac{C'_8}{|b|^{3s}} .$
Combining the estimates of the two terms, we get
$$|\lab^{2mN} h'| \leq  |b|^s (C'_8/|b|^{3s})^{1/2} \leq \frac{C_8}{|b|^{s/2}} 
= \frac{C_8}{|b|^{1/\beta_2}} < \frac{C_8}{|b|},$$
since $0 < \beta_2 < 1$.
Thus,
$$\di \|\lab^{2mN } h\|_0 \leq \|\lab^{2mN} h'\|_0 + \|\lab^{2mN} (h-h')\|_0 
\leq \frac{C_8}{|b|}  + \frac{1}{|b|} = \frac{C'_9}{|b|} .$$

Next, using (7.9) and Lemma 5.2 with $B = 2A_0 |b|$ and $H = 1$, we get
\begin{eqnarray*}
|\lab^{3mN} h|_\theta 
&    =   & |\lab^{mN} ( \lab^{2mN} h) |_\theta \\
& \leq  & A_0 \left[ 2A_0 |b| \, \theta^{mN} 
+ |b| \, \| \lab^{2mN} h\|_0 \right] \leq C''_9 .
\end{eqnarray*}
Hence
$ \|\lab^{3mN} h\|_{\theta,b} \leq \frac{C_9}{|b|}  \leq \frac{1}{|b|^{3/4}},$
assuming $|b| \geq b_0$ and $b_0 = b_0(\ep)$ is sufficiently large. This gives
$$\|\lab^{3mN} h\|_{\theta,b} \leq \frac{1}{|b|^{3/4}}\, \|h\|_{\theta,b}$$
for all $h \in \ff_\theta(\hU)$. 

Let $n \geq 3mN$ be an arbitrary integer. Writing $n = r(3mN)+ k$ for some
$k = 0,1, \ldots, 3mN-1$, and using the above $r$ times we get
$\|\lab^{r3mN} h\|_{\theta,b} \leq \frac{1}{|b|^{3r/4}}\, \|h\|_{\theta,b} .$
As in (7.9), using Lemma 5.2 with $H = 1$ and $B = |\lab^{r3mN} h|_\theta$, implies
\begin{eqnarray*}
|\lab^{n} h|_\theta = |\lab^k ( \lab^{r3mN} h) |_\theta 
 \leq   A_0 \left[ |\lab^{r3mN} h|_\theta \, \theta^{k} + |b| \, \| \lab^{r3mN} h\|_0 \right] ,
\end{eqnarray*}
so 
$$\frac{1}{|b|} |\lab^{n} h|_\theta \leq 2A_0 \|\lab^{r3mN} h\|_{\theta,b} \leq \frac{2A_0}{|b|^{3r/4}}\, \|h\|_{\theta,b} .$$
This and $\|\lab^{n} h\|_0 \leq \|\lab^{r3mN} h|\_0 \leq \frac{1}{|b|^{3r/4}}\, \|h\|_{\theta,b}$ give
\begin{eqnarray*}
\|\lab^{n} h\|_{\theta,b} 
& \leq & \frac{3A_0}{|b|^{3r/4}}\, \|h\|_{\theta,b} = 3A_0 e^{-(3r/4) \log |b|} \|h\|_{\theta,b} .
\end{eqnarray*}
We have $3r/4 \geq (r+1)/4$ for all $r \geq 1$, so the above implies
\begin{eqnarray} 
\|\lab^{n} h\|_{\theta,b}  \leq  3A_0  e^{-\frac{(r+1) \log |b|}{4}} \|h\|_{\theta,b}
 \leq  3A_0 e^{-\frac{(r+1)3mN}{12LN}} \|h\|_{\theta,b} 
& \leq& 3A_0  \rho_6^n \|h\|_{\theta,b} ,
\end{eqnarray}
where $\rho_6 = e^{-1/(12LN)} \in (0,1)$. 

Thus, (7.13) holds for all $h \in \ff_\theta(\hU)$ and all integers $n \geq 3mN = 3LN \lceil\log |b|\rceil$. 
Finally, recall the eigenfunction $h_a \in \ff_\theta(\hU)$ for the operator $L_{f- (P_f+a)\tau}$
from Sect. 5.1. It is known that $\|h_a\|_\theta \leq \Con$ for bounded $a$, e.g. for $|a|\leq a_0$.
It now follows from
$$\lab^n(h/h_a) = \frac{1}{\lambda_a^n h_a}\, L^n_{f- (P+a+\i b) \tau} h$$
and the above estimate that there exist constants $0 < \rho < 1$, $a_0 > 0$, 
$b_0 \geq 1$ and  $C > 0$ such that if $a,b\in \R$  satisfy $|a| \leq a_0$ and $|b| \geq b_0$, 
then 
\be
\|L_{f -(P_f+a+ \i b)\tau}^n h \|_{\theta,b} \leq C \;\rho^n \; \| h\|_{\theta,b}
\ee
for any integer $n \geq 3LN \log |b|$ and any  $h\in \ff_\theta (\hU)$.

\bs

This completes the proof of Theorem 1.3 under the assumption that $f \in \ff_{\theta_1}(\hU)$. The case $f \in \ff_{\theta}(\hU)$ 
follows by using an approximation procedure. To our knowledge this has not been done anywhere in details, and the argument 
involved is  not trivial, so we will sketch it for completeness.

\bs

\noindent
{\it Sketch of the proof of Theorem 1.3 for arbitrary $f \in \ff_\theta(\hU)$:}\\
We will use again the constants from the beginning of Sect. 6.1, including $\theta_1$, $\theta_2$, etc.
Fix $B$, $L$, $N$, $c'$, $c_4$ as before and define $\hb$ by (\ref{eq:hbdef}). 
Let $|a| \leq a_0$ and $|b| \geq b_0$, where $b_0$ is given by (\ref{eq:b0cond}). 

Let $f \in \ff_\theta(\hU)$ be an arbitrary real-valued function. Take an integer $t = t(b)  > 0$ so that
\be
\theta^t = \left\lceil \frac{8A_0}{\log|b|} \right\rceil ,
\ee
where $A_0$ is the constant from Lemma 5.2. There exists a real-valued $\ftt$ depending only on $t$ coordinates
such that 
$$\|f- \ftt\|_0 \leq |f|_\theta \, \theta^t \leq T\, \theta^t $$
(see the end of Ch. 1 in \cite{PP}), where $T$ is as in (5.3). Then $\ftt \in \ff_{\theta_1}(\hU)$,
$\|\ftt\|_0 \leq 2T$ and 
$$|\ftt|_\theta \leq \frac{4}{\theta^t} \leq \frac{\log |b|}{2A_0} \quad, \quad |\ftt|_{\theta_1} \leq \frac{4}{\theta_1^t}.$$

Let $\lambda_{at}$ be the {\it largest eigenvalue} of
$\Fat =  \ftt - (P_t + a)\tau,$
where $P_t = P_{\ftt}$, and let $h_{at} \in \ff_\theta(\hU)$ be a corresponding (positive) eigenfunction such that 
$\int h_{at} \, d\hnu_{at} = 1$, where $\hnu_{at}$ is the unique regular probability measure on $\hU$ with 
$(\Fat)^* \hnu_{at} = \hnu_{at}$. 


\def\nnt{\nn^{(t)}}

For $|a|\leq a'_0$,  as in \cite{D}, consider the function
$$\fat(u) = \ftt (u) - (P_t + a) \tau(u) + \ln h_{at}(u) -  \ln h_{at}(\sigma(u)) - \ln \lambda_{at}$$
and the operators 
$$ \labt = L_{\fat - \i\,b\tau} : \ff_\theta(\hU) \longrightarrow \ff_\theta (\hU)\:\:\: , 
\:\:\: \mat = L_{\fat} : \ff_\theta(\hU) \longrightarrow \ff_\theta(\hU) .$$
Then $\mat \; 1 = 1$ and $\di |(\labt^m h)(u)| \leq (\mat^m |h|)(u)$ for all $u\in \hU$.

Using part of the proof of Lemma 4.1 in \cite{PeS5}, one shows that $|h_{at}|_\theta \leq \Con\, |\ftt|_\theta$ for some global constant
$\Con > 0$. Thus, $|\fat|_\theta \leq \Con\, |\ftt|_\theta$, and is also clear that $\|\fat\|_0 \leq \Con$.

Next, define the set $K_0$, cylinders $\cc_m$ and their sub-cylinders $\dd_j$ and the function
$\omega = \omega_J$ as in Sect. 6.1 and consider the operator $\nnt = \nnt_J$ on $\ff_\theta(\hU)$
defined by
$$\nnt(h) = \mat^N (h) = L^N_{\fat} h .$$
It is important to notice that 
\be
e^{2|\ftt|_\theta} \; \diam_\theta(\cc_m) \leq \frac{1}{|b|^{1/A_0} },
\ee
provided we took the constants $A_0$ in Lemma 5.2 and  $B$ in Sect. 6.1 so that $A_0 \geq \frac{2B}{|\log \theta|}$.
Indeed, for the length $\ell_m$ of $\cc_m$ we have (\ref{eq:Bcond}), so
\begin{eqnarray*}
e^{2|\ftt|_\theta}\; \diam_\theta(\cc_m) 
& \leq &  e^{\frac{\log |b|}{A_0}} \theta^{\ell_m} =
|b|^{1/A_0} e^{-\ell_m |\log\theta| } \leq  |b|^{1/A_0} e^{-(|\log\theta|/B) \log|b| }\\
& =    & |b|^{1/A_0 - \frac{1}{B} |\log \theta|} \leq |b|^{-1/A_0} ,
\end{eqnarray*}
which proves (7.16).


Then we define the metric $\dd(u,u')$ on $\hU$ and the class of positive functions $\kk_E$ as in Sect. 6.2.
Now with the above one easily shows that Lemma 6.5 is valid in the form $\nnt (\kk_E) \subset \kk_E$.
Indeed, the main observation to make to prove this is that, given $u, u' \in \hU$ such that  there exists
an integer $p \geq 0$ with $\sigma^p(Y(u, u')) \subset \cc_m$ for some $m \leq m_0$ and $\ell(u,u') \geq p$,
then for any integer $k \geq 1$, if $v, v'(v) \in \hU$ satisfy $\sigma^{k}(v) = u$, $\sigma^k(v') = u'$
and belong to the same cylinder of length $k$, then
\begin{eqnarray}
|\fat_k(v) - \fat_k(v')| 
& \leq & \sum_{j=0}^{m-1} |\fat|_\theta \,\theta^{m-j}\, \dte (u,u') \leq \Con\; |\ftt|_\theta\,D_\theta(u,u') \nonumber \\
& \leq & \Con\; |\ftt|_\theta\, \diam_\theta(\cc_m) \leq \Con .
\end{eqnarray}
With this observation, a simple modification of the proof of Lemma 6.5 gives $\nnt H \in \kk_E$ for
every $H \in \kk_E$.

Next, we define the class of functions $\kk_0$ as in Sect. 6.3 and prove the analogue of Lemma 6.7:
$L_{\ftt - s \tau}(\kk_0) \subset \kk_0$ for all $s$ with $|s| \leq M_0$ and all integers $q\geq N$.
To prove this, the choice of $\theta_1$ is important; it implies 
$$\diam_{\theta_1}(\cc_m) = \theta_1^{\ell_m} \leq \theta_1^{\ell_m/2} \; (\theta^{\ell_m})^{1/\alpha_2}
\leq \theta_1^{\ell_m/2} \; \diam(\tcc_m) .$$
Then, assuming $u, u', v, v', w, w'$ are as in the proof of Lemma 6.7, we derive
\begin{eqnarray}
|\fat_N(w) - \fat_N(w')| 
& \leq & \Con\, |\ftt|_{\theta_1} \theta_1^{p+N} \diam_{\theta_1}(\cc_m)
 \leq  \Con \; \frac{4}{\theta_1^t} \theta_1^{p+N} \, \theta_1^{\ell_m/2} \; \diam(\tcc_m) \nonumber\\
& \leq & \Con \, \theta_2^{p+N} \theta_1^{\ell_m/2 - t} \; \diam (\tcc_m) 
\leq  \Con \, \theta_2^{p+N} \; \diam (\tcc_m)  \leq 1,
\end{eqnarray}
since $t < < \ell_m/2$. Now the rest of the proof of Lemma 6.7 is the same, and as a consequence one
gets (as in Corollary 6.8) that the eigenfunctions $h_{at}$ belong to $ \kk_0$.

Finally, the arguments in Sect. 6.4 can be repeated with very little change -- the main one is that
in the first estimate of $|(\labt^Nh)(v) - (\labt^N)(v')|$ one has to use (7.17) again. This proves the
analogue of Lemma 6.9, where the operator $\lab^N$ is replaced by $\labt$.

We will now prove Lemma 6.9 in its original form under the present assumption that $f \in \ff_\theta(\hU)$.

\bs

\noindent
{\bf Lemma 7.3.}  {\it Assume $f\in \ff_\theta(\hU)$.
Choosing $E > 1$ and $\mu_0$ as in Sect. {\rm 5.2} and assuming $N$ is sufficiently large, 
for any $|a| \leq a'_0$, any $|b| \geq b_0$ and any $(h,H) \in \kk_b$ 
there exists $J \in \J(b)$ such that $(\lab^{N} h, \nn_J H ) \in \kk_b$.}

\bs

\noindent
{\it Proof.} Consider the function
$$\zeta = \fa_N - \fat_N \in \ff_\theta(\hU) .$$
Notice that for any $u \in \hU$ and any function $h$ on $\hU$ we have
$$(\ma^N h)(u) = \sum_{\sigma^Nv = u} e^{\fa_N(v)} h(v)
=   \sum_{\sigma^Nv = u} e^{\fat_N(v)} e^{(\fat_N - \fa_N)(v)} h(v)
= (\mat^N (e^{\zeta} h))  (u) .$$
Thus, $\ma^N h = \mat^N (e^\zeta h)$, and similarly one observes that
 $\lab^N h = \labt^N (e^\zeta h)$.
 
 We will now repeat the argument from the proof of Lemma 6.9.

Let $|a| \leq a'_0$, $|b| \geq b_0$ and $(h, H) \in \kk_b$. We will construct a representative set 
$J \in \J(b)$ such that  $(\lab^N h , \nn_J H) \in \kk_b$.  Given  an arbitrary representative set $J$,
we will first show that  $(\lab^N h , \nn_J H) = (\labt^n (e^\zeta h), \nnt_J (e^\zeta H)) $ has property (ii). 

Assume that the points $u,u'$, the cylinder $\cc_m$ in $U_1$, the integer $p \geq 0$ and the points 
$v,v'\in \hU$ satisfy (\ref{eq:u-cond})  for some $i = 1,2$ and $\ell = 1, \ldots, \ell_0$. 
Since $h_{at} \in \kk_0$, we have
$$|\ln h_a(w) - \ln h_a(w') | \leq \frac{| h_a(w) - h_a(w')|}{\min \{ |h_a(w)|, |h_a(w')|\}} 
\leq E_1 \,\theta_2^{p+2N} \, \diam(\tcc_m) , $$
and similarly, $|\ln h_{at}(v) - \ln h_{at}(v') | \leq E_1 \,\theta_2^{p+N} \, \diam(\tcc_m)$.
By (7.18),
\be
|\fat_N(w) - \fat_N(w')| \leq 
\Con\, \theta_2^{p+N}\, \diam(\tcc_m)  \leq 1,
\ee
assuming $N$ is chosen appropriately. Thus, using (7.19) and a 
similar but simpler estimate for $| \fa_N(w) - \fa_N (w')|$, we get
\begin{eqnarray}
|\zeta(w) - \zeta (w')|
& =    & | (\fa - \fat)_N(w) - (\fa - \fat)_N (w')|\nonumber\\
& \leq &  C_{10} \, \theta_2^{p+N} \; \diam (\tcc_m) < C_{10}
\end{eqnarray}
for some global constant $C_{10} > 0$.
This implies
\begin{eqnarray}
|e^{\zeta (w) - \zeta (w')} - 1|
& \leq & \Con\, \left| \zeta (w) - \zeta (w') \right|  \leq  C_{11} \, \theta_2^{p+N} \; \diam (\tcc_m) < C_{11}
\end{eqnarray}
for some global constant $C_{11} > 0$.

Now using (7.19), (7.20) and (7.21), as in the proof of Lemma 6.9,
 for any $a$ and $b$ with $|a| \leq a'_0$ and $|b| \geq b_0$   we derive
\begin{eqnarray*}
&        &|(\lab^N h)(v) - (\lab^N  h)(v')| = |(\labt^N (e^\zeta h))(v) - (\labt^N  (e^\zeta h))(v') |\\
& \leq & [2 e^{1+ C_{10}}\, E \theta_2^N  + 2 C_{12} ]  \,|b|\, \theta_2^{p+N}\,    \diam(\tcc_m)\, \nnt_J (e^\zeta H)(v') 
\end{eqnarray*}
for some global constant $C_{12} > 0$.
Assuming $2 e^{1+C_{10}} \theta^N \leq 1/2$ and $2 C_{12}  \leq E/2$, we get
$$|(\lab^N h)(v) - (\lab^N  h)(v')|
\leq  E \,|b|\, \theta_2^{p+N}\,    \diam(\tcc_m)\, \nnt_J (e^\zeta H)(v') 
= E \,|b|\, \theta_2^{p+N}\,    \diam(\tcc_m)\, (\nn_J H)(v') ,
$$
so, $(\lab^N h , \nn_J H)$ has property (ii).

Now  we will construct $J$ so that  $|\lab^N h|(u) \leq (\nn_J H) (u)$ for all $u \in \hU$, which is equivalent to
\be
|\labt^N (e^\zeta h)|(u) \leq (\nnt_J (e^\zeta H)) (u)
\ee
for all $u \in \hU$. 

Define the functions $\psi_\ell, \gao_\ell, \gat_\ell : \hU  \longrightarrow \C$  as in the proof of Lemma 6.9.
Notice that
$$\di \psi_\ell(u) = e^{(\fat_{N}+\i b\tau_{N})(\vl_1(u))} (e^{\zeta} h) (\vl_1(u)) 
+ e^{(\fat_{N}+\i b\tau_{N})(\vl_2(u))} (e^\zeta h) (\vl_2(u)) ,$$
$$\di \gao_\ell(u) = (1-\mu_0)\, e^{\fat_{N} (\vl_1(u))} (e^\zeta H) (\vl_1(u)) + e^{\fat_{N}(\vl_2(u))} (e^\zeta H) (\vl_2(u)) ,$$
and similarly for $\gat_\ell(u)$.
We will use again the functions $\varphi_\ell(u) = \varphi_\ell(\hz_0, u)$, $u \in U$,  from Sect. 5.3. 

As before (7.22) is trivially satisfied for $u \notin V_b$ for any choice of $J$.

Consider an arbitrary  $m = 1, \ldots, m_0$. We will construct  $j \leq j_0$ with $\dd_{j} \subset \cc_m$,  
and a pair  $(i, \ell)$ for which $(i,j,\ell)$ will be included in $J$.

\ms

\noindent
{\bf Case 1.}
There exist $j\leq j_0$ with $\dd_j \subset \cc_m$, $i = 1,2$ and $\ell \leq \ell_0$   such that the 
first alternative in Lemma 6.10(b) holds for $\hdd_j$,  $i$ and $\ell$. This case is dealt with exactly
as in the proof of Lemma 6.9.
 

\medskip

\noindent
{\bf Case 2.}
For all $j\leq j_0$ with $\dd_j \subset \cc_m$, $i = 1,2$ and $\ell \leq \ell_0$ the 
second alternative in Lemma 6.10(b) holds for $\hdd_j$,  $i$ and $\ell$, i.e.
\be
|h(\vl_i(u))|\geq \frac{1}{4}\, H(\vl_i(u)) > 0
\ee
for any  $u \in \hcc'_m$.

Let $u,u' \in \hcc'_m$, and let $i = 1,2.$ 
Using (6.20) and the assumption that $(h,H) \in \kk_b$, and in particular property 
(ii) with $p = 0$, $v = \vl_i(u)$ and $v' = \vl_i(u')$, and also (7.20) and (7.21) with $p =0$, 
and assuming e.g.  $\min\{ |e^{\zeta(\vl_i(u))} h(\vl_i(u))| , |e^{\zeta(\vl_i(u'))} h(\vl_i(u'))|\}  = |e^{\zeta(\vl_i(u'))} h(\vl_i(u'))|$, we get
\begin{eqnarray*}
&        & \frac{|e^{\zeta(\vl_i(u))} h(\vl_i(u)) - e^{\zeta(\vl_i(u'))} h(\vl_i(u'))|}{\min\{ |e^{\zeta(\vl_i(u))} h(\vl_i(u))| , |e^{\zeta(\vl_i(u'))} h(\vl_i(u'))| \}}\\
& \leq &   \frac{|e^{\zeta(\vl_i(u))}  - e^{\zeta(\vl_i(u'))}|\; | h(\vl_i(u'))|}{|e^{\zeta(\vl_i(u'))} h(\vl_i(u'))|}
+  \frac{e^{\zeta(\vl_i(u))} | h(\vl_i(u)) -  h(\vl_i(u'))|}{ |e^{\zeta(\vl_i(u'))} h(\vl_i(u'))| }\\
& \leq &   |e^{\zeta(\vl_i(u))  - \zeta(\vl_i(u'))} - 1|
+  C_{13} \frac{ | h(\vl_i(u)) -  h(\vl_i(u'))|}{ | h(\vl_i(u'))| }\\
& \leq & C_{13} \theta_2^N \diam(\tcc_m) + \frac{E|b|\,\theta_2^N H(\vl_i(u'))}{|h(\vl_i(u'))| } \diam (\tcc_m) 
  \leq   (C_{13} + 4 E|b|) \,   \theta_2^{N} \frac{\ep_1}{|b|} < 5  E \theta_2^N \ep_1 ,
\end{eqnarray*}
assuming $E\geq C_{13}$. So, the angle between the vectors
$e^{\zeta(\vl_i(u)} h(\vl_i(u))$ and $e^{\zeta(\vl_i(u')} h(\vl_i(u'))$ in $\R^2$  is  
$< 10 E \theta_2^N \ep_1 < \pi/6$ by (\ref{eq:Ncond}).  

Since $e^{\zeta(\vl_i(u)} $ and $e^(\zeta{\vl_i(u)} $ are real numbers, the arguments of the complex numbers\\
$e^{\zeta(\vl_i(u)} h(\vl_i(u))$ and $e^{\zeta(\vl_i(u')} h(\vl_i(u'))$ are the same as those of 
$h(\vl_i(u))$ and $h(\vl_i(u'))$. 

As before,  for each $i = 1,2$ we can choose a real continuous
function $\thetam_i(u)$, $u \in  \cc'_m$, with values in $[0,\pi/6]$  and a constant $\lambdam_i$ such that
$$\di h(\vl_i(u)) = e^{\i(\lambdam_i + \thetam_i(u))}|h(\vl_i(u))| \quad , \quad u\in \cc'_m .$$
Fix an arbitrary $u_0\in \cc'_m$ and set $\lambdam = |b| \varphi_\ell(u_0)$, and assume again 
 that $|\lambdam_2 - \lambdam_1 + \lambdam | \leq \pi$. Then
$|\thetam_i(u) - \thetam_i(u')|\leq 2 \sin |\thetam_i(u) - \thetam_i(u')| < 16 E \theta_2^N \ep_1$ for all
$u,u' \in \cc'_m$.

As in the proof of Lemma 6.9, 
the difference between the arguments of the complex numbers
$e^{\i \,b\,\tau_N(\vl_1(u))} (e^\zeta h) (\vl_1(u))$ and $e^{\i \,b\, \tau_N(\vl_2(u))} (e^\zeta h) (\vl_2(u))$
is given by the function
\begin{eqnarray*}
\Gl(u) 
& = & [b\,\tau_N(\vl_2(u)) + \thetam_2(u) + \lambdam_2] -  
[b\, \tau_N(\vl_1(u)) + \thetam_1(u) + \lambdam_1]\\
& = & (\lambdam_2-\lambdam_1) + |b| \varphi_\ell(u) + (\thetam_2(u) - \thetam_1(u)) ,
\end{eqnarray*}
and as before we prove that there exist $j \leq j_0$ and $\ell \leq \ell_0$ such that
$\ep_3 \leq |\Gl(u)| <  \frac{3\pi}{2}$ for all $u \in \dd'_j$.

\def\halpha{\hat{\alpha}}
\def\hgamma{\hat{\gamma}}

Now, exactly as in the proof of Lemma 6.9, following arguments from \cite{D} (see also \cite{St2})
we prove that, assuming 
$\mu_0 \leq t$, we have $|\psi_\ell(u)| \leq \gamma^{(1)}_\ell (u)$  
for all $u \in \dd'_j$.  Now set $j_m = j$, $\ell_m = \ell$ and $i_m = 1$,
and include  $(i_m,j_m,\ell_m)$ in the set $J$.  
Then $\dd_{j_m} \subset \cc_m$ and we deduce that  (7.22) holds on $\dd'_{j_m}$.
\endofproof

\bs

Next, we proceed with what is done in Sect. 7. First, we prove parts (a) and (b)  of Lemma 7.1 assuming $f \in \ff_\theta(\hU)$.
Part (a) goes without a change. In part (b) one proves that 
\be
\int_{V_b} (\nn_J H)^2 \, d\nu \leq \rho_3\, \int_{V_b} L^N_{\f0} (H^2) \, d\nu ,
\ee
for any $H\in \kk_E$  and any $J\in \J(b)$,  where $\di \rho_3 = \rho_3(N) < 1$ is possibly a slightly larger constant, and again
$a_0 = a_0(N) > 0$ is chosen sufficiently small. The proof of this uses the same lines as the ones in the proof of Lemma 7.1(b)
combined with the fact that $\|\zeta_N\|_0 = \|\fa_N - \fat_N\| \leq N \, \Con$ for some global constant $\Con > 0$. 

Then, using the analogue of Lemma 7.1 (with $f \in \ff_\theta(\hU)$) and  Lemma 7.3 one proves Lemma 7.2 in the same form -- the 
difference is that now $f \in \ff_\theta(\hU)$ compared to the previous stronger assumption $f \in \ff_{\theta_1}(\hU)$. This gives
the estimate (7.4) in exactly the same form under this more general assumption.  And then one just needs to repeat the argument
from the proof of Theorem 1.3 (the same as under the assumption $f \in \ff_{\theta_1}(\hU)$). 
\endofproof

\bigskip

\noindent
{\it Proof of Theorem} 1.1. This follows from the procedure described in \cite{D} (see Sect. 4 and Appendix 1  there).
\endofproof

\bs

\noindent
{\it Proof of Corollary 1.4.} Let again $\hat{\theta}$ be as in (5.1). Given $\ep > 0$,
choose the constants $C > 0$, $\rho \in (0,1)$, $a_0 > 0$ and $b_0 \geq 1$ as in Theorem 1.3.
Let $\hat{\theta} \leq \theta < 1$. As in the proof of Lemma 5.1,
$(d(x,y))^{\alpha} \leq \Con \dte (x,y)$ will always hold assuming $1/\gamma^{\alpha} \leq \theta$, 
i.e. $\alpha \geq \frac{|\log\theta|}{\log \gamma}$. Here 
$1 < \gamma < \gamma_1$ are the constants from (2.1).
Then for such $\alpha$ we have $|h|_\theta \leq \Con |h|_{\alpha}$. 

Set $\alpha_0 = \frac{|\log \hat{\theta}|}{\log \gamma} > 0$. 
Let again $\alpha_1 \in (0,1]$ be such  that the local stable holonomy maps on $\tR$ are uniformly 
$\alpha_1$-H\"older, i.e. there exists a constant $C_{11} > 0$ such that for any $z,z' \in \tR_i$ 
for some $i = 1, \ldots, k_0$ and any 
$x,y \in W^u_{\tR}(z)$ for the projections $x',y'\in W^u_{\tR}(z')$ of $x, y$ along stable leaves we have
$d(x',y') \leq C_{11} \, (d(x,y))^{\alpha_1}$.

Let $\alpha \in (0,\alpha_0]$; then  $\alpha = \frac{|\log \theta|}{\log \gamma} $
for some $\theta \in [\hat{\theta},1)$. As above this gives $|h|_\theta \leq C'_{12} |h|_\alpha$ for any
$h \in C^\alpha(\hU)$. 

Assume that for a given $h \in C^\alpha(\hU)$ we have $\|h\|_{\alpha,b} \leq 1$; then $\|h\|_0 \leq 1$
and $|h|_\alpha \leq |b|$, so $|h|_\theta \leq C'_{12}|b|$ and therefore $\|h\|_{\theta,b} \leq C'_{12}+1$.
By Theorem 1.3, 
$$\|\lab^{n} h\|_{\theta,b} \leq 2C'_{12} C_{10} |b|^\ep \rho_{13}^n \quad, \quad n \geq 0 ,$$
so in particular 
\be
|\lab^{n} h|_0 \leq C_{12} |b|^\ep \rho_{13}^n
\ee
for all $n \geq 0$.

Next, one needs to repeat part of the arguments from the proof of Theorem 1.3 above. 

First, one needs a version of Lemma 5.2(b) for functions $w\in C^\alpha(\hU)$. Given an integer
$m \geq 0$ and $u,u' \in U_i$ for some $i = 1, \ldots,k_0$, notice that of $\sigma^m(v) = u$,
$\sigma^m(v') = u'$ and $v' = v'(v)$ belongs to the cylinder of length $m$ containing $v$, then
\begin{eqnarray}
|w(\sigma^jv) - w(\sigma^j(v')| 
& \leq & |w|_\alpha (d(\sigma^jv), \sigma^j(v'))^\alpha 
 \leq  C'_{13} \frac{|w|_\alpha}{\gamma^{\alpha(m-j)}} (d(u,u'))^{\alpha \alpha_1} .
\end{eqnarray}
This implies 
\be
|w_m(v) - w_m(v')| \leq C_{13}'' |w|_\alpha (d(u,u'))^{\alpha \alpha_1} .
\ee
This is true for $w = f$, $w = \tau$. Now repeating the argument in the proof of Lemma 5.2(b), 
for $|a| \leq a_0$ and $w \in C^\alpha(\hU)$ we get
\begin{eqnarray*}
        |(L^m_{f-(P+a)\tau} w)(u) - (L^m_{f- (P+a)\tau} w)(u')|
 \leq   C_{13} \left[ \frac{|w|_\alpha}{\gamma^{\alpha m}} 
   + \|w\|_0 \right]\, (L^m_{f- (P+a)\tau} 1)(u)\, (d(u,u'))^{\alpha \alpha_1} .
\end{eqnarray*}
In particular this shows that $L^m_{f-(P+a)\tau} w \in C^{\alpha \alpha_1}(\hU)$ for all $w\in C^\alpha(\hU)$
and all integers $m \geq 0$. Since $w = 1 \in C^\alpha(\hU)$, it now follows from Perron-Ruelle-Frobenius
Theorem that the eigenfunction $h_a \in C^{\alpha \alpha_1}(\hU)$ and so $\fa\in C^{\alpha\alpha_1}(\hU)$
for all $|a| \leq a_0$. Moreover, taking $a_0$ sufficiently small, we may assume that 
$\|h_a\|_{\alpha\alpha_1} \leq C'_{14} = \Con$ for all $|a| \leq a_0$. 
Using (7.27) with $w = f_a$ and $\alpha$ replaced by $\alpha\alpha_1$, we get
$|\fa_m(v) - \fa_m(v')| \leq C_{13}''  (d(u,u'))^{\alpha \alpha_1^2}, $
and also $|\fa(v) - \fa(v')| \leq C_{14}'' \, \rho_{14}^m (d(u,u'))^{\alpha \alpha_1^2}$.

Now, using standard arguments, for $h \in C^\alpha(\hU)$ we get
\be
|\lab^m h(u) - \lab^m h(u')| \leq C_{14} \left[\rho_{14}^m\, |h|_\alpha  
+ |b|\, \|h\|_0\right]\, (d(u,u'))^{\alpha \alpha_1^2} .
\ee
Since $|h|_\alpha \leq |b|$ and $\|h\|_0 \leq 1$, this gives 
$|\lab^mh|_{\alpha \alpha_1^2} \leq \Con\, |b|$ for all $m \geq 0$. Using (7.26) and (7.29) 
with $h$ replaced by $\lab^m h$ and $\alpha$ replaced by $\alpha \alpha_1^2 \leq \alpha_0$, we get
\begin{eqnarray*}
|(\lab^{2m} h) (u) - (\lab^{2m} h)(u')| 
& \leq & \Con \left[\rho_{14}^m\, |\lab^m h|_{\alpha \alpha_1^2} 
+ |b|\, \|\lab^m h\|_0\right]\, (d(u,u'))^{\alpha \alpha_1^4} \\
& \leq & C'_{15} \left[\rho_{14}^m\, |b| + |b|\,|b|^\ep \rho_{13}^m \right]\, (d(u,u'))^{\alpha \alpha_1^4} 
\end{eqnarray*}
Thus, $\|\lab^{2m} h\|_{\alpha \alpha_1^4,b} \leq C_{15} |b|^\ep \rho_{14}^m$ for all $m \geq 0$ and all
$h \in C^\alpha(\hU)$ with $\|h\|_\alpha \leq 1$.  
Since 
$$L^m_{f- (P_f+a + \i b)\tau} h = \frac{1}{h_a}\, \lab^m (h_a h) ,$$ 
it is now easy to get
$$\|L^{m}_{f-(P+a - \i b)} h\|_{\alpha \alpha_1^4,b} \leq C_{16} |b|^\ep \rho_{14}^m\, \|h\|_{\alpha,b}$$
for all $m \geq 0$ and all $h \in C^\alpha(\hU)$.  
Setting $\hbeta = \alpha_1^4$ proves the assertion.
\endofproof

\def\tyj{\tilde{y}^{(j)}}
\def\yl{y^{(\ell)}}
\def\wo{w^{(1)}}
\def\vo{v^{(1)}}
\def\vi{v^{(i)}}
\def\vj{v^{(j)}}
\def\vk{v^{(\tk)}}
\def\uo{u^{(1)}}
\def\wt{w^{(2)}}
\def\xio{\xi^{(1)}}
\def\xit{\xi^{(2)}}
\def\xii{\xi^{(i)}}
\def\xij{\xi^{(j)}}
\def\hxio{\hxi^{(1)}}
\def\hxit{\hxi^{(2)}}
\def\hxii{\hxi^{(i)}}
\def\hxij{\hxi^{(j)}}

\def\cxi{\check{\xi}}
\def\cxio{\cxi^{(1)}}
\def\cxit{\cxi^{(2)}}
\def\cet{\check{\eta}}
\def\ceto{\cet^{(1)}}
\def\cett{\cet^{(2)}}
\def\cv{\check{v}}
\def\cvo{\cv^{(1)}}
\def\cvt{\cv^{(2)}}
\def\cu{\check{u}}
\def\cuo{\cu^{(1)}}
\def\cut{\cu^{(2)}}
\def\cj{c^{(j)}}
\def\fj{f^{(j)}}
\def\gji{g^{(j,i)}}

\def\ao{a^{(1)}}
\def\bo{b^{(1)}}
\def\zi{z^{(i)}}
\def\hGa{\widehat{\Gamma}}

\def\tvo{\tv^{(1)}}
\def\tetao{\teta^{(1)}}



\section{Temporal distance estimates on cylinders}
\setcounter{equation}{0}

Here we prove Lemmas 4.3 and 4.4.

\subsection{A  technical lemma}

Notice that in Lemma 4.1 the exponential maps are used to parametrize $W^u_\ep(z)$ and $W^s_\ep(z)$.
The particular choice of the exponential maps is not  important, however it is important that these
maps are $C^2$. So, we cannot use the maps $\Phi^u_z$ and $\Phi^s_z$ defined in Sect. 3. 
In order to use Lemma 4.1 we will need in certain places to replace the  local lifts $\hf^p_z$ of the 
iterations $f^p$ of the map $f$ by slightly different maps.

For any $x\in \ll$ consider the $C^{2}$ map (assuming $r(x)$ is chosen small enough)
$$\tf_x = (\exp^u_{f(x)})^{-1} \circ f\circ \exp^u_x : E^u (x; r(x))  \longrightarrow E^u (f(x), \tr (f(x)))\;.$$
As with the maps $\hf$, for $y \in \ll$ and an integer $j \geq 1$ we will use the notation
$$\tf_y^j = \tf_{f^{j-1}(y)} \circ \ldots \circ \tf_{f(y)} \circ \tf_y\quad,
\quad \tf_y^{-j} = (\tf_{f^{-j}(y)})^{-1} \circ \ldots \circ (\tf_{f^{-2}(y)})^{-1} 
\circ (\tf_{f^{-1}(y)})^{-1} \;,$$
at any point where these sequences of maps are well-defined.
In a similar way one defines the maps $\tf_x$ and their iterations on $E^s(x;r(x))$.

Following the notation in Sect. 3 and using the fact that the flow $\phi_t$ is contact,
the negative Lyapunov exponents over $\ll$ are 
$-\log \lambda_1 > -\log \lambda_{2}  > \ldots > - \log \lambda_{\tk} .$
Fix $\hep > 0$ as in Sect. 3, assuming in addition that 
$$\hep \leq \frac{\log \lambda_1}{100}\, \min\{ \beta, \vartheta\} 
\quad , \quad \hep \leq \frac{\log \lambda_1 \, (\log \lambda_2 - 
\log \lambda_1)}{4 \log \lambda_1  + 2\log \lambda_2} .$$
For $x\in \ll$ we have an $f$-invariant decomposition
$E^s(x) = E^s_1(x) \oplus E^s_2(x) \oplus \ldots \oplus E^s_{\tk}(x)$
into subspaces of dimensions $n_1, \ldots, n_{\tk}$, where $E^s_i(x)$ $(x\in \ll$) is the
$df$-invariant subbundle corresponding to the Lyapunov exponent $- \log \lambda_i$.
For the {\it Lyapunov $\hep$-regularity function}  $R = R_{\hep} : \ll \longrightarrow (1,\infty)$,
chosen as in in Sect. 3 (see also Sect. 4), we have
\be
\frac{1}{R (x)\, e^{m\hep}} \leq \frac{\|df^m(x)\cdot v\|}{\lambda_i^{-m}\|v\|} 
\leq R (x)\, e^{m\hep} \quad , \quad x\in \ll  \;, \; v\in E^s_i(x)\setminus \{0\} \;, \; m \geq 0.
\ee

We will also assume that the set $P_0$ is as in (4.4), and the regularity functions
$R_{\hep}(x)$, $r(x)$, $\Gamma(x)$, $L(x)$, $D(x)$ satisfy (4.5). 

For the contact form $\omega$ it is known (see e.g. Sect. in \cite{KH} or Appendix B in
\cite{L1}) that $\omega$ vanishes on every stable/unstable manifold of a point on $M$,
while $d\omega$ vanishes on every weak stable/unstable manifold.
For Lyapunov regular points we get a bit of extra information.

\bs

\noindent
{\bf Lemma 8.1.} {\it For every $x\in \ll$ and every
$u = (\uo, \ldots, u^{(\tk)}) \in E^u(x;r(x))$ and $v = (\vo, \ldots,\vk)\in E^s(x;r(x))$
we have}
\be
d\omega_x(u,v) = \sum_{i=1}^{\tk} d\omega_x(\ui,\vi) .
\ee

\noindent
{\it Proof.} It is enough to show that $d\omega_x(\ui, \vj) = 0$ if $i \neq j$. Let e.g. $i < j$.  
Using (3.3), (4.6), (8.1) and  the fact that $d\omega$ is $df$-invariant, for $m \geq 0$ and $x_m = f^m(x)$ we get
\begin{eqnarray*}
 |d\omega_x(\ui,\uj)|  
& =    & |d\omega_{x_m}(df^m(x)\cdot \ui, df^m(x)\cdot \vj)|\\
& \leq & C \|df^m(x)\cdot \ui\| \,\|df^m(x)\cdot \vj)\|
\leq C R^2(x) \|\ui\|\, \|\vj\|\, \frac{(\lambda_i e^{2\hep})^m}{\lambda_j^m} .
\end{eqnarray*}
Since $\lambda_i e^{2\hep} < \lambda_j$, the latter converges to $0$ as $m \to \infty$, 
so $d\omega_x(\ui,\vj) = 0$.

The case $i > j$ is considered similarly by taking $m \to -\infty$.
\endofproof

\def\tyj{\tilde{y}^{(j)}}
\def\cvj{\check{v}_j}
\def\heta{\hat{\eta}}

\def\etao{\eta^{(1)}}
\def\etai{\eta^{(i)}}
\def\zetao{\zeta^{(1)}}
\def\hbeta{\hat{\beta}}

\subsection{ Proof of Lemma 4.3 (a)} 

We will consider cylinders $\cc$ of length $m \geq 1$  in $\check{R}$ with 
$\cc \cap P_0 \cap \Xi_m \neq \e$ (instead of considering cylinders $\cc$ in $R$) 
with corresponding obvious changes in the estimates we need to prove.

Let $\cc$ be a cylinder of length $m$ in  $\check{R}$. 
Fix an arbitrary $z_0 \in \cc \cap P_0 \cap \Xi_m$. Given $x_0\in \cc$, write 
$x_0 = \Phi^u_{z_0}(\xi_0) = \exp^u_{z_0}(\txi_0)$
for some $\xi_0, \txi_0 \in E^u(z_0)$ with $\txi_0 = \Psi^u_{z_0}(\xi_0)$.  
Then $ \|\xi_0\| \leq R_0 \, \diam(\cc)$. 
Set $\cc' = \tPsi \circ \Psi^{-1}(\cc) \subset \tR$, $T = \ttau_m(z_0)$ and $p = [T]$, so that $p \leq T < p+1$.

Since $m$ is the length of  $\cc'$, $\tpp^m(\cc')$ contains a whole unstable leaf
of a proper rectangle $\tR_j$.  Moreover, $z_0\in \cc \cap \Xi_m$ shows that there exists an integer
$m'$ with $m (1-\delta_0)  \leq m' \leq m $ such that  $z = \tpp^{m'}(z_0) \in P_0$. 
Let $z \in \tR_i$. By the choice of the constant $r_1 > 0$ (see the end of Sect. 4.1), there exists $y \in R_i$ such that
$B^u(y,r_1) \subset W^u_{\tR_i}(z)$ and $d(z,y) < r_0/2$. In particular, for every point $b'\in B^u(y,r_1)$
there exists $b \in \cc$ with $\tpp^{m'}(b) = b'$. Set $p' = [\ttau_{m'}(z_0)]$.
Since $\ttau$ takes values in $[0,1]$, the definition of the set $\tP_0$ shows that  $p (1-\delta_0)  \leq p' \leq p$ and
$z_{p'} = f^{p'}(z_0) \in  \phi_{[-1,1]} (P_0) = \tP_0$, so $r(z_{p'}) \geq r_0$ by (4.4). Clearly, $p' \geq \ttau_{m'}(z_0)$.
Then for every $b\in W^u_{r_1}(z_{p'})$ there exists $b \in \cc$ with $f^{p'}(b) = b'$. Consider an arbitrary
$\zeta_{p'} \in E^u(z_{p'}; r_1/R_0)$ such that $\|\zetao_{p'}\| \geq r_1/R_0$, 
and set $\zeta = \hf^{-p'}_{z_{p'}}(\zeta_{p'})$. 
Then $x = \Phi^u_z(\zeta) \in \cc$, so $\diam(\cc) \geq d(z_0,x) \geq \frac{\|\zeta\|}{R_0} \geq \frac{\|\zetao\|}{\Gamma_0 R_0}$. 
On the other hand, Lemma 3.5 in \cite{St4} (see Lemma 9.1 below)  gives
$$\|\zetao\|
 \geq  \frac{1}{\Gamma_0} \|\zetao\|'_{z_0} \geq
\frac{\|\zetao_{p'}\|}{\Gamma_0 \mu_1^{p'}}
 \geq  \frac{r_1/R_0}{\Gamma_0 \mu_1^{(1+\delta_0)p}}
\geq  \frac{r_1}{R_0\Gamma_0 \mu_1^{p} e^{\hep p}}
= \frac{r_1}{R_0\Gamma_0 \lambda_1^{p} e^{2 \hep p}} ,$$
hence $\diam(\cc) \geq  \frac{c_3}{\lambda_1^{p} e^{2 \hep p}}$,
where $c_3 =  \frac{r_1}{R^2_0\Gamma^2_0} \geq 1$. 

This proves the left-hand-side inequality in (\ref{eq:cylbound}) with $C_1 = 1/c_3$. 
The other inequality in (\ref{eq:cylbound}) follows by a similar (in fact, easier) argument. We omit the details.

\def\hb{\hat{b}}

\subsection{Proof of Lemma 4.3(b)}

Assume as in Sect. 4.2 that $L$ is a fixed constant  with $L > 3/\tau_0$.
Let $\cc$ be a cylinder  of length $m$ in $R$ such that there exists $\hz \in \cc \cap P_0 \cap \Xim_L$.
Set $\tcc = \psi(\cc)$.

Let  $\hx_0,\hz_0 \in \cc$, $\hy_0, \hb_0 \in W^s_{R_1}(\hz_0)$.
We can assume that $\cc$ is the {\bf smallest cylinder
containing} $\hx_0$ and $\hz_0$; otherwise we will replace $\cc$ by a smaller cylinder.

It is enough to consider the case when $z_0 = \hz$. Indeed, assuming 
the statement is true with $\hz_0$ replaced by $\hz$, consider arbitrary points $\hx_0, \hz_0\in \cc$.
Set $\{ y \} = W^u_R(\hy_0) \cap W^s_R(\hz)$ and $\{ b \} = W^u_R(\hb_0) \cap W^s_R(\hz)$.
Since the local unstable holonomy maps are
uniformly H\"older, there exist (global) constants $C' > 0$ and $\beta' > 0$ 
such that $d(y,b) \leq C' (d(\hy_0,\hb_0))^{\beta'}$. Thus, using the assumption, 
$$|\Delta(\hx_0,y) - \Delta(\hx_0,b)| \leq C_1 \diam(\cc) (d(y,b))^{\beta_1} 
\leq C_1 (C')^{\beta_1} \diam(\cc) (d(\hy_0,\hb_0))^{\beta'\beta_1} .$$
A similar estimate holds for $|\Delta(\hz_0,y) - \Delta(\hz_0,b)|$, so
\begin{eqnarray*}
|\Delta(\hx_0,\hy_0) - \Delta(\hx_0,\hb_0)| 
& =    &  |(\Delta(\hx_0,y) - \Delta(\hz_0,y) - (\Delta(\hx_0,b) - \Delta(\hz_0,b))|\\
& \leq & |\Delta(\hx_0,y)  - \Delta(\hx_0,b)| + |\Delta(\hz_0,y) - \Delta(\hz_0,b)|\\
& \leq & 2 C_1 (C')^{\beta_1} \diam(\cc) (d(\hy_0,\hb_0))^{\beta'\beta_1} .
\end{eqnarray*}

So, from now on we will assume that $\hz_0 = \hz \in \cc \cap P_0 \cap \Xim_L$.
Then $R(\hz_0) \leq R_0$, $r(\hz_0) \geq r_0$, etc.
Set $x_0 = \Psi(\hx_0)$, $z_0 = \Psi(\hz_0)$, $y_0 = \Psi(\hy_0) \in \tR$, $b_0 = \Psi(\hb_0)$, and then write 
$x_0 = \Phi^u_{z_0}(\xi_0) = \exp^u_{z_0}(\txi_0)$
for some $\xi_0, \txi_0 \in E^u(z_0)$ with $\txi_0 = \Psi^u_{z_0} (\xi_0)$.  
Then $ \|\xi_0\|, \|\txi_0\| \leq R_0 \diam(\tcc)$. 
Similarly, write $y_0 = \exp^s_{\tz_0}(\tv_0) = \Phi^s_{z_0}(v_0)$ 
and also $b_0 = \exp^s_{z_0}(\teta_0) = \Phi^s_{z_0}(\eta_0)$ for some 
$v_0, \tv_0, \eta_0, \teta_0 \in E^s(z_0)$ with $\tv_0 =  \Psi^s_{z_0} (v_0)$ and
$\teta_0 =  \Psi^s_{z_0} (\eta_0)$. By (3.6),
\be
\|\tv_0 - v_0\| \leq R_0 \|v_0\|^{1+\beta} \quad , \quad \|\txi_0 - \xi_0\|\leq R_0\|\xi_0\|^{1+\beta}\quad
, \quad \|\teta_0 - \eta_0\| \leq R_0\| \eta_0\|^{1+\beta} .
\ee

\subsubsection{Pushing forward}

Set $p = [\ttau_m(z_0)]$; then (\ref{eq:cylbound}) holds.
Set $q = [p/2]$. We will in fact assume that $q = p/2$; the difference with the case
when $p$ is odd is insignificant. For any integer $j \geq 0$ set
$z_j = f^j(z_0)$, $x_j = f^j(x_0)$, $y_j = f^j(y_0)$ and also
$$\hxi_j =  d\hf^j_{z_0}(0)\cdot \xi_0 \:\: , \:\: \xi_j = \hf^j_{z_0}(\xi_0) \:\: , 
\:\: \txi_j = \tf^j_{z_0}(\txi_0)\:\: , \:\:
\hv_j =  d\hf^j_{z_0}(0)\cdot v_0 \:\: , \:\: v_j = \hf^j_{z_0}(v_0) \:\: , \:\: \tv_j = \tf^j_{z_0}(\tv_0) ,$$
$$b_j = f^j (b_0) \:\:, \:\: \heta_j =  d\hf^j_{z_0}(0)\cdot \eta_0 \:\: ,\:\: 
\eta_j = \hf^j_{z_0}(\eta_0) \:\: , \:\: \teta_j = \tf^j_{z_0}(\teta_0) .$$

Since $p \geq 4n_0$, we have $q \geq 2n_0$. Notice that $\txi_0 = \Psi^u_{z_0}(\xi_0)$, 
$\tv_0 = \Psi^s_{z_0}(v_0)$, and also
$$\txi_j = \Psi^u_{z_j} (\xi_j) \:\: , \:\: \Phi^u_{z_j}(\xi_j) = x_j \quad, \quad \tv_j = \Psi^s_{z_j} (v_j) 
\quad, \quad \teta_j = \Psi^s_{z_j} (\eta_j) ,$$
so by (3.6),
\be
\|\xi_j - \txi_j \| \leq R(z_j) \|\xi_j\|^{1+\beta} \:\: , \:\: \|v_j - \tv_j\| \leq R(z_j) \|v_j\|^{1+\beta} 
\:\: , \:\: \|\eta_j - \teta_j\| \leq R(z_j) \|\eta_j\|^{1+\beta} .
\ee
Moreover, $\exp^u_{z_j}(\txi_j) = f^j(\exp^u_{z_0}(\xi_0)) = f^j(x_0) = x_j$, $\exp^s_{z_j}(\tv_j) = y_j$
and $\exp^s_{z_j}(\teta_j) = b_j$, so Lemma 4.2 implies
\be
|\Delta(x_j,y_j)  - d\omega_{z_j}(\txi_j,\tv_j)| 
\leq  C_0 \,\left[ \|\txi_j\|^2\, \|\tv_j\|^\vartheta + \|\txi_j\|^\vartheta \|\tv_j\|^2 \right]
\ee
and similarly
$$|\Delta(x_j,b_j)  - d\omega_{z_j}(\txi_j,\teta_j)| 
\leq  C_0 \,\left[ \|\txi_j\|^2\, \|\teta_j\|^\vartheta + \|\txi_j\|^\vartheta \|\teta_j\|^2 \right] $$
for every integer $j \geq 0$. From (8.4) one gets
$$|d\omega_{z_j} (\txi_j, \tv_j) - d\omega_{z_j}(\xi_j,v_j)| 
\leq 2 C_0 R(z_j)\, \|\xi_j\|\, \|v_j\| (\|\xi_j\|^\beta + \|v_j\|^\beta) ,$$
$$|d\omega_{z_j} (\txi_j, \teta_j) - d\omega_{z_j}(\xi_j,\eta_j)| 
\leq 2 C_0 R(z_j)\, \|\xi_j\|\, \|\eta_j\| (\|\xi_j\|^\beta + \|\eta_j\|^\beta) ,$$
and also\footnote{Indeed, from (8.4), 
$\|\txi_j\| \leq \|\xi_j\| (1+ R(z_j) \|\xi_j\|^\beta) \leq \|\xi_j\| (1 + R_0 e^{(p'-j)\hep} r_0/\mu_1^{p'-j} )
\leq \|\xi_j\| (1+  \frac{R_0 r_0}{(e^{-\hep} \mu_1)^{n_0}} ) \leq 2 \|\xi_j\|$, 
assuming $n_0 \geq 1$ is sufficiently large. Similarly,
$\|\tv_j\| \leq 2 \| v_j\|$ and $\|\teta_j\| \leq 2 \| \eta_j\|$.}
$\|\txi_j\|\leq 2\|\xi_j\|$, $\|\tv_j\| \leq 2 \| v_j\|$ and $\|\teta_j\| \leq 2 \| \eta_j\|$.

Using these, it follows from (8.5) that
\begin{eqnarray}
|\Delta(x_j,y_j)  - d\omega_{z_j}(\xi_j,v_j)| 
& \leq &  2 C_0 R(z_j)\, \|\xi_j\|\, \|v_j\| (\|\xi_j\|^\beta + \|v_j\|^\beta) \nonumber\\ 
&       & + 8 C_0\,\left[ \|\xi_j\|^2\, \|v_j\|^\vartheta + \|\xi_j\|^\vartheta \|v_j\|^2 \right] . 
\end{eqnarray}
and similarly
\begin{eqnarray}
|\Delta(x_j,b_j)  - d\omega_{z_j}(\xi_j,\eta_j)| 
& \leq &  2 C_0 R(z_j)\, \|\xi_j\|\, \|\eta_j\| (\|\xi_j\|^\beta + \|\eta_j\|^\beta) \nonumber\\ 
&       & + 8 C_0\,\left[ \|\xi_j\|^2\, \|\eta_j\|^\vartheta + \|\xi_j\|^\vartheta \|\eta_j\|^2 \right] . 
\end{eqnarray}
for every integer $j \geq 0$.


We will be estimating $|\Delta(x_0,y_0) - d\omega_{z_0}(\xi_0,v_0)| $.
Since $\Delta$ is $f$-invariant and $d\omega$ is $df$-invariant  we have
$\Delta(x_0,y_0) = \Delta (x_j,y_j)  \quad , \quad  d\omega_{z_0}(\xi_0,v_0) = d\omega_{z_j}(\hxi_j , \hv_j) ,$
and also $\Delta(x_0,b_0) = \Delta_{z_j}(x_j,b_j)$ and   
$d\omega_{z_0}(\xi_0,\eta_0) = d\omega_{z_j}(\hxi_j , \heta_j)$
for all $j$. (Notice that $d\hf_{x}(0) = df(x)$ for all $x\in M$.)

Since $L > 3/\tau_0$, we have   $p = [\ttau_m(z_0)] > 3m/L$, so $q= p/2 > m/L$. 
Now $z_0\in \Xim$ implies $z_0 \notin\Xi_q(p_0,\hep_0,\hd_0)$, so  
there exist at least $q-\hd_0 q$ numbers $j = 1, \ldots,q$ with $f^j(z_0)\in P_0$. 
Fix an arbitrary integer $\ell$ with
\be
(1 - \hd_0) q \leq \ell \leq q \quad , \quad z_\ell = f^\ell(z_0) \in P_0  .
\ee
It then follows from Lemma 3.1, the choice of $L_0$ and $\|\xi_\ell\| \leq r(z_\ell)$ (since $\ell \leq q= p/2$; see also Sect. 8.3.2 below) that
\be
\|\hxio_\ell - \xio_\ell\| \leq L_0 \|\xi_\ell\|^{1+\beta} .
\ee
Apart from that, using Lemma 9.7(b) below, backwards for stable manifolds, with 
$a = d\hf^{-\ell}_{z_\ell}(0)\cdot v_\ell \in E^s(z_0)$,  $b = d\hf^{-\ell}_{z_\ell}(0)\cdot \eta_\ell \in E^s(z_0)$, since
$v_0 = \hf^{-\ell}_{z_\ell}(v_\ell)$ and $\eta_0 = \hf^{-\ell}_{z_\ell}(\eta_\ell)$, it follows that
$$\|(\ao - \bo) - (\vo_0 - \etao_0)\| \leq L_0 \left[ \|v_0 - \eta_0\|^{1+\beta} 
+ \|\eta_0\|^\beta \|v_0-\eta_0\|\right]
\leq 2L_0 \| v_0 - \eta_0\| .$$
Thus,
\be
\| d\hf^{-\ell}_{z_\ell}(0) \cdot (\vo_\ell - \etao_\ell) - (\vo_0 - \etao_0)\| \leq 2L_0 \|v_0 - \eta_0\| .
\ee

In what follows we denote by $\Con$ a global constant (depending on constant like 
$C_0$, $L_0$, $R_0$ however independent of
the choice of the cylinder $\cc$, the points $x_0,z_0,y_0,b_0$, etc.) which may change from line to line.

Using  (8.9), (8.10) and the above remarks,  we obtain
\begin{eqnarray}
&      & |d\omega_{z_\ell} (\xi_\ell, v_\ell - \eta_\ell)| \nonumber\\
& \leq & |d\omega_{z_\ell} (\xio_\ell, \vo_\ell - \etao_\ell)| 
+ C_0 \sum_{i=2}^k \|\xi_\ell^{(i)}\| \, (\|v_\ell^{(i)}\| + \| \eta_\ell^{(i)}\|)\nonumber\\
& \leq &  |d\omega_{z_\ell} (\hxio_\ell, \vo_\ell - \etao_\ell)| + \Con \|\xi_\ell\|^{1+\beta} 
\|\vo_\ell - \etao_\ell\| +  C_0 \sum_{i=2}^k \|\xi_\ell^{(i)}\| \, (\|v_\ell^{(i)}\| 
+ \| \eta_\ell^{(i)}\|)\nonumber\\
& \leq & |d\omega_{z_\ell} (d \hf^{\ell}_{z_0}(0) \cdot \xio_0, \vo_\ell - \etao_\ell)| 
+  \Con \|\xi_\ell\|^{1+\beta} \|\vo_\ell - \etao_\ell\| +  C_0 \sum_{i=2}^k \|\xi_\ell^{(i)}\| \, 
(\|v_\ell^{(i)}\| + \| \eta_\ell^{(i)}\|)\nonumber\\
& =     & |d\omega_{z_0} (\xio_0, d\hf^{-\ell}_{z_\ell}(0)\cdot (\vo_\ell - \etao_\ell))|
+ \Con \|\xi_\ell\|^{1+\beta} \|\vo_\ell - \etao_\ell\| 
+  C_0 \sum_{i=2}^k \|\xi_\ell^{(i)}\| \, (\|v_\ell^{(i)}\| + \| \eta_\ell^{(i)}\|)\nonumber\\
& \leq &  |d\omega_{z_0} (\xio_0, \vo_0 - \etao_0)| + 2C_0 L_0 \|\xi_0\|\, \|v_0 - \eta_0\| 
+ \Con \|\xi_\ell\|^{1+\beta} \|\vo_\ell - \etao_\ell\| \nonumber\\
&      & +  C_0 \sum_{i=2}^k \|\xi_\ell^{(i)}\| \, (\|v_\ell^{(i)}\| + \| \eta_\ell^{(i)}\|)\nonumber\\
& \leq & \Con \, \diam(\cc) \, \|v_0 - \eta_0\| 
+ \Con \|\xi_\ell\|^{1+\beta} \|\vo_\ell - \etao_\ell\| 
+  C_0 \sum_{i=2}^k \|\xi_\ell^{(i)}\| \, (\|v_\ell^{(i)}\| + \| \eta_\ell^{(i)}\|) .
\end{eqnarray}

\subsubsection{Estimates for $\|\xi_\ell\|$, $\|v_\ell\|$ and $\| \eta_\ell\|$}

We will now use the choice of $\ell$ to estimate $\|\xi_\ell\|$, $\|v_\ell\|$ and $\| \eta_\ell\|$ by means of $\|\xi_0\|$,
$\|v_0\|$ and $\|\eta_0\|$. We will first estimate $\|\xi_q\|$, $\|v_q\|$ and $\|\eta_q\|$.

Using the definition of $\xi_j$, $p = 2q$, $z_0\in P_0$  and (3.11) we get
$\|\xi_{q}\| \leq \|\xi_q\|'_{z_q} \leq \frac{\|\xi_p\|'_{z_p}}{\mu_1^{p-q}}
\leq  \frac{\Gamma(z_p) e^{q \hep} \|\xi_p\|}{\lambda_1^{q}}
\leq  \frac{\Gamma_0 e^{2q \hep} \|\xi_p\|}{\lambda_1^{q}}\;.$
Since $\Phi^u_{z_p}(\xi_p) = x_p$ and  $d(x_p,z_p) \leq \diam(\tR_i)$, we get
$\|\xi_p\| \leq R(z_p) d(x_p,z_p) \leq R_0 e^{p \hep} r_1 < R_0 e^{p \hep}$. Thus,
\be
\|\xi_{q}\| \leq  \frac{R_0 \Gamma_0 e^{4q \hep}}{\lambda_1^{q}}\;.
\ee
Using (3.11) again (on stable manifolds) and $\|v_0\| \leq 2\delta'/R_0 < 1$, we get
\be
\|v_{q}\| = \|v_q\|'_{z_q} \leq  \frac{\|v_0\|'_z}{\mu_1^q}  \leq \frac{\Gamma_0 e^{q\hep} \|v_0\|}{\lambda_1^q} 
\leq \frac{\Gamma_0 e^{q\hep}}{\lambda_1^q} .
\ee
Similarly, $\|\eta_q\| \leq  \frac{\Gamma_0 e^{q\hep}}{\lambda_1^q}$.

Next, it follows from (\ref{eq:cylbound}) that $(\lambda_1 e^{2\hep})^{2q} \geq c_3/\diam(\cc)$, so
\be
q \geq \frac{1}{2\log (\lambda_1 e^{2\hep})}\, \log \frac{c_3}{\diam(\cc)} .
\ee
This and (8.12) give
\begin{eqnarray}
\|\xi_{q}\|
& \leq & R_0 \Gamma_0 (\lambda_1 e^{-4\hep})^{-q} = R_0 \Gamma_0 e^{-q \log (\lambda_1 e^{-4\hep})} \nonumber
 \leq  R_0 \Gamma_0 \, e^{-\frac{\log (\lambda_1 e^{-4\hep})}{2\log ( \lambda_1 e^{2 \hep})}
\log \left(\frac{c_3}{\diamf(\cc)}\right) } \nonumber \\
& =    &  R_0 \Gamma_0 \,\left(\frac{ c_3}{\diam(\cc)}\right)^{- \frac{\log\lambda_1 - 4 \hep}{2\log \lambda_1 + 4\hep}}
\leq \frac{R_0 \Gamma_0}{c_3} \, \left(\diam(\cc)\right)^{\frac{\log\lambda_1 - 4 \hep}{2\log \lambda_1 + 4\hep}} ,
\end{eqnarray}
since $\frac{\log\lambda_1 - 4 \hep}{\log \lambda_1 + 2\hep} < 1$.
Similarly, (8.13) yields
\begin{eqnarray*}
\|v_{q}\|
& \leq & \Gamma_0 (\lambda_1 e^{-\hep})^{-q}
 \leq  \Gamma_0 \, e^{-\frac{\log (\lambda_1 e^{-\hep})}{2\log ( \lambda_1 e^{2 \hep})} 
 \log \left(\frac{c_3}{\diamf(\cc)}\right) } 
\leq \frac{\Gamma_0}{c_3} \, \left(\diam(\cc)\right)^{\frac{\log\lambda_1 -  \hep}{2\log \lambda_1 + 4\hep}}\;.
\end{eqnarray*}
The same estimate holds for $\|\eta_q\|$.

We need similar estimates, however with $q$ replaced by $\ell$. 
Since $q-\ell \leq \hd_0 q$ by (8.8), as in (8.15) one obtains
\begin{eqnarray*}
\|\xi_\ell\| 
& \leq &  \|\xi_\ell\|'_{z_\ell} \leq \|\xi_q\|'_{z_q} \leq \Gamma_0 \|\xi_q\|
 \leq  \frac{R_0 \Gamma^2_0}{c_3} \, \left(\diam(\cc)\right)^{\frac{\log\lambda_1  - 4 \hep}{2\log \lambda_1 + 4\hep}} .
\end{eqnarray*}
Since $\lambda_{\tk}^{\delta_0} < e^{\hep}$ by the choice of $\hd_0$ in Sect. 4.2, we have
$\lambda_{\tk}^{q-\ell} \leq \lambda_{\tk}^{\hd_0 q} < e^{\hep q}$, and therefore
\begin{eqnarray*}
\|v_\ell\| 
 \leq  \Gamma(z_\ell) e^{(q-\ell)\hep} \lambda_{\tk}^{q-\ell}\|v_{q}\|
\leq \Gamma_0 e^{3q \hep} \|v_{q}\| \leq \Gamma^2_0 (\lambda_1 e^{-4\hep})^{-q}
\leq  \frac{\Gamma^2_0}{c_3} \, \left(\diam(\cc)\right)^{\frac{\log\lambda_1 -  4\hep}{2\log \lambda_1 + 4\hep}} ,
\end{eqnarray*}
and again the same estimate holds for $\|\eta_\ell\|$.
Thus, taking the constant $C'' > 0$ so large that  $C'' \geq R_0 \Gamma_0^2/c_3$, we get
$\|v_{\ell}\|, \|\eta_\ell\|  , \| \xi_\ell\|
\leq C'' \left(\diam(\cc)\right)^{\frac{\log\lambda_1 -  4\hep}{2\log \lambda_1 + 4\hep}} .$
Using these we get the following estimates for the terms in (8.11):
\begin{eqnarray*}
\|\xi_\ell\|\, \|v_\ell\| (\|\xi_\ell\|^\beta + \|v_\ell\|^\beta)
 \leq  2  (C'')^3 \left(\diam(\cc)\right)^{(2+\beta)\frac{\log\lambda_1 
 -  4\hep}{2\log \lambda_1 + 4\hep}}  \leq  2 (C'')^3 (\diam(\cc))^{1+\hbeta} ,
\end{eqnarray*}
where we choose 
\be
0 < \hbeta  = \min\left\{  \frac{1}{4} \, \min \{\beta , \vartheta\} \; ,\:
\frac{\log\lambda_2 - \log \lambda_1}{2 \log \lambda_1}  \right\} ,
\ee
and we use the assumption $\hep  \leq \frac{\log \lambda_1}{100} \min\{ \beta , \vartheta \}$.
Then $(2+\beta)\frac{\log\lambda_1 -  4\hep}{2\log \lambda_1 + 4\hep} \geq 1+ \hbeta$
and also\\ $(2+\vartheta)\frac{\log\lambda_1 -  4\hep}{2\log \lambda_1 + 4\hep} \geq 1+ \hbeta$
which is used in the next estimate. Similarly,
$$\|\xi_\ell\|^{1+\beta} \|v_\ell\| \leq (C'')^3 (\diam(\cc))^{1+\hbeta} ,$$
and 
\begin{eqnarray}
\|\xi_\ell\|^2\, \|v_\ell\|^\vartheta +  \|\xi_\ell\|^\vartheta \|v_\ell\|^2  
\leq 2 (C'')^3 (\diam(\cc))^{1+\hbeta} .
\end{eqnarray}

Next, for any $\xi = \xio + \xit + \ldots + \xi^{(\tk)} \in E^u(z)$ or $E^s(z)$ for some $z\in M$ set
$\cxit = \xit + \ldots + \xi^{(\tk)}$, so that $\xi = \xio + \cxit$.
Using Lemma 3.5 in \cite{St4} (see Lemma 9.1 below), $p-\ell = 2q - \ell \geq q$
and the fact that $\|\xi_{\ell}\| \leq \|\xi_p\| \leq R_0 r_1 \leq R_0$, we get
$\| \cxit_{\ell}\|'_{z_\ell}  \leq \frac{\Gamma_0 \|\cxit_{\ell}\|}{\mu_2^{q}} 
\leq  \frac{\Gamma_0 \|\xi_{\ell}\|}{\mu_2^{q}} \leq  \frac{\Gamma_0 R_0 }{\mu_2^{q}} . $
Similarly, using Lemma 3.5 in \cite{St4} (backwards for the map $f^{-1}$ on stable manifolds),
$z_0 \in P_1\subset P$, $v_0 = v_{j,1}(z_0) \in E^s(z_0,r'_0)$ and the fact that 
$\|v_0\| \leq \delta' < 1$, we get
$\| \cvt_{\ell}\|'_{z_\ell}  \leq \frac{\Gamma_0 \|v_0\|}{\mu_2^{q (1-\hd_0)}} 
\leq  \frac{\Gamma_0 }{\mu_2^{q(1-\hd_0)}} . $
Hence for $i \geq 2$ we have
$\|\xii_\ell\| \leq |\cxit_\ell| \leq \|\cxit_\ell\| \leq \frac{\Gamma_0 R_0 }{\mu_2^{q}} , $
and similarly $\di \|\vi_\ell\| \leq \frac{\Gamma_0}{\mu_2^{q(1-\hd_0)}}$. 
Using these estimates, (8.14), $\mu_2 = \lambda_2 e^{-\hep}$, and the assumtions about $\hep$,
we get
\begin{eqnarray*}
\|\xii_\ell\|\, \| \vi_\ell\|
& \leq & \Gamma_0^2 R_0\,  (\lambda_2 e^{-2\hep})^{-  2q} 
= \Gamma_0^2 R_0\, e^{-2q \log(\lambda_2e^{-2\hep})}
 \leq  \Gamma_0^2 R_0\, e^{\frac{- \log(\lambda_2e^{-2\hep})}{\log(\lambda_1 e^{2\hep})}\, 
\log\frac{c_3}{\diamf(\cc)} }\\
& \leq & \Gamma_0^2 R_0\,\left( \frac{\diam(\cc)}{c_3}\right)^{\frac{\log \lambda_2 - 2\hep}{ \log\lambda_1 + 2\hep} } 
\leq C'' (\diam(\cc))^{1+\hbeta} ,
\end{eqnarray*}
using $\hbeta \leq \frac{\log\lambda_2 - \log \lambda_1}{2 \log \lambda_1}$ by (8.16) and assuming
$C'' \geq \Gamma_0^2 R_0/ (c_3)^{\log\lambda_2/\log\lambda_1}$. Then
\begin{eqnarray*}
\frac{\log \lambda_2 - 2\hep}{ \log\lambda_1 + 2\hep}  - 1
 =       \frac{\log \lambda_2  -  2\hep  - \log \lambda_1 - 2\hep}{ \log\lambda_1 + 2\hep} 
 \geq  \frac{\log \lambda_2 - \log \lambda_1 - 4\hep}{\log \lambda_1 + 2\hep} \geq \hbeta . 
\end{eqnarray*}

\subsubsection{Final estimate}

Using  (8.11) and the above estimates for $\|\xi_\ell\|$, $\|v_\ell\|$,
$\|\xii_\ell\|\, \| \vi_\ell\|$,  we obtain
\begin{eqnarray*}
|d\omega_{z_\ell} (\xi_\ell, v_\ell - \eta_\ell)|
& \leq & \Con \, \diam(\cc) \, \|v_0 - \eta_0\| + \Con \, (\diam(\cc))^{1+ \hbeta} .
\end{eqnarray*}

Next, using (8.6) and (8.7) with $j = \ell$  and the previous estimate we get
\begin{eqnarray}
|\Delta(x_0, y_0) - \Delta (x_0, b_0)|
& =     & |\Delta(x_\ell, y_\ell) - \Delta (x_\ell, b_\ell)|\nonumber\\
& \leq & |d\omega_{z_\ell} (\xi_\ell, v_\ell - \eta_\ell)| + \Con \, (\diam(\cc))^{1+ \hbeta}\nonumber\\
& \leq & \Con \, \diam(\cc) \, \|v_0 - \eta_0\| + \Con \, (\diam(\cc))^{1+ \hbeta} .
\end{eqnarray}

Next, we consider two cases.

\ms

\noindent
{\bf Case 1.} $\diam(\cc) \leq \| v_0 -\eta_0\|^{\vartheta/2}$. Then (8.18) immediately implies
$$|\Delta(x_0, y_0) - \Delta (x_0, b_0)| \leq \Con \, \diam(\cc) \, \|v_0 - \eta_0\|^{\hbeta \vartheta/2} .$$


\noindent
{\bf Case 2.} $\diam(\cc) \geq \| v_0 -\eta_0\|^{\vartheta/2}$. 
Set $\{X'\} = W^u_R(y_{0}) \cap W^s_{R}(x_{0})$ and $X = \phi_{\Delta(x_{0},y_{0})}(X')$. 
Then $X \in W^u_{\ep_0}(y_{0})$ and it is easy to see that
$|\Delta(x_0, y_0) - \Delta (x_0, b_0)|  = |\Delta(X,b_{0})| .$
We have $X = \exp^u_{y_{0}}(\tt)$ and $b_{0} = \exp^s_{y_{0}}(\ts)$ for some $\tt \in E^u(y_{0})$ and 
$\ts \in E^s(y_{0})$. Clearly $\|\tt\| \leq \Con$.
Using Liverani's Lemma (Lemma 4.1) we get
$$|\Delta(X,b_{0})| \leq C_0 [ |d\omega_{y_{0}}(\tt, \ts)| 
+ \|\tt\|^2 \|\ts\|^{\vartheta} + \|\tt\|^{\vartheta} \|\ts\|^2]
\leq \Con \, \|\ts\|^{\vartheta} .$$
However, $\|\ts\| \leq \Con d(y_{0},b_{0}) \leq \Con \|v_{0}-\eta_{0}\|$, so
$$|\Delta(X,b_{0})| \leq \Con \|v_0 - \eta_0\|^{\vartheta} 
\leq \Con \diam(\cc) \|v_0 - \eta_0\|^{\vartheta/2} .$$
This proves the lemma.
\endofproof



\subsection{ Proof of Lemma 4.4} 

\subsubsection{Set-up -- choice of some constants and initial points}

Choosing a constant $\ep' \in (0,r_0/2)$ sufficientky small, for any $z \in M$ and any 
$z'\in B^u(z,\ep')$ the local unstable holonomy map $\hh_{z}^{z'} : W^s_{\ep'}(z) \longrightarrow W^s_{\ep_0}(z')$ 
is well defined and uniformly H\"older continuous. Replacing $\ep'$ by a smaller constant if necessary, by (3.7) 
for $z\in P_0$ and $z'\in P_0\cap B^u(z,\ep')$ the {\it pseudo-holonomy map} 
$$\hhh_{z}^{z'} = (\Phi^s_{z'})^{-1} \circ \hh_{z}^{z'} \circ \Phi^s_{z} : 
E^s (z; \ep') \longrightarrow E^s (z'; r_0) $$
is uniformly H\"older continuous, as well. Thus, there exist constants
$C' > 0$ and $\beta'' > 0$ (depending on the set $P_0$) so that for $z$, $z'$ as above we have
\be
\|\hhh_z^{z'}(u) - \hhh_z^{z'}(v)\| \leq C' \|u-v\|^{\beta''} \quad , \quad u,v\in E^s(z;\ep') .
\ee
We will assume $\beta'' \leq \beta$, where $\beta \in (0,1]$ is the constant from Sect. 3.


Fix arbitrary constants $\delta' > 0$ with
\be
(\delta')^{\beta''} < \frac{\beta_0 \kappa \theta_0}{128 L_0C_3 R_0\Gamma_0^2}  ,
\ee
$s_0$ with $0 < s_0 < \delta'/(2R^2_0)$ and $\delta''$ with 
\be
0 < \delta''  < \min\left\{\frac{\delta'}{3R_0} , \frac{\beta_0 \delta_0 \kappa}{100 R^3_0 L_0 C^2_3}  ,
 \frac{s_0 \theta_0 c_0}{4C_3 \gamma_1 R^2_0}  \right\}   ,
\ee  
Then set 
$\di \delta_0 = \frac{s_0 \theta_0}{16R_0} > 0 .$

Next, assuming $\beta'' > 0$ is taken sufficiently small and $C' > 0$ sufficiently large, for any 
$j = 1, \ldots, \ell_0$ there exists a Lipschitz\footnote{Uniform continuity is enough.} family of 
unit vectors\footnote{E.g. define
$\eta_j(Z,z) = \frac{\left((\Phi^u_z)^{-1}\circ \Phi^u_Z(r_0\eta_j(Z)/2)\right)^{(1)}}
{\|\left((\Phi^u_z)^{-1}\circ \Phi^u_Z(r_0\eta_j(Z)/2)\right)^{(1)}\|}$.}
 $\eta_j(Z,z) \in E^u_1(z)$ ($Z\in \tP_0$, $z \in B^u(Z,r_0/2)\cap \tP_0$)
such that $\eta_j(Z,Z) = \eta_j(Z)$ and for any $v \in E^s(Z)$ we have
$$| \omega_z(\eta_j(Z,z), \hhh_Z^z(v)) - \omega_Z(\eta_j(Z), v)| \leq C' d(Z,z) \|v\|^{\beta''} .$$
Fix a constant $\ep'' \in (0,\ep'/2)$ so small that $C' (\ep'')^{\beta''} < \delta'$. Then
\be
| \omega_z(\eta_j(Z,z), \hhh_Z^z(v)) - \omega_Z(\eta_j(Z), v)| \leq  \delta' \|v\|^{\beta''} 
\quad, \quad Z \in \tP_0\: , \: z \in B^u(Z,r_0)\cap \tP_0 .
\ee

Using the symbolic coding provided by the Markov family $\{R_i\}$ it is easy to see that
there exists an integer $N_0\geq 1$ such that for any integer $N \geq N_0$  we have
$\pp^N(B^u_{\ep'}(z)) \cap B^s(z', \delta'') \neq \e$ for any $z,z'\in R$ (see the notation in the beginning of Sect. 4).

Fix for a moment $Z\in P_0 $. 
Given $j = 1, \ldots, \ell_0$, since $\eta_j(Z) \in E^u_1(Z)$, by 
Lemma 8.1 and the choice of $\theta_0 > 0$ (see Sect. 4.2), there exists  $\cvj(Z)\in E^s_1(Z)$ with 
$d\omega_{Z}(\eta_j(Z) ,\cvj(Z)) \geq \theta_0$ 
and  $\|\cvj(Z)\| = 1 .$ {\bf Fix a vector $\cvj(Z)$} with the above property
for every $j$.

Set
\be
v_j (Z) = \frac{s_0}{R_0}  \cvj(Z) \in E^s_1(Z)\quad , \quad 
y_j(Z) = \Phi^s_{Z}(v_j (Z)) \in W^s_{s_0}(Z)  .
\ee
Then $s_0/R_0^2 \leq d(Z,y_j(Z)) \leq s_0$.
Since  $d\omega_{Z}(\eta_j(Z) , v_j(Z)) \geq s_0\theta_0/R_0$,  by (4.2),
\be
|d\omega_{Z}(\eta_j(Z) , v ) | \geq \frac{s_0 \theta_0}{2R_0} \quad, \quad v\in E^s(Z)\:, 
\: \|v - v_j(Z)\| \leq \frac{s_0 \theta_0}{2C_0 R_0} .
\ee

{\bf Fix an arbitrary  $N \geq N_0$}. It follows from the above that  
for each  $Z\in P_0$, each $i = 1,2$ and each $j = 1,\ldots, \ell_0$ there exists 
\be
\yjo(Z) \in \pp^N(B^u(Z,\ep')) \cap B^s(y_j (Z), \delta'')  \: \mbox{\rm and } \:
\yjt(Z) \in \pp^N(B^u(Z,\ep')) \cap B^s(Z, \delta'') .
\ee
Fix points $\yji(Z)$ with these properties; then $\yji(Z) \in W^s_{\ep_0}(Z)$. 
We have
$$\yji(Z) = \Phi^s_{Z}(\wji(Z))\quad \mbox{\rm  for some } \quad   \wji(Z) \in E^s(Z)$$
such that $\wji(Z) \in (\Phi^s_{Z})^{-1}(B^s(\yji (Z), \delta''))$.  
For $z\in B^u(Z,\ep')$ set 
\be
\wji(Z,z) = \hhh_{Z}^{z} (\wji(Z)) \in E^s(z) .
\ee
Notice that  
\be
\Phi^s_{z}(\wji(Z,z)) = \pi_{\yji(Z)}(z) . 
\ee

Given $Z\in P_0$ and $z\in B^u(Z,\ep') \cap \tP_0$,  $d(y_j(Z), \yjo(Z)) \leq \delta''$ implies
$\|\wjo(Z) - v_j (Z)\| \leq \delta'' R_0  ,$
In particular, 
$\di \frac{s_0}{2R_0} \leq \|\wjo(Z)\| \leq \frac{2s_0}{R_0} .$
Apart from that, $\|\wjt(Z)\| \leq \frac{\delta''}{R_0}$. Now (8.22) gives
\be
\|\wjt(Z,z)\| = \|\hhh_Z^z(\wjt(Z)) - \hhh_Z^z(0)\| \leq C' \|\wjt(Z)\|^{\beta''} 
\leq C' \left(\frac{\delta''}{R_0}\right)^{\beta''} < \frac{s_0}{4 R_0^3} .
\ee
A similar estimate holds for $\wjo(Z,z)$, so we get
\be
\|\wjt(Z,z)\| \leq  \frac{s_0}{2R_0^3} \leq \|\wjo(Z,z)\| \leq 2 s_0  R_0
\quad , \quad  Z\in P_0 \;,\;  z\in B^u(Z,\ep') \cap \tP_0 .
\ee

Next, (8.24) implies
$|d\omega_{Z}(\eta_j(Z) , \wjo(Z) ) | \geq \frac{s_0 \theta_0}{2R_0} ,$
while (8.22) yields 
$|d\omega_{Z}(\eta_j(Z,z) , \wjo(Z,z) ) | \geq 
\frac{s_0 \theta_0}{2R_0} - \delta' s_0 \left(\frac{4s_0}{R_0} \right)^{\beta''}\geq \frac{s_0 \theta_0}{4R_0} $
and therefore
\be
|d\omega_{Z}(\eta_j(Z,z) , \wjo(Z,z) ) | \geq 4 \delta_0 \quad, \quad Z\in P_0 \:,\: z\in B^u(Z, \ep'')\cap \tP_0 .
\ee


To finish with this preparatory section, let $\cc$ be a cylinder of length $m$ in $R$ such that 
$\cc \cap P_0 \cap \Xim_L \neq \e$,
let $Z \in \cc \cap P_0 \cap \Xim_L$, $Z_0 = \Psi(Z)$, and let $x_0  \in \Psi(\cc)$,  $z_0  \in \Psi(\cc)$   
have the form $x_0 = \Phi^u_{Z_0}(u_0)$,  $z_0 = \Phi^u_{Z_0}(w_0)$, where
\be
d(x_0,z_0) \geq \kappa \, \diam (\Psi(\cc)) 
\ee
for some $\kappa \in (0,1]$, and
\be
\left\langle \frac{w_0 - u_0}{\|w_0 - u_0\|} , \eta_j(Z_0) \right\rangle \geq \beta_0 
\ee
for some $j = 1, \ldots, \ell_0$. Fix $\kappa$ and $j$ with these properties. Set $\tcc = \Psi(\cc)$.
Then $Z_0, z_0 \in \tcc \cap \tP_0$. By the assumption on $m$, $\diam(\tcc) < \ep''$, so
$z_0\in B^u(Z_0, \ep'')$. Let $z_0 = \phi_{t_0}(z)$ for some $z\in \cc$ and $t_0 \in (-\chi, \chi)$. Set
$$x_0 = \Phi^u_{z_0}(\xi_0) \quad , \quad v_0 = d\phi_{t_0}(z) \cdot \wjo(Z, z_0) \in E^u(z_0; r_0/R_0) ,$$ 
for some $\xi_0 \in E^u(z_0; r_0/R_0)$; then $\|\xi_0\| \leq R_0\, \diam(\tcc)$.

\def\qq{\mathcal Q}

\subsubsection{Estimates for $|d\omega_{z_0} (\xio_0 , \vo_0)| $} 



Since $Z_0, z_0  \in \tP_0$ and $\|w_0\| \leq \ep'' < < r_0/R_0$, the map
$$\qq = (\Phi^u_{z_0})^{-1}\circ \Phi^u_{Z_0} : E^u(Z_0;r_0/ R_0^2) \longrightarrow E^u(z_0)$$
is well-defined and $C^{1+\beta}$. Using $d(\Phi^u_{z_0})^{-1}(z_0) = \id$,  $\qq (w_0) = 0$ and $\qq (u_0) = \xi_0$, we get
$d\qq (w_0) =   d(\Phi^u_{z_0})^{-1}(z_0) \circ d\Phi^u_{Z_0}(w_0) = d\Phi^u_{Z_0}(w_0) .$
Now (3.8) implies\footnote{{\it Proof of } (8.33): Using $C^2$ coordinates in $W^u_{r_0}(Z_0)$, we can identify $W^u_{r_0}(Z_0)$ with an
open subset  $V$ of $\R^{n^u}$ and regard $\Phi^u_{Z_0}$ and $\Phi^u_{z_0}$ as $C^{1+\beta}$ maps on $V$
whose derivatives and their inverses are bounded by $R_0$. By Taylor's formula (3.8),
$\Phi^u_{Z_0}(u_0) - \Phi^u_{Z_0}(w_0) = d\Phi^u_{Z_0}(w_0)\cdot (u_0-w_0) + \eta ,$
for some $\eta \in \R^{n^u}$ with $\|\eta\| \leq R_0 \|u_0 - w_0\|^{1+\beta}$. Hence
$d (\Phi^u_{z_0})^{-1}(z_0) \cdot (\Phi^u_{Z_0}(u_0) - \Phi^u_{Z_0}(w_0)) 
= d\Phi^u_{Z_0}(w_0)\cdot (u_0-w_0) + \eta .$
Since $Z_0 \in P_0$, by (3.9),
$$\|d\Phi^u_{Z_0}(w_0) - \id\| = \|d\Phi^u_{Z_0}(w_0) -  d\Phi^u_{Z_0}(0)\| 
\leq R_0 \|w_0\|^{1+\beta} ,$$
so $\| d\Phi^u_{Z_0}(w_0)\| \leq 2R_0$.
Using Taylor's formula again,
$$\qq (u_0) - \qq (w_0) = (\Phi^u_{z_0})^{-1} (\Phi^u_{Z_0}(u_0)) - (\Phi^u_{z_0})^{-1} (\Phi^u_{Z_0}(w_0)) =
d (\Phi^u_{z_0})^{-1}(z_0) \cdot (\Phi^u_{Z_0}(u_0) - \Phi^u_{Z_0}(w_0)) + \zeta $$
for some $\zeta$ with 
$\|\zeta\| \leq R_0 \|\Phi^u_{Z_0}(u_0) - \Phi^u_{Z_0}(w_0)\|^{1+\beta}
\leq R_0\left( 2R_0 \|w_0 - u_0\| + R_0 \|w_0 - u_0\|^{1+\beta} \right)^{1+\beta} 
\leq 9 R_0^{3} \|u_0-w_0\|^{1+\beta}.$
Thus, 
$\xi_0 = \qq (u_0) - \qq (w_0) = d\Phi^u_{Z_0}(w_0)\cdot (u_0-w_0) + \eta + \zeta ,$
where $\|\eta + \zeta\| \leq (R_0+ 9R_0^3) \|u_0-w_0\|^{1+\beta} \leq 10 R_0^3 \|u_0-w_0\|^{1+\beta}$.}
\be
\|\xi_0 - d\Phi^u_{\hz_0}(w_0) \cdot (u_0 - w_0) \| \leq 10 R_0^3 \|u_0 - w_0\|^{1+\beta}  .
\ee

\ms

Next, by (4.10) the direction of $w_0 - u_0$ is close to $\eta_j(Z_0)$. More precisely, let
$w_0 - u_0 = t \eta_j(Z_0) + u$ for some $t \in \R$ and $u \perp \eta_j(Z_0)$. 
Then for $s = t/\|w_0 - u_0\|$ we have
$$\frac{w_0 - u_0}{\|w_0 - u_0\|} = s \eta_j(Z_0) + \frac{u}{\|w_0 - u_0\|} \, ,$$
so $\di s= \left\la \frac{w_0 - u_0}{\|w_0 - u_0\|} , \eta_j(Z_0) \right \ra  \geq \beta_0$, and therefore
$t = s \|w_0 - u_0\| \geq \beta_0 \|w_0 - u_0\|$. Moreover,
\begin{eqnarray*}
\|u\|^2 
& =    & \|w_0 - w_0 - t \eta_j(Z_0)\|^2 = \|w_0 - u_0\|^2 - 2t \la w_0 - u_0 , \eta_j(Z_)) \ra +t^2\\
& =    & \|w_0 - u_0\|^2 \left(1- 2s \left\la \frac{w_0 - u_0}{\|w_0 - u_0\|}, \eta_j(Z_0) \right\ra 
+ s^2\right)  = \|w_0 - u_0\|^2 (1-2s^2 +s^2)\\
& =    & \|w_0 - u_0\|^2 (1-s^2) \leq (1 - \beta_0^2) \|w_0 - u_0\|^2\;,
\end{eqnarray*}
and therefore $\|u\| \leq \sqrt{1-\beta_0^2}\,  \|w_0 - u_0\|$. 

Since $v_0 = d\phi_{t_0}\cdot \wjo(Z, z_0) = \hhh_{Z}^{z_0} (\wjo(Z))$, it follows from (8.24) 
with $z = z_0$ and $w = v_0$ that
$|d\omega_{z_0}(\eta_j(Z,z_0), v_0)| \geq 4 \delta_1$, while (8.29) gives 
$s_0/(2R_0^3) \leq \|v_0\| \leq 2 s_0 R_0  \leq 2\delta'/R_0$.
Using $d\Phi^u_{Z_0}(0) = \id$ and (3.9), we have
$\|d\Phi^u_{Z_0}(w_0) - \id\| \leq R_0 \|w_0\|^\beta \leq R_0 (R_0\ep'')^\beta \leq R_0^2 (\ep'')^\beta .$
Moreover, $\beta_0^2  (1+ \theta_0^2/(64C_0)^2) = 1$, so
$ \beta_0^2 \theta_0^2 = (64C_0)^2 (1-\beta_0^2)$, and therefore 
$4 C_0 \sqrt{1-\beta_0^2} = \beta_0 \theta_0/16 .$
The above, (8.29), (8.23), (8.21),
(8.22), $\|\vo_0\| \leq |v_0| \leq \|v_0\|$, Lemma 8.1 and the fact 
that $\eta_j(z_0) \in E^u_1(z_0)$ imply 
\begin{eqnarray*}
&        & |d\omega_{z_0} (\xio_0 , \vo_0)| \nonumber\\
& =      &  |d\omega_{z_0} (\xi_0 , \vo_0)| 
\geq |d\omega_{z_0} (d \Phi^u_{Z}(w_0)\cdot (u_0 - w_0) , \vo_0 )|
           - |d\omega_{z_0} (\xi_0 - d \Phi^u_{Z}(w_0)\cdot (u_0 - w_0) , \vo_0 )|\nonumber\\
& \geq & t  |d\omega_{z_0} (d \Phi^u_{Z}(w_0)\cdot \eta_j(z_0) , \vo_0 )| 
-  |d\omega_{z_0} (d \Phi^u_{Z}(w_0)\cdot u , \vo_0 )|
          - 10 C_0 R_0^3 \|u_0 - w_0\|^{1+\beta}\|\vo_0\|\nonumber\\
& \geq &  \beta_0 \|u_0 - w_0\| \, [\;  |d\omega_{z_0} (\eta_j(z_0) , \vo_0 )| 
- |d\omega_{z_0} (d \Phi^u_{Z}(w_0)\cdot \eta_j(z_0) - \eta_j(z_0), \vo_0 )| \; ]\nonumber\\
&         & - C_0 (1+  R_0^2 (\ep'')^\beta)\sqrt{1-\beta_0^2}\,  \|u_0 - w_0\| \|\vo_0\|  
- 10 C_0 R_0^3 \|u_0 - w_0\|^{1+\beta}\|\vo_0\|\nonumber\\
& \geq & \|u_0 - w_0\| \, [ \; \beta_0  |d\omega_{z_0} (\eta_j(z_0) , v_0 )|  
-  \beta_0 C_0 R_0^2 (\ep'')^\beta \|v_0\| 
           - 2 C_0 \sqrt{1-\beta_0^2}\,  \|v_0\|  - 10C_0 R_0^3 (2\ep'')^{\beta}\|v_0\| \;]\nonumber\\
& \geq & \|u_0-w_0\| \, [ \; 4 \beta_0 \delta_0 -  2 \beta_0 C_0 R_0^2  \delta''  s_0
          - 4 C_0 \sqrt{1-\beta_0^2}\,  s_0  - 20 C_0 R_0^3 \delta''  s_0 \;]\nonumber\\
& \geq & \|u_0-w_0\| \, [ \; 4 \beta_0 \delta_0 -  \beta_0 \delta_0
          - \beta_0 \delta_0 - \beta_0 \delta_0  \;] = \|u_0-w_0\| \, \beta_0\delta_0 . 
\end{eqnarray*}
Combining this with  (4.9) and (3.7) gives
\be
|d\omega_{z_0} (\xio_0 , \vo_0)|  \geq \frac{\beta_0 \delta_0 \kappa }{R_0} \, \diam(\tcc)\, .
\ee

Next, set $\txi_0 = \Psi^u_{z_0} (\xi_0) \in E^u(z_0)$.  Then 
\be
\exp^u_{z_0}(\txi_0) = \Phi^u_{z_0}(\xi_0) = x_0 ,
\ee
and
\be
\frac{\kappa}{R_0}\, \diam(\tcc) \leq \|\xi_0\| \leq R_0  \diam(\tcc) . 
\ee


Next, set $\tv_0 = \Psi_{z_0}^s (v_0)\in E^s(z_0)$ and $y_0 = \exp^s_{z_0}(\tv_0)$; 
then using $v_0 = \wjo (Z,z_0)$, (8.25) and (8.27), we get
\be 
y_0 = \exp^s_{z_0}(\tv_0) = \Phi^s_{Z}(\wjo(Z,z_0)) = \pi_{\yjo(Z)}(z_0)
\in  B^s(z_0, \delta'')  .
\ee

We will now prove that
\be
|\Delta(x_0,y_0)| \geq \frac{\beta_0 \delta_0\kappa}{2 R_0} \, \diam(\tcc) .
\ee
From this and Lemma 4.3(b), (4.11) follows easily for
$d_1 \in B^s(\yj_1(Z), \delta'')$ and $d_2 \in B^s (Z, \delta'')$,
using the choice of $\delta''$.


It follows from (3.6), $\|v_0\|\leq r_0/R_0$ and $\|\xi_0\|\leq r_0/R_0$ that
$\|\tv_0 - v_0\| \leq R_0 \|v_0\|^{1+\beta}$ and $\|\txi_0 - \xi_0\|\leq R_0\|\xi_0\|^{1+\beta} ,$
and in particular $\|\tv_0\| \leq 2\|v_0\|$ and $\|\txi_0\| \leq 2\|\xi_0\| \leq 2 R_0 \diam(\tcc)$.

As in Sect. 8.3.1, set $p = [\ttau_m(z_0)]$, $q = [p/2]$, and for $j \geq 0$ define
$z_j = f^j(z_0)$, $x_j = f^j(x_0)$, $y_j = f^j(y_0)$,
$\hxi_j =  d\hf^j_{z_0}(0)\cdot \xi_0$, etc. in the same way.
By the choice of $\ep'' > 0$  all estimates in Sect. 
8.3.1 hold without change. Choosing an arbitrary $z\in \cc \cap P_0 \cap \Xim$, as before we find
$j \geq 0$ with $\pp^j(z) \in P_0$ such that (8.8) holds for $\ell = [\ttau_j(\Psi(z))]$ and
$r(z_\ell) \geq r_0$. Fix $\ell$ with these properties; then (8.9) and (8.10) hold again.

We need an estimate from below for $|d\omega_{z_\ell}(\xi_\ell, v_\ell)|$ similar to (8.11).
Instead of using Lemma 9.7 this time it is enough to use Lemma 3.1.
Since $v_\ell = \hf^{\ell}_{z_\ell}(v_0)  \in E^s(z_\ell)$ and $z_0\in P$ implies  $L(z_0) \leq L_0$, for
$w = d\hf^{-\ell}_{z_\ell}(0)\cdot v_\ell$, using Lemma 3.1, we get
\be
\|\vo_0 - \wo \| \leq L_0(z) |v_0|^{1+\beta} \leq L_0  \|v_0\|^{1+\beta} .
\ee
As in the proof of (8.11) we will now use the estimates in Sect. 8.3.2.
It  follows from Lemma 8.1, (8.9) and (8.35) that
\begin{eqnarray*}
|d\omega_{z_\ell}(\xi_\ell, v_\ell)| 
& \geq & |d\omega_{z_\ell}(\xio_\ell, \vo_\ell)| - \sum_{i=2}^k |d\omega_{z_\ell}(\xii_\ell, \vi_\ell)| \\
& \geq & |d\omega_{z_\ell}(\hxio_\ell, \vo_\ell)| - C_0 L_0 \|\xi_\ell\|^{1+\beta} \|v_\ell\| 
- C_0 \sum_{i=2}^k \|\xii_\ell\|\,  \|\vi_\ell\| \\
& =      &  |d\omega_{z_0}(d\hf^{-\ell}_{z_\ell}(0)\cdot \hxio_\ell, d\hf^{-\ell}_{z_\ell}(0)\cdot\vo_\ell)| 
               - C_0 L_0  \|\xi_\ell\|^{1+\beta} \|v_\ell\| - C_0 \sum_{i=2}^k \|\xii_\ell\|\,  \|\vi_\ell\| \\
& =      &  |d\omega_{z_0}(\xio_0, \wo)| - C_0 L_0  \|\xi_\ell\|^{1+\beta} \|v_\ell\|  
- C_0 \sum_{i=2}^k \|\xii_\ell\|\,  \|\vi_\ell\| \\
& \geq &  |d\omega_{z_0}(\xio_0, \vo_0)| - C_0 L_0 R_0  \, \diam(\tcc) \|v_0\|^{1+\beta}
        - \Con \, (\diam(\tcc))^{1+\hbeta}  .
\end{eqnarray*}
Combining this with (8.6) and (8.30) gives
\begin{eqnarray*}
|\Delta(x_0,y_0)|
& =      &  |\Delta(x_\ell, y_\ell)| \geq |d\omega_{z_\ell}(\xi_\ell,v_\ell)| 
-  8 C_0 R_0\, \|\xi_\ell\|\, \|v_\ell\| (\|\xi_\ell\|^\beta + \|v_\ell\|^\beta)\nonumber\\ 
&         & - 8 C_0\,\left[ \|\xi_\ell\|^2\, \|v_\ell\|^\vartheta 
+ \|\xi_\ell\|^\vartheta \|v_\ell\|^2 \right]\nonumber\\
& \geq & |d\omega_{z_0}(\xio_0, \vo_0)| - C_0 L_0 R_0 \,  \diam(\tcc) \|v_0\|^{1+\beta}
- \Con\, (\diam(\tcc))^{1+\hbeta} \nonumber\\
& \geq & \frac{\beta_0 \delta_0\kappa}{R_0} \, \diam(\tcc) - C_0 L_0 R_0  \, \diam(\tcc) \|v_0\|^{1+\beta}
- C''' (\diam(\tcc))^{1+\hbeta}          
\end{eqnarray*}
for some constant $C''' > 0$. Now assume
$(2\ep'')^{\hbeta} \leq \frac{\beta_0 \delta_0 \kappa }{4 R_0 C''' } ,$
and recall that $\|v_0\| \leq \delta'$ and $\diam(\tcc) \leq 2\ep''$. By (8.28), $\|v_0\|\leq 2s_0$, while
(8.20) implies $\|v_0\|^\beta \leq (\delta'')^\beta  < (\delta')^\beta < \frac{\beta_0\kappa \theta_0}{128 L_0 C_0 R_0}$. 
Thus, using (8.21),
$ C_0 L_0 R_0 \kappa\,  \diam(\tcc) \|v_0\|^{1+\beta} \leq 
C_0 L_0 R_0 \diam (\tcc) \, 2s_0 \,\frac{\beta_0\kappa \theta_0}{128 L_0 C_0 R^2_0} 
\leq \diam (\tcc) \, \frac{\beta_0 \delta_0 \kappa}{4R_0} ,$
and therefore $\Delta(x_0, y_0 ) \geq  \frac{\beta_0 \delta_0 \kappa}{2R_0} \, \diam(\tcc)$.
This proves (8.38).
\endofproof




\def\Bmt{\overline{B_{\ep_0}(\mt)}}
\def\chBo{\check{B}^{u,1}}
\def\tBo{\tB^{u,1}}
\def\hBo{\hB^{u,1}}
\def\hpi{\hat{\pi}}
\def\tU{\widetilde{U}}
\def\tr{\tilde{r}}
\def\chBo{\check{B}^{u,1}}
\def\tBo{\tB^{u,1}}
\def\hBo{\hB^{u,1}}
\def\hpi{\hat{\pi}}
\def\tL{\widetilde{L}}

\def\pij{\pi^{(j)}}
\def\wo{w^{(1)}}
\def\vo{v^{(1)}}
\def\uo{u^{(1)}}
\def\wt{w^{(2)}}
\def\vt{v^{(2)}}
\def\twt{\tw^{(2)}}
\def\tvt{\tv^{(2)}}
\def\xio{\xi^{(1)}}
\def\xit{\xi^{(2)}}
\def\txio{\txi^{(1)}}
\def\txit{\txi^{(2)}}
\def\ut{u^{(2)}}
\def\tut{\tu^{(2)}}
\def\tuo{\tu^{(1)}}
\def\tut{\tu^{(2)}}
\def\Wo{W^{(1)}}
\def\Vo{V^{(1)}}
\def\Uo{U^{(1)}}
\def\tUo{\tU^{(1)}}
\def\Wt{W^{(2)}}
\def\Vt{V^{(2)}}
\def\Ut{U^{(2)}}
\def\tWt{\tW^{(2)}}
\def\Vt{\tV^{(2)}}
\def\tUt{\tU^{(2)}}

\def\tde{\tilde{\delta}}
\def\tSt{\widetilde{S}^{(2)}}
\def\So{S^{(1)}}
\def\hr{\hat{r}}

\section{Regular distortion for Anosov flows}
\setcounter{equation}{0}

In this section we prove Lemma 4.1. Here we do not need to assume that the flow $\phi_t$ is contact.

\subsection{Expansion along $E^u_1$}
\setcounter{equation}{0}

Let again $M$ be a $C^2$ complete Riemann manifold and $\phi_t$  be a  $C^2$ Anosov flow on $M$.
Set
$$\hmu_2 = \lambda_1 + \frac{2}{3} (\lambda_2-\lambda_1)  \quad , \quad 
\hnu_1 = \lambda_1 + \frac{1}{3} (\lambda_2 - \lambda_1) .$$
Then $\hmu_2 < \mu_2 e^{-\hep}$ and $\lambda_1 < \nu_1 < \hnu_1 < \hmu_2 < \mu_2 < \lambda_2$.  
For $\hep > 0$, apart from (3.1), we assume in addition that 
$$e^{\hep} \leq \frac{2\lambda_2}{\lambda_2 + \hmu_2} .$$

For a non-empty set $X\subset E^u (x)$  set $\ell(X) = \sup \{ \|u\|  : u \in X\} .$
Given $z\in \ll$ and $p \geq 1$, setting $x = f^p(z)$, define
$$\hB^u_p (z, \delta) = \{ u \in E^u (z) : \|\hf^p_z(u) \| \leq\delta\} .$$

Fix for a moment $x\in \ll$ and an integer $p \geq 1$, set $z = f^{-p}(x)$ and given $v \in E^u (z; r(z))$, set 
$$z_j = f^j(z)\quad, \quad v_j = \hf^j_z(v) \in E^u (z_j) \quad , \quad w_j = d\hf_z^j(0)\cdot v \in E^u (z_j)$$
for any $j = 0,1,\ldots, p$ (assuming that these points are well-defined).

For any $v = \vo + \vt + \ldots+ v^{(\tk)} \in E^u(x)$ with $v^{(j)} \in E^u_j$, set
$\tvt =  \vt + \ldots+ v^{(\tk)} \in \tE_2^u(x)$.

\bs

\noindent
{\bf Lemma 9.1.}  {\it Assume that the regularity function $\hr \leq r$ satisfies 
\be
\hr(x) \leq  \min \left\{  \left( \frac{1/\hmu_2 - 1/\lambda_2}{6 
\Gamma^2 (x) D (x)} \right)^{1/\beta}  \;, \;
\left( \frac{1/\lambda_1 - 1/\hnu_1}{6 e^{3\hep} \Gamma^2(x) D(x)} \right)^{1/\beta} \right\}\;
\ee
for all $x\in \ll$. Then for any $x\in \ll$ and any $V = \Vo + \tVt \in E^u(x; \hr(x))$, setting $y = f^{-1}(x)$
and $U = \hf_x^{-1}(V)$, we have
\be
\|\tUt\|'_{y} \leq \frac{\|\tVt\|'_x}{\hmu_2} ,
\ee
and
\be
\|\Uo\|'_{y} \geq \frac{\|\Vo\|'_x}{\hnu_1} \;.
\ee
Moreover, if $V, W \in E^u(x; \hr(x))$ and $\Wo = \Vo$, then for $S = \hf_x^{-1}(W)$ we have
\be
\|\tUt - \tSt\|'_{y} \leq \frac{\|\Vt - \tWt\|'_x}{\hmu_2} ,
\ee
and, if $\tWt =   \Vt \in E^u(x; \hr(x))$ and  $S = \hf_x^{-1}(W)$ again, then}
\be
\|\Uo - \So\|'_{y} \geq \frac{\|\Vo - \Wo\|'_x}{\hnu_1}  .
\ee

\ms

\noindent
{\it Proof.} The estimates (9.2) and (9.3) follow from Lemma 3.5 in \cite{St4}, while the proofs of (9.4)
and (9.5) are similar,  so we omit the details.
\endofproof

\bs

Next,  for any $y \in \ll$, $\ep\in (0,r(y)]$ and $p \geq 1$ set 
$\hBo_p(y,\ep) =   \hB^u_p(u;\ep) \cap E^u_1(y) .$

Replacing the regularity function with a smaller one, we may assume that
\be
L(x) (\hr(x))^\beta \leq \frac{1}{100 n_1} \quad , \quad x\in \ll ,
\ee
where $n_1 = \dim (E^u_1(x))$.

The proof of the following lemma is similar to the proof of Proposition 3.2 in \cite{St3}. We omit the details.

\bs

\noindent
{\bf Lemma 9.2}  {\it 
Let  $z\in \ll$ and $x = f^p(z)$ for some integer $p \geq 1$, and let $\ep \in (0, \tr(x)]$. Then
\be
\ell (\hB^{u}_p(z,\ep)) 
\leq  2 k \Gamma^3(x) \,  \ell (\hBo_p(z,\ep) ) .
\ee
Moreover for any $\ep' \in (0, \ep]$ there exists $u \in \hBo_p(z,\ep')$ with}
\be
\|u\| \geq \frac{\ep'}{2 k \ep \Gamma^2(x)}\ell (\hB^{u}_p(z,\ep))\quad \mbox{\rm and} \quad \|\hf^p_z(u)\| \geq \ep'/2 .
\ee

\bs

To prove the main result in this section, it remains to compare diameters of sets of the form $\hBo_p(y,\ep)$. 

\bs  

\noindent
{\bf Lemma 9.3.} {\it There exist a regularity function $\hr(x) < 1$ ($x\in \ll$)  such that:}

\ms

(a) {\it For any  $x\in \ll$ and any $0 < \delta \leq  \ep \leq \hr(x)$ we have
\be
\ell \left( \hBo_{p}(f^{-p}( x),\ep)  \right) 
\leq  16 n_1 \frac{\ep}{\delta} \, \ell \left( \hBo_p (f^{-p} ( x) , \delta) \right) \;
\ee
for any integer $p \geq  1$.}

\ms

(b) {\it For any  $x\in \ll$ and any $0 < \ep \leq \hr(x)$  and any $\rho \in (0,1)$, for any $\delta$ with
$0 < \delta \leq  \frac{\rho \, \ep}{16 n_1} $
we have
\be
\ell \left( \hBo_p(f^{-p}(x),\delta)  \right) \leq  \rho \,  \ell \left( \hBo_p(f^{-p} (x) , \ep) \right) \;
\ee
for any  integer $p \geq 1$.}

\bs  

\noindent
{\bf Theorem 9.4.} {\it There exist a regularity function $\hr(x) < 1$ ($x\in \ll$) such that:}

\ms

(a) {\it For any  $x\in \ll$ and any $0 < \delta \leq  \ep \leq \hr(x)$ we have
$\di \ell \left( \hB^{u}_p(z,\ep)  \right) \leq  \frac{32 \tk n_1 \Gamma^3(x)  \ep}{\delta}  
\, \ell \left(\hB^{u}_p(z,\delta) \right) \;$
for any integer $p \geq  1$, where $z = f^{-p}(x)$.}

\ms

(b) {\it For any  $x\in \ll$, any $0 < \ep \leq \hr(x)$, any  $\rho \in (0,1)$ and any $\delta$ with
$0 < \delta  \leq \frac{\rho \ep}{32 \tk n_1 \Gamma^3(x) }$
we have
$\di \ell \left( \hB^{u}_p(z ,\delta) \right) \leq  \rho \, \ell \left( \hB^{u}_p( z ,\ep)  \right) $
for all integers $p \geq 1$, where $z = f^{-p}(x)$.}

\ms

(c) {\it For any $x\in \ll$, any  $0 < \ep' < \ep \leq \hr(x)/2$, any 
$0 < \delta < \frac{\ep'}{100n_1}$ and any integer $p \geq 1$,
setting $z = f^{-p}(x)$, there exists $u \in \hBo_p(z,\ep')$ such that for every $v \in E^u(z)$ with
$\|\hf^p_z(u) - \hf^p_z(v)\| \leq \delta$ we have} 
$\di \|v\| \geq \frac{\ep'}{4 \ep \tk \Gamma^3(x)} \, \ell(\hB^u_p(z,\ep)) .$

\bs

Using Lemma 9.3, we will now prove Theorem 9.4. The proof of Lemma 9.3 is given
in the next sub-section. In fact, part (c) above is a consequence of Lemmas 3.1 and 
9.2 and does not require Lemma 9.3.

\bs

\noindent
{\it Proof of Theorem} 9.4. Choose the function $\hr(x)$ as in Lemma 9.3.

\ms

(a) Let $0 < \delta < \ep \leq \tr(x)$. 
Given an integer $p \geq 1$, set $z = f^{-p}(x)$. Then Lemmas 9.2 and 9.3 and (9.7) imply
\begin{eqnarray*}
\ell (\hB^{u}_p(z,\ep)) 
 \leq  2 \tk \Gamma^3(x) \, \ell (\hBo_p(z,\ep))
\leq  2 \tk \Gamma^3(x)  \, 16 n_1  \frac{\ep}{\delta}\, \ell (\hBo_p(z,\delta))
 \leq   32 \tk n_1 \Gamma^3(x)  \frac{\ep}{\delta} \, \ell (\hB^{u}_p(z,\delta)) .
\end{eqnarray*}

(b) Let $x\in \ll$ and $0 < \ep \leq \hr(x)$. Given $\rho \in (0,1)$, set
$\rho' = \frac{\rho}{2\tk \Gamma^3(x)} < \rho$.  By Lemma 9.3(b), if
$0 < \delta \leq \frac{\rho' \ep}{16 n_1}$ then (9.10) holds with $\rho$ replaced by $\rho'$
for any  integer $p \geq 1$ with $z = f^{-p}(x)$. Using this and Lemma 9.2 we get
\begin{eqnarray*}
\ell (\hB^{u}_p(z,\delta) ) \leq   2 \tk \Gamma^3(x) \, \ell (\hBo_p(z,\delta))
 \leq    2 \tk \Gamma^3(x) \,\rho'\,  \ell (\hB^{u}_p(z,\ep)) =  \rho\,  \ell (\hB^{u}_p(z,\ep)) ,
\end{eqnarray*}
which completes the proof.

\ms

(c) 
Given $x\in \ll$, $z = f^{-p}(x)$, let $\ep'$, $\ep$ and  $\delta$ be as in the assumptions. 
Let $u \in \hBo_p(z,\ep')$ be such that $\|u\|$ is the maximal possible.
By Lemma 9.2, for $U = \hf^p_z(u) \in E^u_1(x)$ we have $\ep'/2 \leq \|U\| \leq \ep'$. 
Setting $W = d\hf^p_z(0) \cdot u \in E^u_1(x)$, Lemma 3.1 and (9.6) give
$\| W - U\| \leq L(x) |U|^{1+\beta} \leq \frac{\|U\|}{100 n_1} ,$
so $\|W\| \leq \frac{101 \ep'}{100}$. 

Let $v = (\vo, \tvt) \in E^u(z)$ be such that for $V = \hf^p_z(v)$ we have $\|V - U\| \leq \delta$.
Then $|V-U| \leq \delta$, so $\|\Vo - \Uo\| \leq \delta$ and $\| \tVt \| \leq \delta$.  

Set $S = d\hf^p_z(0) \cdot v$; then $\So = d\hf^p_z(0) \cdot \vo$. By Lemma 3.1 and (9.7),
$\| \So - \Vo\| \leq \frac{|V|}{100 n_1} \leq \frac{\|V\|}{100 n_1} \leq \frac{\ep' + \delta}{100 n_1} ,$
so
$$\| \So - \Wo\| 
 \leq  \|\So - \Vo\| + \|\Vo - \Uo\| + \| \Uo - \Wo\|
 \leq  \frac{\ep' + \delta}{100 n_1} + \delta + \frac{\ep'}{100 n_1} < \frac{\ep'}{30 n_1} .$$

Choose an orthonormal basis $e_1, \ldots, e_{n_1}$ in $E^u_1(x)$ such that $W = \Wo = c_1 e_1$ 
for some $c_1 \in [\ep'/3,\ep']$.
Let $\So = \sum_{i=1}^{n_1} d_i e_i$. Then the above implies $|d_1-c_1| \leq \frac{\ep'}{30 n_1}$ and
$|d_i| \leq \frac{\ep'}{30 n_1}$ for all $ i = 2,\ldots,n_1$.

Notice that for any $i = 1, \ldots, n_1$, $u' = d\hf^{-p}_x(0) \cdot (\ep' e_i/2) \in \hBo_p(z,\ep')$. 
Indeed, by Lemma 3.1 and (9.6),
$\|\hf^p_z(u') - d \hf^p_z(0) \cdot u' \| \leq \frac{\|\ep' e_i/2\|}{100n_1} = \frac{\ep'}{200n_1}$, so
$\| \hf^p_z(u')\| \leq \|d \hf^p_z(0) \cdot u' \| + \frac{\ep'}{200n_1} 
= \frac{\ep'}{2} + \frac{\ep'}{200 n_1} < \ep'$. By the choice of
$u$, this implies $\|u'\| \leq \|u\|$, so $\| d\hf^{-p}_x(0) \cdot e_i\| \leq \frac{2 \|u\|}{\ep'}$ 
for all $i = 1, \ldots, n_1$.

The above yields
$$\|d_1 \, d\hf^{-p}_z(0) \cdot e_1 \| \geq \|c_1 d\hf^{-p}_z(0) \cdot e_1 \| 
- \|(d_1- c_1) d\hf^{-p}_z(0) \cdot e_1 \|
\geq \|u\| - \frac{\ep'}{30 n_1} \cdot \frac{2\|u\|}{\ep'}  = \|u\| \left(1 - \frac{1}{15n_1}\right) .$$
Moreover, for $i \geq 2$ we have
$\|d_i \, d\hf^{-p}_z(0) \cdot e_i \| \leq \frac{\ep'}{30 n_1} 
\cdot \frac{2\|u\|}{\ep'}  = \frac{\|u\|}{15n_1} .$
Hence
\begin{eqnarray*}
\|\vo\|
& =      & \| d \hf^{-p}_x(0) \cdot \So\| = \left\| \sum_{i=1}^{n_1} d_i \, d\hf^{-p}_z(0) \cdot e_i \right\|\\
& \geq & \|d_1 \, d\hf^{-p}_z(0) \cdot e_1 \|  - \sum_{i=2}^{n_1} \|d_i \, d\hf^{-p}_z(0) \cdot e_i \| 
\geq   \|u\| \left(1 - \frac{1}{15n_1}\right) - n_1 \, \frac{\|u\|}{15n_1} > \frac{\|u\|}{2} .
\end{eqnarray*}
Combining this with Lemma 9.2 gives, 
$\|v\| \geq |v| \geq \|\vo\| > \frac{\|u\|}{2} 
\geq \frac{\ep'}{4 \tk \ep \Gamma^3(x) }\, \ell( \hB^u_p(z,\ep))$.
\endofproof

\bs

What we actually need later  is the following immediate consequence of Theorem 9.4 which 
concerns sets of the form
$$B^u_T(z,\ep) = \{ y \in W^u_\ep(z) : d(\phi_T(y), \phi_T(z)) \leq \ep\}\;,$$
where $z\in \ll$, $\ep > 0$ and $T > 0$. 

\bs

\noindent
{\bf Corollary 9.5.} {\it 
There exist an $\hep$-regularity function $\hr(x) < 1$ ($x\in \ll$) and a global constant $L_1 \geq 1$ such that:}

\ms

(a) {\it We have
$\diam \left( B^{u}_T(z,\ep)\right)  \leq  L_1\,\Gamma^3(x) 
\frac{\ep}{\delta}  \, \diam \left(B^{u}_T (z,\delta) \right) $
for any $x\in \ll$, any $0 < \delta \leq \ep \leq \hr(x)$ and any  $T > 0$, where $z = \phi_{-T}(x)$.}

\ms

(b) {\it For any  $x\in \ll$, any $0 < \ep \leq \hr(x)$, any $\rho \in (0,1)$ and any $\delta$ with
$0 < \delta \leq \frac{\rho \ep}{L_1 \Gamma^3(x)}$ we have
$\diam \left( B^{u}_T (z ,\delta)\right) \leq  \rho \, \diam \left( B^{u}_T ( z ,\ep)\right)$
for all $T > 0$, where $z = \phi_{-T}(x)$.}

\ms

(c) {\it For any  $x\in \ll$, any $0 < \ep' < \ep \leq \hr(x)$, any $0 < \delta \leq \frac{\ep'}{100 n_1}$ and any
$T > 0$, for $z = \phi_{-T}(x)$ there exists $z' \in B^u_T(z,\ep')$ such that 
$d(z,y) \geq \frac{\ep'}{L_1 \ep \Gamma^3(x)} \, \diam(B^u_T(z,\ep)) .$
for every $y \in B^u_T(z', \delta)$.}

\def\hq{\hat{q}}

\subsection{Linearization along $E^u_1$}

Here we prove Lemma 9.3 using arguments similar to these in the proofs of Theorem 3.1 and Lemma 3.2
in \cite{St4}.
 
We use the notation from Sect. 9.1. Let $\hr(x)$, $x\in \ll$, be as in Lemma 9.1.

\bs 

\noindent
{\bf Proposition 9.6.} {\it There exist regularity functions $\hr_1(x) \leq \hr(x)$ 
and $L(x)$, $x\in \ll$, such that:}

\ms

(a) {\it For every $x\in \ll$ and every $u\in E^u_1(x; \hr_1(x))$ there exists
$$F_x(u) = \lim_{p\to\infty} d\hf^p_{f^{-p}(x)}(0)\cdot \hf_x^{-p}(u) \in E^u_1(x; \hr(x)) .$$
Moreover, $\| F_x(u) - u\| \leq L(x)\, \|u\|^{1+\beta}$ 
for any $u\in E^u_1(x,\hr_1(x))$ and any integer $p \geq 0$.}

\ms 

(b) {\it The maps $F_x : E^u_1( x; \hr_1(x)) \longrightarrow F_x (E^u_1( x; \hr_1(x))) 
\subset E^u_1(x; \hr(x))$ ($x\in \ll$)
are uniformly Lipschitz. More precisely, 
$$\|F_x(u) - F_x(v) - (u-v)\| \leq C_1 \, [\|u-v\|^{1+\beta} + \| v\|^\beta\cdot \|u-v\|] 
\quad, \quad x\in \ll\:\:, \: u,v\in E^u_1(x; \hr_1(x)) .$$
Assuming that $\hr_1(x)$ is chosen sufficiently small, this yields}
$$\frac{1}{2} \|u-v\| \leq \|F_x(u) -F_x(v)\| \leq 2\|u-v\| \quad, 
\quad x\in \ll\:\:, \: u,v\in E^u_1(x; \hr_1(x)) .$$

(c) {\it For any $x \in M$ and any integer $q\geq 1$, setting $x_q = f^{-q}(x)$, we have
$$d\hf_{x_q}^q(0) \circ F_{x_q} (v) = F_x \circ \hf_{x_q}^q (v)$$ 
for any $v \in E^u_1 (x_q; \hr_1(x_q))$ with $\|\hf^q_{x_q}(v)\| \leq \hr_1(x)$.}

\bs

\def\ao{a^{(1)}}
\def\bo{b^{(1)}}

As in \cite{St4} this is derived from the following lemma. Part (b) below is a bit stronger than 
what is required here, however
we need it in this form for the proof of Lemma 4.2 in Sect. 8.

\bs

\noindent
{\bf Lemma 9.7.}  {\it There exist regularity functions $\hr_1(x)$ and $L(x)$, $x\in \ll$ with the following
properties:}

\ms

(a) {\it If $x\in M$, $z = f^{p}(x)$ and $\|\hf^p_z(v)\| \leq r(x)$ for some $v\in E^u_1 (z; \hr_1(z))$ and some
integer $p \geq 1$, then $\|d\hf^p_z(0)\cdot v\| \leq 2 \|\hf^p_z(v)\|$ and
$\|d\hf_z^p(0)\cdot v - \hf^p_z(v)\| \leq L(x)\, \|\hf^p_x(v)\|^{1+\beta} .$
Similarly, if $\|d\hf^p_z(0)\cdot v\|\leq \hr_1(x)$ for some $v\in E^u_1 (z)$ and some integer $p \geq 1$, then
$\|\hf^p_x(v) \| \leq 2 \|d\hf^p_x(0)\cdot v\|$ and}
$\|\hf_x^p(v) - d\hf^p_x(0)\cdot v\| \leq L(x)\, \|d\hf^p_x(0)\cdot v\|^{1+\beta} .$

\ms

(b)  {\it For any $x \in \ll$ and any integer $p \geq 1$, setting $z = f^{-p}(x)$,  the map
$$F^p_x = d \hf^p_z(0) \circ (\hf^p_x)^{-1}:  E^u (x; \hr_1(x)) \longrightarrow E^u (x; \hr(x)) ,$$
is such that
\be
\left\| \left[ (F^p_x (a))^{(1)} - (F^p_x (b))^{(1)} \right] - [\ao - \bo ]\right\|  
\leq L(x) \, \left[ \|a - b \|^{1+\beta} + \| b \|^\beta \cdot \|a - b\| \right]  
\ee
for all $a,b \in E^u (x ; \hr_1(x))$. Moreover,}
\be
\frac{1}{2} \| a - b\| \leq \left\| d\hf^p_z(0)\cdot \left[ (\hf^p_x)^{-1}(a) - (\hf^p_x)^{-1}(b)\right] \right\|
\leq 2  \| a - b\| \quad, \quad a,b \in E^u_1 (x ;  \hr_1(x)) .
\ee

\ms

\noindent
{\it Proof of Lemma} 9.7. Set $\hr_1(x) = \hr(x)/2$, $x\in \ll$.

Part (a) follows from Lemma 3.1 (see also the Remark after the lemma).
The proofs of the other parts are almost one-to-one repetitions of arguments from the proof of
Lemma 3.2 in \cite{St4}, so we omit them.
\endofproof

\bs

\noindent
{\it Proof of Proposition 9.6.} This is done following the arguments from the proof of Theorem 3.1
in \cite{St3}. We omit the details again.
\endofproof

\bs

\def\tBo{\widetilde{B}^{u,1}}

For $z\in \ll$, $\ep \in (0, \hr_1(z)]$ and an integer $p \geq 0$ set
$$\tB^{u,1}_p(z,\ep) = F_z(\hB^{u,1}_p(z,\ep)) \subset E^u_1(z; \hr(z)) .$$
Then, using  Proposition 9.6(c) we get
\be
d\hf^{-1}_x(0) (\tBo_{p+1}(x, \delta )) \subset \tBo_p({f^{-1}(x)}, \delta) 
\quad , \quad x\in \ll\;, \; p \geq 1 .
\ee
Indeed, if $\eta \in \tBo_{p+1}(x, \delta)$, then $\eta = F_x(v)$ for some 
$v \in \hBo_{p+1}(x,\delta)$, and then clearly
$w = \hf^{-1}_x(v) \in \hBo_p(x, \delta)$. Setting $y = f^{-1}(x)$, by Proposition 9.6(c), 
$\eta = F_x(v) = F_x (\hf_y(w)) = d\hf_y(0) \cdot (F_y(w))$, so $d\hf_x^{-1}(0)\cdot \eta 
= F_y(w) \in \tBo_p(y,\delta)$.
Moreover, locally near $0$ we have an equality in (9.13), i.e. if $\delta' \in (0,\delta)$
is sufficiently small, then $d\hf^{-1}_x(0) ( \tBo_{p+1}(x,\delta)) \supset \tBo_p(f^{-1}(x),\delta')$.

To prove part (a) of Lemma 9.3 we have to establish the following lemma which is similar to Lemma 4.4 
in \cite{St4} (see also the Appendix in \cite{St4}), and the proof uses almost the same argument.

\bs

\noindent
{\bf Lemma 9.8.} {\it Let $x\in \ll$ and let  $0 < \delta \leq  \ep \leq \hr_1(x)$. Then 
$$\ell \left( \tB^{u,1}_p( f^{-p}(x),\ep)  \right) 
\leq  4 n_1  \frac{\ep}{\delta}   \, \ell \left(\tB^{u,1}_p ( f^{-p}(x) , \delta) \right) $$
for any  integer $p \geq 0$, where $n_1 = \dim (E^u_1(x))$.}
\endofproof

\bs

Lemma 9.3(b) is a consequence of the following.

\bs

\noindent
{\bf Lemma 9.9.} {\it Let $x\in \ll$ and let  $0 <  \ep \leq \hr_1(x)$ and $\rho \in (0,1)$. 
Then for any $\delta$ with
$0 < \delta \leq  \frac{\rho \, \ep}{4 n_1} $
 we have
$\di \ell \left( \tB^{u,1}_p( f^{-p}(x),\delta)  \right) 
\leq \rho  \, \ell \left(\tB^{u,1}_p ( f^{-p}(x) , \ep) \right) $
for any  integer $p \geq 0$.}

\ms

\noindent
{\it Proof of Lemma 9.9.} As in the proof of Lemma 4.1(b) in \cite{St4},
we have to repeat the argument in the proof of Lemma 9.8. We omit the details.
\endofproof

\def\htau{\hat{\tau}}
\def\tcc{\tilde{\cc}}

\subsection{Consequences for cylinders in Markov partitions}

Here we prove Lemma 4.1 using arguments similar to these in Sect. 4 in \cite{St3}. 
We sketch the argument for completeness.
We use the notation from Sect. 4.

Let $\hr(x)$ be the canonical $\ep$-regularity function from Theorem 9.4 and Corollary 9.5. Here
$\ep \in (0,\hep]$ is some constant depending on $\hep$. Then (see the end of Sect. 3.2) there exists 
a constant $\hr'_0 > 0$ such that $\hr(x) \geq \hr'_0$ for all $x\in P_0$. {\bf Fix $\ep$ and
$\hr'_0$ with these properties.}

Let $S > 0$ be a {\it Lipschitz constant} for the projection 
$\psi : \cup_{i=1}^{k_0} \phi_{[-\ep,\ep]}(D_i) \longrightarrow \cup_{i=1}^{k_0} D_i $
along the flow,
i.e. for all $i = 1,\ldots,k_0$ and all $x\in \phi_{[-\ep,\ep]}(D_i) $ we have $\psi(x) = \pr_{D_i}(x)$.  
Let $c_0$, $\gamma$ and $\gamma_1$ be the constants from (2.1). Next, assuming that the constant $\ep > 0$ 
is chosen so that $e^{\ep}/\gamma < 1$, fix an integer $d_0 \geq 1$ such that 
\be
\frac{2 k \Gamma^3_0 e^{2\ep} r_1}{\hr'_0} < (\mu_1 e^{\ep})^{d_0} 
\quad , \quad  \frac{1}{c_0 (\gamma e^{-\ep})^{d_0}} < \frac{\hr'_0}{2} .
\ee
Set
\be
r'_0 = \hr'_0 e^{- (d_0 + 1) \ep} .
\ee

\ms

\noindent
{\it Proof of Lemma} 4.1. 
First note the following. Let $z\in \tR_j$ be such that $\tpp^{d_0+1}(z) \in \tP_0$. Then
$z\in C_V[\ii']$ for some $\ii' = [i_0, \ldots,i_{d_0+1}]$ with $i_0 = j$, where
$V = W^u_{\tR}(z)$.  Set $\ii = [i_0, \ldots,i_{d_0}]$. We claim that
\be
C_V[\ii] \subset B_V(z, r'_0) \quad \mbox{\rm and } \quad r(z) \geq r'_0 .
\ee
Indeed, by (2.1) and (9.14), $\diam(C_V[\ii]) \leq \frac{1}{c_0\gamma^{d_0+1}} < r'_0/2$. 
On the other hand, $\hr(x)$ is a Lyapunov  $\hep$-regularity function and $y = \tpp^{d_0+1}(z) \in \tP_0$ 
and the definition of $\tP_0$ show that $\hr(y) \geq r_0 $. Also recall that $0 < \tau(x) \leq 1$ for all 
$x\in \tR$ by the choice of the Markov family.  Now using (9.15), we get
$\hr(z) \geq \hr(y) e^{-\tau_{d_0+1}(z) \ep} \geq \hr'_0 e^{-(d_0+1)\ep } = r_0' >   2\, \diam(C_V[\ii]) .$
This proves (9.16).

\medskip

(a)  Assume that $m > d_0$, and let $\ii = [i_0,i_1, \ldots,i_m]$ and 
$\ii' = [i_0,i_1, \ldots,i_m, i_{m+1}]$ be admissible sequences. Let  $\cc = C[\ii]$ and $\cc' = C[\ii']$
be the corresponding cylinders in $\tR$. Assume that  there exists $z \in \cc' \cap P_0$
with $\pp^{m+1}(z) \in P_0$.

Fix such a point $z \in \cc'$; then  $y = \tpp^{m+1}(z)  \in \tP_0$ and $\tpp^j(z) \in \tR_{i_j}$ 
for all $j = 0,1,\ldots,m+1$. Set $\tpp^{m-d_0}(z) = x$, $V = W^u_{\tR}(x)$, Since $\tpp^{d_0+1}(x) = y\in \tP_0$, 
we have $\hr(y) \geq \hr'_0$, so
$\hr(x) \geq r'_0$.

Consider the cylinders
$$\tcc' = C_{V}[i_{m-d_0}, i_{m-d_0+1}, \ldots, i_m, i_{m+1}] \subset 
\tcc = C_{V}[i_{m-d_0}, i_{m-d_0+1}, \ldots, i_m]  \subset V .$$
Since $\tpp^{d_0+1}(x) = y$, using (9.20) we get $\tcc \subset B_{V}(x, r'_0)$.
On the other hand it is easy to see using (2.1) that 
$\tcc' \supset B_{V}(x,c_0 \hr'_0/\gamma_1^{d_0+1})$.
Corollary 9.5(a) with $x$ and $z$ as above, $T = \tau_{m-d_0}(z) > 0$ and
$0 < \delta = \delta_3 = \frac{c_0 \hr'_0}{B \gamma_1^{d_0+1}} < \ep = r'_0$,
combined with (9.16), gives
$\diam(B^u_T(z,\delta_3)) \geq \frac{\delta_3}{B L_1 \Gamma_0^3 r'_0}\, \diam(B^u_T(z, r'_0)) .$
However, using the above information about $\tcc$ and $\tcc'$, as in the proof of 
Proposition 3.3 in \cite{St2},  one easily observes that $\cc' \supset B^u_T(z,\delta_3)$ 
and $\cc \subset B^u_T(z, B r'_0)$. Thus,
$\diam (\cc') \geq \frac{\delta_3}{B L_1 \Gamma_0^3 r'_0}\, \diam(\cc)$. 

This proves part (a) for $m > p_0$. Since there are only finitely many cylinders of length $\leq p_0$,
it follows immediately that there exists $\rho_1 \in (0, \frac{\delta_3}{B L_1 \Gamma_0^3 r'_0}\,]$ 
which satisfies the requirements of part (a).

\ms

(b) Fix an integer $q' \geq 1$ so large that
$\frac{1}{c_0 \gamma^{q'}} \leq \delta_0 $
and set $r''_0 = r_0 e^{-q' \hep}$.
Let $\rho'\in (0,1)$. 
It follows from Corollary 9.5(b) that for $z\in \tR \cap \ll$ with $\Phi_T(z) \in \tP_0$ for some $T > 0$  we have
$$\diam( B^u_T(z, B\delta)) \leq \rho' \, \diam( B^u_T(z, r''_0/B)) ,$$
provided $ 0 < \delta \leq \delta_0 = \frac{\rho' r''_0}{B^2 L_1 \Gamma^3_0}$. The rest of the proof is now very similar to the
proof of Proposition 3.3(b) in \cite{St2}, and we omit the details.

\ms

(c) Take the integer $q_0 \geq 1$ so large that 
$\frac{1}{c_0 \gamma^{p_0+q_0}} < \delta/B = \frac{\ep'}{100 B n_1}$, where $n_1 = \dim(E^u_1)$.

Let again $m > d_0$,  let $\ii = [i_0,i_1, \ldots,i_m]$ be an admissible sequence, let $\cc = \cc_W[\ii]$
be the corresponding cylinder in an unstable leaf $W$ in $\tR$. 
Let $z\in \cc \cap \tP_0$ and let $\tpp^m(z) = z'$.
Set $z''= \tpp^{m-d_0}(z)$, $V = W^u_{\tR}(z'')$.
If $z' = \phi_T(z)$ and $z'' = \phi_t(z)$; 
then $\phi_{T-t}(z'') = z'$, so $T-t = \ttau_{d_0}(z'') < d_0$.
Thus, $\hr(z'') \geq \hr(z') e^{-d_0 \ep} \geq \hr'_0 e^{-d_0\ep} > r'_0$. As in part (a), for the cylinder
$\tcc = C_{V}[i_{m-d_0}, i_{m - d_0+1}, \ldots, i_m]$ in $V$, we have 
$$z'' \in B_V(z'', c_0\hr'_0/ \gamma_1^{d_0}) \subset \tcc = \tpp^{m-d_0}(\cc) \subset B_V (z'', r'_0) .$$
Setting $\ep' = c_0\hr'_0/ \gamma_1^{d_0} < \ep = B r'_0$, it follows from Corollary 9.5(c) that 
for $0 < \delta = \frac{\ep'}{100 n_1}$ there exists $x \in B^u_t(z,\ep')$ such that for every $y \in W^u_\delta (z)$ with
$d(\phi_t(y),\phi_t(x)) \leq  \delta$ we have 
\be
d(z,y) \geq \frac{\ep'}{L_1 \ep \Gamma_0^3} \, \diam(B^u_t (z, B r'_0)) 
\geq \frac{c_0r_0}{L_1 B r'_0\Gamma_0^3} \, \diam(\cc),
\ee
since $\cc \subset B^u_t(z, B r'_0)$.

Let  $x\in \cc$ and let
$\cc' = C[\ii'] = C[i_0,i_{1}, \ldots, i_{m+1}, \ldots, i_{m+q_0}]$
be the sub-cylinder of $\cc$ of co-length $q_0$ containing $x$.
Then for the cylinder
$$\tcc' = C_{V}[i_{m-d_0}, i_{m-d_0+1}, \ldots, i_{m}, i_{m+1}, \ldots, i_m, i_{m+q_0}] \subset V $$
we have $\tpp^{m-d_0}(x) \in \tcc'$ and $\diam(\tcc') < \frac{1}{c_0 \gamma^{d_0+q_0}} < \delta/B$.
Since for any $y \in \cc'$ we have $\tpp^{m-d_0} (y) \in \tcc'$, it follows that 
$d(\tpp^{m-d_0}(x), \tpp^{m-d_0}(y)) < \delta/B$
and therefore $d(\phi_t(x), \phi_t(y)) < \delta$. Thus, $y$ satisfies (9.17). This proves the assertion with
$\rho_1 =  \frac{c_0\hr'_0}{L_1 B r'_0 \gamma_1^{d_0} \Gamma_0^3} $.
\endofproof


\bs


\bs

{\sc University of Western Australia, Crawley WA 6009, Australia}

{\sc\it E-mail address:} luchezar.stoyanov@uwa.edu.au

\end{document}